\newif\ifdraftmode
\draftmodefalse
\documentclass[10pt]{article}

\ifdraftmode
\fi

\newcommand\thisShorttitle{Optimal sampling for least squares approximation with general dictionaries}
\newcommand\thisTitle{\thisShorttitle}
\newcommand\thisSubject{math.NA, cs.LG}
\newcommand\thisAuthor{Philipp Trunschke, Anthony Nouy}
\newcommand\thisKeywords{least squares, optimal sampling, leverage score sampling, Christoffel function}

\usepackage{arxiv}

\usepackage[utf8]{inputenc}      % allow utf-8 input
\usepackage[T1]{fontenc}         % use 8-bit T1 fonts
\usepackage[dvipsnames]{xcolor}  % \color command
\usepackage{hyperref}            % hyperlinks
\usepackage{url}                 % simple URL typesetting
\usepackage{doi}                 % hyperlink doi numbers
\usepackage{microtype}           % microtypography
\usepackage{graphicx}
\hypersetup{
    breaklinks,
    raiselinks=true,
    colorlinks,
    citecolor=gray,
    filecolor=gray,
    linkcolor=gray,
    urlcolor=gray,
    bookmarksopenlevel=1,    % Gliederungstiefe der Bookmarks
    bookmarksopen=true,      % Expandierte Untermenues in Bookmarks
    bookmarksnumbered=true,  % Nummerierung der Bookmarks
    bookmarkstype=toc,       % Art der Verzeichnisses
    pdfdisplaydoctitle=true, % Display document title instead of file name in title bar.
    pdfstartview=FitH,       % Starting view of PDF document: Fits the width of the page to the window.
    pdfpagemode=UseOutlines, % Bookmarks im Viewer anzeigen
    pdfpagelayout=OneColumn, % Displays the document in one column; continuous scrolling.
    pdftitle={\thisTitle},
    pdfsubject={\thisSubject},
    pdfauthor={\thisAuthor},
    pdfkeywords={\thisKeywords},
}

\usepackage{lipsum}              % Can be removed after putting your text content

\usepackage{enumitem}            % provides user control over the layout of list environments
\usepackage[ruled,linesnumbered]{algorithm2e}  % algorithm environment
\usepackage{caption}             % Allow numbered figures without caption.
\usepackage{subcaption}          % Subfigues
\usepackage{float}               % H option for figure environment

\usepackage{csquotes}  % Needed for biblatex+babel
\usepackage[
    bibstyle=alphabetic,
    citestyle=alphabetic,
    sorting=nyt, % ynt,
    sortcites=true,
    giveninits=true,
    maxbibnames=99,
    backend=biber,
    maxalphanames=5,
    backref=true,
    natbib=true,     % \citet, \citep
    url=false,
    date=year, % terse,
    abbreviate=false,
]{biblatex}

\DeclareFieldFormat{eprint:hal}{%
    \mkbibacro{HAL}\addcolon\space
    \ifhyperref
        {\href{https://hal.archives-ouvertes.fr/#1}{%
        \nolinkurl{#1}
        \iffieldundef{eprintclass}
            {}
            {\addspace\UrlFont{\mkbibbrackets{\thefield{eprintclass}}}}}}
        {%
        \nolinkurl{#1}
        \iffieldundef{eprintclass}
            {}
            {\addspace\UrlFont{\mkbibbrackets{\thefield{eprintclass}}}}}}
\DeclareFieldAlias{eprint:HAL}{eprint:hal}

\DeclareSourcemap{
  \maps[datatype=bibtex]{
    \map{
      \step[fieldset=issn, null]
    }
  }
}
\renewbibmacro{in:}{}  % Do not print "In."

\def\keywordname{{\bfseries \emph{Key words.}}}%

\def\mscclassname{{\bfseries \emph{AMS subject classifications.}}}%
\def\mscclasses#1{\par\addvspace\medskipamount{\rightskip=0pt plus1cm
\def\and{\ifhmode\unskip\nobreak\fi\ $\cdot$
}\noindent\mscclassname\enspace\ignorespaces#1\par}}
\def\codename{{\bfseries \emph{Code.}}}%
\def\code#1{\par\addvspace\medskipamount{\rightskip=0pt plus1cm\noindent\codename\enspace\ignorespaces\url{#1}\par}}

\ifdraftmode
    \definecolor{amaranth}{rgb}{0.9, 0.17, 0.31}%
    \definecolor{americanrose}{rgb}{1.0, 0.01, 0.24}
    \usepackage[notref,notcite]{showkeys}
    
    \definecolor{bleudefrance}{rgb}{0.19, 0.55, 0.91}
    \newcommand{\numRevisions}{1}
    \newcommand{\revision}[2][0]{%
    \begingroup%
        \newcount\colorRatio%
        \colorRatio=\numexpr(100*(#1+1))/\numRevisions\relax%
        \colorlet{revisionColor}{RoyalBlue!\the\colorRatio!black}\color{revisionColor}#2%
    \endgroup}
    \newenvironment{revisione}[1][0]{\par%
        \begingroup%
        \newcount\colorRatio% 
        \colorRatio=\numexpr(100* (#1+1))/\numRevisions\relax%
        \colorlet{revisionColor}{RoyalBlue!\the\colorRatio!black}\color{revisionColor}
    }{\endgroup\par}
\else
    \newcommand{\note}[1]{}
    \newcommand{\revision}[2][0]{%
    \begingroup%
        #2%
    \endgroup}
    \newenvironment{revisione}[1][0]{\par\begingroup}{\endgroup\par}
\fi

\usepackage{pgf,tikz}
\usetikzlibrary{positioning}
\usetikzlibrary{calc}
\usetikzlibrary{intersections}
\usetikzlibrary{decorations.pathreplacing}  % braces

\tikzset{core/.style={inner sep=0pt}}
\tikzset{contraction/.style={line width=0.75}}
\tikzset{contractionDots/.style={contraction, dotted}}

\usepackage{amsmath}             % basic math symbols, operators and equation environments
\usepackage{amsfonts}            % blackboard math symbols
\usepackage{amssymb}             % \lesssim
\usepackage{amsthm}              % an enhanced version of LATEX’s \newtheorem command
\usepackage{thmtools}      % frontend for amsthm, Erst nach hyperref laden...
\usepackage[fixamsmath,disallowspaces]{mathtools} % extends amsmath and fixes some bugs, defines DeclareMathOperator and DeclarePairedDelimiters
\mathtoolsset{showonlyrefs=true}
\usepackage{newtxtext}
\usepackage[varvw]{newtxmath}  % mathptmx is very old and does not support bold fonts.

\colorlet{dimgray}{black!35!white}
\colorlet{lightgray}{dimgray!35!white}
\declaretheoremstyle[bodyfont=\itshape, mdframed={backgroundcolor=lightgray, linecolor=dimgray, linewidth=0.75pt, innertopmargin=1.5ex}]{claim}
\declaretheorem[style=claim]{theorem}
\declaretheorem[style=claim, numberlike=theorem]{lemma}
\declaretheorem[style=claim, numberlike=theorem]{proposition}

\declaretheoremstyle[mdframed={backgroundcolor=lightgray, linecolor=dimgray, linewidth=0.75pt, innertopmargin=1.5ex}]{definition}

\declaretheoremstyle[bodyfont=\itshape, mdframed={backgroundcolor=white, linecolor=dimgray, linewidth=0.75pt, innertopmargin=1.5ex}]{remark}
\declaretheorem[style=remark, numberlike=theorem]{remark}
\declaretheoremstyle[mdframed={backgroundcolor=white, linecolor=dimgray, linewidth=0.75pt, innertopmargin=1.5ex}]{example}

\newcommand{\indep}{\perp\kern-0.6em\perp}

\newcommand*{\mbb}[1]{\mathbb{#1}}
\newcommand*{\mcal}[1]{\mathcal{#1}}
\newcommand*{\mfrak}[1]{\mathfrak{#1}}
\newcommand*{\dd}{\ensuremath{\mathrm{d}}}
\newcommand*{\dx}[1][x]{\ensuremath{\,\dd{#1}}}

\let\inf\relax  % remove the definition of \inf before redeclaring it
\DeclareMathOperator*{\inf}{inf\vphantom{\sup}}
\DeclareMathOperator*{\argmin}{arg\,min}

\DeclarePairedDelimiter{\pars}{\ensuremath{(}}{\ensuremath{)}}
\DeclarePairedDelimiter{\bracs}{\ensuremath{[}}{\ensuremath{]}}
\DeclarePairedDelimiter{\braces}{\ensuremath{\{}}{\ensuremath{\}}}

\DeclarePairedDelimiter{\inner}{\langle}{\rangle}
\DeclarePairedDelimiter{\norm}{\|}{\|}
\DeclarePairedDelimiter{\abs}{\lvert}{\rvert}

\makeatletter
\newcommand{\opnorm}{\@ifstar\@opnorms\@opnorm}
\newcommand{\@opnorms}[1]{%
  \left|\mkern-1.5mu\left|\mkern-1.5mu\left|
   #1
  \right|\mkern-1.5mu\right|\mkern-1.5mu\right|
}
\newcommand{\@opnorm}[2][]{%
  \mathopen{#1|\mkern-1.5mu#1|\mkern-1.5mu#1|}
  #2
  \mathclose{#1|\mkern-1.5mu#1|\mkern-1.5mu#1|}
}
\makeatother

\mathchardef\mhyphen="2D

\usepackage{dsfont}              % \mathds{1}

\let\oldbullet\bullet
\newlength{\raisebulletlen}
\setbox1=\hbox{$\bullet$}\setbox2=\hbox{\tiny$\bullet$}
\renewcommand\bullet{\raisebox{\raisebulletlen}{\,\tiny$\oldbullet$}\,}

\makeatletter
\newcommand*{\rom}[1]{\expandafter\@slowromancap\romannumeral #1@}
\makeatother

\DeclarePairedDelimiterX\Set[1]\{\}{%
  #1%
}

\def\multiset#1#2{\ensuremath{\left(\kern-.3em\left(\genfrac{}{}{0pt}{}{#1}{#2}\right)\kern-.3em\right)}}

\newcommand{\sfrac}[2]{#1\hspace{-0.3ex}/\hspace{-0.15ex}#2}

\newcommand*{\xhat}[2][0.3em]{#2\kern-#1\hat{\vphantom{#2}}\kern#1}
\title{\thisTitle} % \thanks{A preliminary version of this manuscript is available on \texttt{arXiv}. The manuscript is neither published nor submitted elsewhere.}}
\renewcommand{\undertitle}{}
\renewcommand{\shorttitle}{\thisShorttitle}
\date{}

\author{
\href{https://orcid.org/0000-0002-2995-126X}{\includegraphics[height=0.7em]{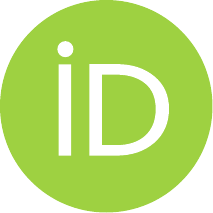}\hspace{1mm}\textcolor{black}{Philipp Trunschke\thanks{Website: \url{https://ptrunschke.github.io/}}}} \\
    Physikalisch-Technische Bundesanstalt \\
    Department 8.4 --- Mathematical Modelling and Data Analysis \\
    Abbestra\ss{}e 2--12\\
    10587 Berlin \\
    Germany \\
    \href{mailto:philipp.trunschke@ptb.de}{\texttt{philipp.trunschke@ptb.de}} \\
\And
\href{https://orcid.org/0000-0002-2149-2986}{\includegraphics[height=0.7em]{orcid.pdf}\hspace{1mm}\textcolor{black}{Anthony Nouy\thanks{Corresponding Author. Website: \url{https://anthony-nouy.github.io/}}}} \\
    Centrale Nantes \\
    Nantes Universit\'e \\
    Laboratoire de Math\'ematiques Jean Leray \\
    CNRS UMR 6629 \\
    France \\
    \href{mailto:anthony.nouy@ec-nantes.fr}{\texttt{anthony.nouy@ec-nantes.fr}}
}

\begin{document}
\maketitle
% \vspace{-5em}

\begin{abstract}
We consider the problem of approximating an unknown function from point evaluations.
This problem is a crucial subproblem in many modern (nonlinear) approximation schemes.
When obtaining these point evaluations is costly, minimising the required sample size becomes crucial.
Recently, an increasing focus has been on employing importance sampling strategies to achieve this.
For the approximation in a $d$-dimensional linear space, an optimal i.i.d.\ sampling measure 
% is obtained by minimising a bound for the number of i.i.d.\ samples required to achieve a stable approximation, achieving 
achieves
a sampling complexity of $\mathcal{O}(d\log (d))$.
%We show that this complexity is optimal for i.i.d.\ sampling.
However, the corresponding sampling measure depends on an orthonormal basis of the linear space, which is rarely known (in particular in the context of nonlinear approximation where the linear space arises as a local linearisation of a nonlinear model class like neural networks or tensor networks).
Consequently, sampling from these measures is challenging in practice.
This manuscript presents a strategy for estimating an orthonormal basis.
This strategy can be performed offline and does not require evaluations of the sought function.
We establish convergence and illustrate the practical performance through numerical experiments.
Comparing the presented approach with standard Monte Carlo sampling demonstrates a significant reduction in the number of samples required to achieve a good estimation of an orthonormal basis.
\end{abstract}

%NOTE: Also change the pdfkeywords in \hypersetup!
\vspace{-1em}
\mscclasses{15A03 \and 41A30 \and 62J02 \and 65Y20 \and 68Q25}
% \vspace{-1em}
% \code{https://github.com/ptrunschke/tangent_space_least_squares}

% {\color{red}
% \begin{itemize}
%     \item Theorems should be understandable with minimal context.
%     \item In the plots: Add an axvline where $k=d$.
% \end{itemize}
% }

\section{Introduction}
\label{sec:introduction}

Despite being one of the oldest statistical methods~\cite{legendre1806nouvelles}, least squares estimation remains the default method for regression in practice and an integral part of modern approximation algorithms such as the \emph{iteratively reweighted least squares algorithm} for sparse approximation~\cite{Chartrand2008,Daubechies2009}, the \emph{alternating least squares method} for low-rank approximation~\cite{oseledets_2011_dmrg,holtz_2012_alternating} as well as several methods for neural network approximation~\cite{ainsworth_2021_galerkin_nns,fokina_2023_growing_axons,ni_2023_elm,trunschke_2024_sgd,adcock2022cas4dl}.

In the small data regime or for unfavourably distributed data, it is known that least squares methods can become unstable.
This is classically tackled by regularisation~\cite{oster2021approximating,Trunschke2024,engl1996regularization,tibshirani_lasso_1996}.
Regularisation, however, requires knowledge of the sought function's regularity and a diligent choice of the regularisation parameter.
An alternative approach, inspired by importance sampling for Monte Carlo integration, is provided by optimal sampling~\cite{cohen_2017_optimal}.
Since its recent rediscovery, this approach has been applied in more and more applications, reaching from (multi-)linear and sparse approximation~\cite{haberstich_2020_thesis,haberstich_2022_boosted,cohen_2021_ppde,nobile_2020_optimal_ml,arras_2019_sequential} to more general nonlinear schemes for the optimisation of neural networks~\cite{trunschke_2024_sgd,Bruna_2024,adcock2022cas4dl}.

The optimal sampling method relies on linear spaces and requires knowledge of $L^2$-orthonormal bases of these spaces for practical application.
In practice, however, such bases may not be known, and we may only have access to generating systems.
For nonlinear model classes with differentiable parametrisations (e.g., tensor networks or neural networks with differentiable activation functions), for example, the linear space typically arises from a local linearisation of the model class within an iterative algorithm. 
In this context, it is natural to ask whether optimal sampling is still possible given only knowledge of a general  dictionary of functions that span a linear space.

\subsection{Contributions and Structure}

% {\color{red}
% \begin{itemize}
%     \item Explain the rough idea in this section.(Without introducing notation)
%     \item Describe the goal.
%     \item Describe the steps to arrive at the stated goal.
%     \item Repeat the main results of this paper!
% \end{itemize}
% }

This manuscript considers the problem of $L^2$-orthonormalising potentially overcomplete dictionaries of arbitrary functions with respect to general measures.
%An orthonormalisation algorithm is proposed and analysed, and its effectiveness is demonstrated in numerical experiments.
This orthonormalisation requires the Gramian matrix of the dictionary.
%Numerically stable weighted least squares approximation generally requires a careful selection of sample points.
%Drawing those points at random necessitates knowledge of an $L^2$-orthonormal basis to compute the \emph{optimal sampling measure}~\cite{cohen_2017_optimal}, 
% which requires knowledge of the Gramian matrix of the dictionary.
This manuscript presents a sampling strategy that iteratively refines a sample-based estimate of the Gramian matrix, allowing to obtain a measure that is asymptotically as good as the optimal measure.
This algorithm converges and reduces the orthogonalisation costs in practice.
The motivation for this is that orthogonalisation can not always be considered an offline step, like in the very important cases of iterative algorithms for nonlinear approximation discussed above.

Section~\ref{sec:setting} recalls known results on least squares approximation and (approximate) optimal sampling.
We complement this discussion by Theorem~\ref{thm:dlogd} which shows that the presented optimal sampling measure is indeed \emph{optimal}.
We only focus on approximating the optimal sampling measure and only briefly discuss the problem of sampling.
We prove in Theorem~\ref{thm:weighted_sample_size} that, while subsampling from a larger random sample set seems to work exceptionally well in practice~\cite{adcock2022cas4dl}, it has theoretical limitations.
This shows that sampling remains an important open problem.

Section~\ref{sec:iterative_refinement} introduces our algorithm for estimating the Gramian of a dictionary with respect to a general measure and an initialisation strategy that is based on determinantal point processes (Proposition~\ref{prop:vol_sampling}).
We also show that the initialisation requires at most $d$ steps (Proposition~\ref{prop:increasing_rank}).
%our proposed sampling and Gramian estimation methods.
In Section~\ref{sec:convergence_of_Gk} we analyse the proposed refinement algorithm, showing that the iterates are unbiased estimates with decreasing variance (Theorem~\ref{thm:mean_and_var}) and that the condition number of the approximation problem decreases with high probability (Theorem~\ref{thm:high_probability}).
Section~\ref{sec:experiments_intro} provides numerical experiments and demonstrates an application of the proposed algorithm to Christoffel--Darboux approximation.
The paper concludes with a discussion of the results in Section~\ref{sec:discussion}.

\subsection{Setting and related work}
\label{sec:setting}

\paragraph{Least squares approximation.} Let $\rho$ be a  measure on a set $\mcal{X}$ and consider a dictionary (or frame) of functions $\braces{B_1, \ldots, B_D} \subseteq L^2(\rho)$ represented as the components of the feature map $B : \mcal{X}\to\mbb{R}^D$ for some $D\in\mbb{N}$.
This dictionary spans a  $d$-dimensional subspace ($d\le D$)
$$
    \mcal{V} := \operatorname{span}\pars{B_1, \ldots, B_D} \subseteq L^2(\rho)
$$
and we consider the problem of finding the $L^2(\rho)$-best approximation of $u\in L^2(\rho)$ in $\mcal{V}$, formally,
\begin{equation}
\label{eq:best_approximation}
    u_{\mcal{V}} = \argmin_{v\in \mcal{V}}\ \norm{u - v}_{L^2(\rho)} .
\end{equation}
Commonly, the $L^2(\rho)$-norm cannot be computed exactly and is approximated by the sample-based estimate
\begin{equation}
\label{eq:empirical_norm}
    %\norm{v} \approx
    \norm{v}_{n} := \pars*{\frac{1}{n}\sum_{i=1}^n w\pars{x_i} \abs{v\pars{x_i}}^2}^{1/2} .
\end{equation}
The \emph{weight function} $w$ and the points $x_1,\ldots,x_n$ can be chosen in multiple ways.
When $w$ satisfies $w \ge 0$ and $\int w^{-1} \dx[\rho] = 1$, the measure $w^{-1}\rho$ is a probability measure and $x_1,\ldots,x_n\sim w^{-1}\rho$ can be drawn i.i.d.\ at random.
In this case, the semi-norm~\eqref{eq:empirical_norm} is an unbiased Monte Carlo estimate of the $L^2(\rho)$ norm.
Replacing the true $L^2(\rho)$ norm in~\eqref{eq:best_approximation} by the estimate~\eqref{eq:empirical_norm} yields the estimate
\begin{equation}
\label{eq:empirical_best_approximation}
    u_{\mcal{V},n} = \argmin_{v\in\mcal{V}}\ \norm{u - v}_{n} .
\end{equation}

\begin{minipage}{\textwidth}
\begin{remark}
    Solving~\eqref{eq:best_approximation} with general dictionaries can be numerically challenging due to ill-conditioning.
    % Gramian matrix $G\in\mathbb{R}^{D\times D}$, given by $G_{kl} = (B_k, B_l)_{L^2(\rho)}$.
    % This has been investigated in~\cite{Adcock2019_frames,Adcock2020_frames}.
    In the context of frames,~\cite{Adcock2019_frames} replaces the $L^2(\rho)$-orthogonal projection $u_{\mathcal{V}}$ by a regularised version $u_{\mathcal{V}}^{\varepsilon}$ satisfying
    $$
        \norm{u - u_{\mathcal{V}}^\varepsilon}
        \le \norm{u - u_{\mathcal{V}}} + \sqrt{\varepsilon}\|u_{\mathcal{V}}\| \ .
    $$
    \Cite{Adcock2020_frames}~presents a refined algorithm for which the $\mathcal{O}(\sqrt{\varepsilon})$-term can be replaced by an $\mathcal{O}(\varepsilon)$-term under mild conditions.
\end{remark}
\end{minipage}

\paragraph{Error bounds.} The quality of the approximation~\eqref{eq:empirical_norm} can be measured in terms of the inequality
\begin{equation}
\label{eq:rip}
    \pars{1-\delta}\norm{v}_{L^2(\rho)}^2
    \le \norm{v}_{n}^2,
    % \le \pars{1+\delta}\norm{v}^2,
    \qquad
   \forall v\in \mcal{V},
\end{equation}
and critically influences the quality of the resulting estimate~\eqref{eq:empirical_best_approximation}.
If the condition~\eqref{eq:rip} is satisfied for some $0 < \delta < 1$, it can be shown~\cite{cohen_2017_optimal} that
\begin{equation}
\label{eq:error_bound_n}
    \norm{u - u_{\mcal{V},n}}_{L^2(\rho)}^2
    \le \norm{u - u_{\mcal{V}}}_{L^2(\rho)}^2 + \tfrac{1}{1-\delta} \norm{u - u_{\mcal{V}}}_{n}^2 ,
\end{equation}
i.e.\ that the error of the empirical best approximation~\eqref{eq:empirical_best_approximation} is bounded by the error of the $L^2(\rho)$-best approximation~\eqref{eq:best_approximation}.
Since $\norm{\bullet}_n$ and $\norm{\bullet}_{L^2(\rho)}$ are dominated by $\norm{w^{1/2}\bullet}_{L^\infty(w^{-1}\rho)}$, equation~\eqref{eq:error_bound_n} yields the immediate estimate
\begin{align}
\label{eq:error_bound_Linf}
    \norm{u - u_{\mcal{V},n}}_{L^2(\rho)}^2
    &\le \norm{u - u_{\mcal{V}}}_{L^2(\rho)}^2 + \tfrac{1}{1-\delta} \norm{w^{1/2}(u - u_{\mcal{V}})}_{L^\infty(w^{-1}\rho)}^2 \\
    &\le \big(1 + \tfrac{1}{1-\delta}\big) \norm{w^{1/2}(u - u_{\mcal{V}})}_{L^\infty(w^{-1}\rho)}^2
    .
\end{align}
However, due to the different norms on both sides of the inequality, this bound can not be truly called quasi-optimal.
A bound in expectation avoiding the stronger weighted $L^\infty$-norm is derived in~\cite{haberstich_2022_boosted}.
%A slight adaptation of this bound is presented in the subsequent theorem and proven in Appendix~\ref{app:error_bound_expectation}.

\begin{theorem}
[\cite{haberstich_2022_boosted} Theorem 2.4]
\label{thm:error_bound_expectation}
    Suppose the probability of the random event~\eqref{eq:rip} is lower bounded by $p$.
    Then
    \begin{equation}
    \label{eq:error_bound_L2}
        \mbb{E}\bracs*{\norm{u - u_{\mcal{V},n}}_{L^2(\rho)}^2\;\big|\;\eqref{eq:rip}}
        \le \pars*{1 + \tfrac{1}{(1-\delta) p}} \norm{u - u_{\mcal{V}}}_{L^2(\rho)}^2 .
    \end{equation}
\end{theorem}
Given these results, it is not surprising that the probability of~\eqref{eq:rip} has been the focus of numerous works~\cite{cohen_2017_optimal,haberstich_2022_boosted,chkifa2024randomized,nouy2024weighted,eigel_2022_ls}.

\paragraph{Sample size bounds.} To bound the probability of the random event~\eqref{eq:rip}, it is useful to define the Gramian matrix $G\in\mbb{R}^{D\times D}$ and its associated inverse Christoffel function $\mfrak{K}_G : \mcal{X}\to[0,\infty)$ via
\begin{equation}
\label{eq:gramian}
    G_{kl} := (B_k, B_l)_{L^2(\rho)}
    \qquad\text{and}\qquad
    \mfrak{K}_G\pars{x}
    % := \sup_{v\in\mbb{R}^D} \frac{\abs{B\pars{x}^\intercal v}^2}{v^\intercal G v},
    := B\pars{x}^\intercal G^+ B\pars{x}
    = \inner{G^+, B(x)B(x)^\intercal}_{\mathrm{Fro}},
\end{equation}
where $\inner{\bullet,\bullet}_{\mathrm{Fro}}$ denotes the Frobenius-inner product and $G^+$ denotes the Moore--Penrose pseudo-inverse of $G$.
Proposition~\ref{prop:inverse_christoffel_definition} in Appendix~\ref{app:christoffel} demonstrates that this definition is equivalent to the customary definition
$
    \mfrak{K}_G\pars{x} := \sup_{v\in\mcal{V}\setminus\{0\}} \frac{\abs{v(x)}^2}{\norm{v}_{L^2(\rho)}^2} .
$
This equivalence demonstrates that $\mfrak{K}_G$ depends only on the space $\mcal{V}$ even though it is defined in terms of the dictionary $B$.

% where we assume that the supremum is taken only over $v\in\ker(G)^\perp$ for which the quotient is defined.
The probability of~\eqref{eq:rip} holding is then bounded by the following result, which is a slight adaptation of~\cite[Lemma~2.1]{dolbeault_2021_domains} or~\cite[Theorem~2.5]{haberstich_2022_boosted} and is a direct consequence of Theorem~1.1 in~\cite{Tropp2011}.

\begin{theorem}
\label{thm:sample_size_bounds}
    Let $p, \delta\in(0,1)$, $\gamma := ((1+\delta)\ln(1+\delta)-\delta)^{-1}$ and assume that
    \begin{equation}
    \label{eq:sample_size_bound}
        n \ge \norm{w\mfrak{K}_G}_{L^\infty(w^{-1}\rho)}\gamma \ln(\tfrac{d}{p})
    \end{equation}
    independent points $x_1,\ldots, x_n$ are drawn according to $w^{-1}\rho$.
    Then,
    $$
        \mathbb{P}\bracs*{
                \pars{1-\delta}\norm{v}_{L^2(\rho)}^2
                \le \norm{v}_{n}^2\;,\; \forall v\in\mcal{V}
                %\le \pars{1+\delta}\norm{v}_{L^2(\rho)}^2
        }
        \ge 1 - p .
    $$
\end{theorem}
The idea of \emph{optimal sampling} is choosing the weight function
\begin{equation}
\label{eq:optimal_weight_function}
    w_G = \norm{\mfrak{K}_G}_{L^1(\rho)}\mfrak{K}_G^{-1} ,
\end{equation}
which minimises the lower bound~\eqref{eq:sample_size_bound}
% in equation~\eqref{eq:sample_size_bound}
because
\begin{equation}
     \norm{w\mfrak{K}_G}_{L^\infty(w^{-1}\rho)}
    \ge \norm{w\mfrak{K}_G}_{L^1(w^{-1}\rho)}
    = \norm{\mfrak{K}_G}_{L^1(\rho)}
    = \norm{w_G\mfrak{K}_G}_{L^\infty(w_G^{-1}\rho)} .
\end{equation}
This choice yields the factor
$$
    \norm{w_G \mfrak{K}_G}_{L^\infty(w_G^{-1}\rho)}
    = \norm{\mfrak{K}_G}_{L^1(\rho)}
    = \int \inner*{G^+, B(x)B(x)^\intercal} \dx[\rho]
    = \inner*{G^+, \int B(x)B(x)^\intercal \dx[\rho]}
    = \inner*{G^+, G}
    = d .
$$
Notably, the theory has already been known for discrete domains $\mathcal{X}$, where the values of the optimal weight function $w$ are known as \emph{leverage scores}~\cite{Mahoney2010}.

\revision[0]{
Although there exist spaces $\mathcal{V}$ for which $n=d$ random points are sufficient to guarantee almost surely $\Vert v \Vert_n >0$ for all nonzero $v$ in $\mcal{V}$ (e.g.\ polynomials or trigonometric polynomials), the resulting $\delta$ may be extremely close to $1$,\footnote{\revision[0]{If $\Vert v \Vert_n >0$ holds for all nonzero $v \in \mcal{V}$, then~\eqref{eq:rip} holds for some (random) $\delta$, depending on the sample $x_1, \ldots, x_n$.}} resulting in a huge condition number of the corresponding least squares problem~\eqref{eq:empirical_best_approximation}.  
%Also, satisfying \eqref{eq:rip} with a fixed $\delta>0$ is required to obtain quasi-optimality results in expectation.
Moreover, for general spaces $\mathcal{V}$ the bound~\eqref{eq:sample_size_bound} is optimal, since the probability of~\eqref{eq:rip} can drop sharply when $x_1,\ldots,x_n\sim w^{-1}\rho$ are independent and $n < d\ln(d)$. 
We show this optimality in Theorem \ref{thm:dlogd} in Appendix  ~\ref{app:dlogd}. 
%We show this in the subsequent theorem, which is a generalisation of Proposition~5 in~\cite{Trunschke2024}.
% (a generalisation of Proposition~5 in~\cite{Trunschke2024}). 
% This theorem also shows that the factor $\|\mathfrak{K}_G\|_{L^\infty(\rho)}$ appears not only in the sufficient sample size bound~\eqref{eq:sample_size_bound} but also in the necessary sample size bound~\eqref{eq:sample_size_bound_necessary}.
%A proof is given in Appendix~\ref{app:dlogd}.
}

Recent work has strived to reduce the complexity further to $n = \mathcal{O}(d)$ by subsampling~\cite{haberstich_2022_boosted,Dolbeault2022,pmlr-v99-chen19a}\footnote{The references~\cite{pmlr-v99-chen19a} and~\cite{Dolbeault2022} provide theoretical guarantees for discrete and continuous spaces.
The reference~\cite{haberstich_2022_boosted} provides an easy-to-implement heuristic method, which, however, does not guarantee a sample size $n = \mathcal{O}(d)$.} or dependent sampling~\cite{chkifa2024randomized,nouy2024weighted}.
Both approaches are related to optimal sampling, extending the relevance of the present contribution to these approaches.

% \subsection{Related work}
% \label{sec:related_work}

% {\color{red}TODO: Various works~\cite{dolbeault_2021_domains,adcock_2020_domain,Migliorati2020,adcock_2022_domain} address this problem in the special case of polynomial approximation and uniform measures on sufficiently regular compact domains, but they require the knowledge of specific polynomial bases.}

\paragraph{Approximate optimal sampling.} The optimal weight function for i.i.d.\ sampling in~\eqref{eq:optimal_weight_function} depends on the unknown Gramian matrix $G$.
To bypass the problem of finding the exact Gramian matrix $G$, it seems natural to replace $G$ with an estimate $\hat{G}^{(0)}$ and define
\begin{equation}
\label{eq:w_opt_est}
    w_{\hat{G}^{(0)}} := \norm{\mfrak{K}_{\hat{G}^{(0)}}}_{L^1(\rho)} \mfrak{K}_{\hat{G}^{(0)}}^{-1} .
\end{equation}
If $\mfrak{K}_{\hat{G}^{(0)}}$ satisfies the uniform framing
\begin{equation}
\label{eq:framing}
    c \mfrak{K}_G \le \mfrak{K}_{\hat{G}^{(0)}} \le C\mfrak{K}_G
\end{equation}
for some constants $0 < c < C$, then~\cite{dolbeault_2021_domains}
\begin{equation}
\label{eq:suboptimality_factor}
    % \norm{w_{\hat{G}^{(0)}} \mfrak{K}_G}_{L^\infty(w_G^{-1}\rho)} 
    % = \norm*{\frac{\mfrak{K}_G}{\mfrak{K}_{\hat{G}^{(0)}}}}_{L^\infty(w_G^{-1}\rho)}
    %   \norm{\mfrak{K}_{\hat{G}^{(0)}}}_{L^1(\rho)}
    % \le \norm*{\frac{\mfrak{K}_G}{\mfrak{K}_{\hat{G}^{(0)}}}}_{L^\infty(w_G^{-1}\rho)}
    %     \norm*{\frac{\mfrak{K}_{\hat{G}^{(0)}}}{\mfrak{K}_G}}_{L^\infty(w_G^{-1}\rho)} 
    %     \norm{\mfrak{K}_G}_{L^1(\rho)}
    % \le \frac{C}{c}d .
    \norm{w_{\hat{G}^{(0)}} \mfrak{K}_G}_{L^\infty(w_G^{-1}\rho)} 
    = \norm*{\frac{\mfrak{K}_G}{\mfrak{K}_{\hat{G}^{(0)}}}}_{L^\infty(w_G^{-1}\rho)}
      \norm{\mfrak{K}_{\hat{G}^{(0)}}}_{L^1(\rho)}
    \le \frac{1}{c} \norm{\mfrak{K}_{\hat{G}^{(0)}}}_{L^1(\rho)}
    \le \frac{C}{c} \norm{\mfrak{K}_{G}}_{L^1(\rho)}
    = \frac{C}{c}d .
\end{equation}
% {\color{red}Note that we could also stop after the first equality, showing that we just need the upper bound $\mathfrak{K}_{\hat{G}^{(0)}} \ge C\mathfrak{K}_{G}$ and a bound on $\|\mathfrak{K}_{\hat{G}^{(0)}}\|$.}
This bound is almost optimal since it differs from the theoretically optimal $\norm{w_G\mfrak{K}_G}_{L^\infty(w_G^{-1}\rho)} = d$ just by the factor $\frac{C}{c}$.
The challenge is thus to find an estimate ${\hat{G}^{(0)}}$ that satisfies the uniform framing~\eqref{eq:framing}.
This is simplified by the subsequent generalisation of Lemma~3.2 from~\cite{dolbeault_2021_domains}.
\begin{lemma}
\label{lem:framing_equivalence}
    Let $G, H \in\mbb{R}^{D\times D}$ be \revision[0]{symmetric and positive semi-definite}.
    Then
    $$
        C^{-1} G \preceq H \preceq c^{-1} G
        \qquad\Leftrightarrow\qquad
        \ker(H) = \ker(G)
        \quad\text{and}\quad
        c \mfrak{K}_G \le \mfrak{K}_{H} \le C\mfrak{K}_G .
    $$
\end{lemma}
The proof of this lemma, together with a brief discussion of the condition $\ker(H) = \ker(G)$, can be found in Appendix~\ref{app:matrix_framing}.
We now briefly discuss three approaches to solve the problem of estimating $G$.

\paragraph{Monte Carlo approach.}
Provided that $\rho$ is a probability measure, the arguably simplest and most natural estimate of $G$ is given by
\begin{equation}
\label{eq:hatG_0}
    \hat G^{(0)} := \frac{1}{n} \sum_{i=1}^n B\pars{x^{(0)}_{i}} B\pars{x^{(0)}_{i}}^\intercal
\end{equation}
with $x^{(0)}_{1}, \ldots, x^{(0)}_{n}$ independently drawn from the original measure $\rho$.
This approach was originally proposed, independently, in~\cite{Migliorati2020} and~\cite{Adcock2020_near_optimal}, and Lemma~\ref{lem:framing_equivalence} ensures that a matrix framing for $\hat{G}^{(0)}$ implies the uniform framing for $\mfrak{K}_{{\hat{G}^{(0)}}}$ and hence, that a sufficiently good estimate $\hat{G}^{(0)}$ of $G$ will yield an adequate weight function $w_{\hat{G}^{(0)}}$.
The problem with this approach is that the sample size bound~\eqref{eq:sample_size_bound} from Theorem~\ref{thm:sample_size_bounds}
% that ensures a framing with $c \ge \frac{1}{1+\delta}$ and $C \le \frac{1}{1-\delta}$ with high probability
scales with $n \gtrsim \norm{\mfrak{K}_G}_{L^\infty(\rho)}$.
Even if $G$ can be estimated in an offline phase, the Monte Carlo approach may be infeasible since $\norm{\mfrak{K}_G}_{L^\infty(\rho)}$ can, in principle, grow arbitrarily fast with the dimension of $\mathcal{V}$ and can even be unbounded.\footnote{The constant $\norm{\mfrak{K}_G}_{L^\infty(\rho)}$ is unbounded when $B$ are polynomials and $\rho$ is the Gaussian measure.
An example where $\norm{\mfrak{K}_G}_{L^\infty(\rho)}$ grows exponentially is provided in section~\ref{ex:pw_const}.
In particular, even though $d=18$, the constant $\norm{\mfrak{K}_G}_{L^\infty(\rho)}$ is of the order of $10^{16}$.
The sample size bound~\eqref{eq:sample_size_bound} thus indicates that $n \gtrsim 10^{16}$ samples 
%(as many as are used for training contemporary foundation models) 
would be required to control the probability of~\eqref{eq:rip}.}

Nevertheless, it should be pointed out that the Monte Carlo approach is a generic approach which works in a general black box setting.
A striking example is optimising neural networks~\cite{trunschke_2024_sgd,adcock2022cas4dl} where $\mathcal{V}$ arises from a local linearisation of the model class and the exact computation of the corresponding Gramian $G$ is generally NP-hard~\cite{blaschko_2019_nnnorms}.

% \Cite{adcock_2020_domain} analyses the behaviour of Christoffel sampling when the Christoffel function is associated with a product basis on a tensor product domain containing the actual domain of interest.
\paragraph{Framing approach.}
If $\|\mathfrak{K}_G\|_{L^\infty(\rho)}$ is sufficiently small, we can compute the unweighted Monte Carlo estimate~\eqref{eq:hatG_0} and then use $\hat{G}^{(0)}$ to draw a smaller sample of size $\mathcal{O}(d \ln(d))$ as proposed in~\cite{dolbeault_2021_domains}.
However, bounding $\|\mathfrak{K}_G\|_{L^\infty(\rho)}$ a priori requires intimate knowledge of the basis $B$ and the measure $\rho$.
One approach to do this was originally proposed in~\cite{adcock_2020_domain} and relies on the subsequent fundamental lemma.
% {\color{blue}As with the framing approach, the problem with the multi-level framing approach from~\cite{dolbeault_2021_domains} is that finding a sequence of subspaces satisfying the framing~\eqref{eq:subspace_growth} requires deep knowledge of the basis $B$ and the measure $\rho$.}
% One approach, originally proposed in~\cite{dolbeault_2021_domains}, is to prove the uniform framing~\eqref{eq:framing} for the inverse Christoffel function $\mathfrak{K}_{\tilde{G}}$ for a Gramian $\tilde{G}$ with respect to a different measure $\tilde{\rho}$ and prove that $\tilde{G}$ satisfies the matrix framing $C^{-1}\tilde{G} \preceq G \preceq c^{-1}\tilde{G}$ from Lemma~\ref{lem:framing_equivalence}.
% % to draw the sample points $x^{(0)}_{1}, \ldots, x^{(0)}_{n}$ and estimate $\hat{G}^{(0)}$.
% \color{Purple}
% % and later refined in~\cite{adcock_2022_domain,Adcock2020_near_optimal}
% Both of these bounds rely on the subsequent fundamental lemma originally in~\cite{adcock_2020_domain}.
\begin{lemma}
\begin{revisione}[0]
    For every measurable subset $R\subseteq\mathcal{X}$,  let $G_R$ be the Gramian matrix of $B$ with respect to the measure $\rho_R(\bullet) := \rho(\bullet \cap R)$.
    Let $\mathcal{R}$ be a family of subsets of $\mathcal{X}$ satisfying $\ker(G_R)=\ker(G)$ for every $R\in\mathcal{R}$ and that for every $x$ there exists $R\in\mathcal{R}$ with $x \in R$.
    Then
    $$
        \mathfrak{K}_G(x)
        \le \inf_{\substack{R\in\mathcal{R}\\x\in R}} \mathfrak{K}_{G_{R}}(x).
    $$
\end{revisione}
\end{lemma}
% \begin{proof}
%     Let $R\in\mathcal{R}$ and observe that
%     $$
%         \|v\|_{L^2(\rho)}^2
%         = \int_{\mathcal{X}} |v(x)|^2 \dx[\rho](x)
%         \ge \int_{R} |v(x)|^2 \dx[\rho](x)
%         = \rho(R) \int_{R} |v(x)|^2 \dx[\tfrac{\rho}{\rho(R)}](x)
%         = \rho(R) \|v\|_{L^2(\rho(\bullet\,|\,R))}^2 .
%     $$
%     Inserting this bound into the equivalent definition of $\mathfrak{K}_G$ from Proposition~\ref{prop:inverse_christoffel_definition} yields
%     \begin{align}
%         \mathfrak{K}_{G}(x)
%         &= \sup_{v\in \mathcal{V}\setminus\{0\}} \frac{|v(x)|^2}{\|v\|_{L^2(\rho)}^2}
%         \le \sup_{v\in \mathcal{V}\setminus\{0\}} \frac{|v(x)|^2}{\rho(R) \|v\|_{L^2(\rho(\bullet\,|\,R))}^2}
%         = \tfrac{1}{\rho(R)} \mathfrak{K}_{G_{R}}(x) . \qedhere
%     \end{align}
% \end{proof}
Assuming $\rho$ is the uniform measure on $\mathcal{X}$ and $\mathcal{V}$ is a space of polynomials of bounded degree, the works~\cite{adcock_2020_domain, dolbeault_2021_domains} use this lemma to derive conditions
% ($\lambda$-rectangle property~\cite{adcock_2020_domain}, inner cone condition~\cite{dolbeault_2021_domains})
on $\mathcal{X}$ ensuring $\|\mathfrak{K}_G\|_{L^\infty(\rho)} \lesssim n$.\footnote{
\Cite{dolbeault_2021_domains} assumes $\mathcal{V}$ is a space of polynomials of bounded degree.
\Cite{adcock_2020_domain} assumes $\mathcal{V}$ is spanned by Legendre polynomials for a downward closed index set, like hyperbolic crosses.
Various other lower and upper bounds are derived in~\cite{dolbeault_2021_domains}.}

\paragraph{Multi-level framing approach.} Expanding on the framing approach, the authors of~\cite{dolbeault_2021_domains} assume that the space $\mcal{V}$ can be decomposed into a nested sequence of subspaces $\mcal{V}_1 \subset \ldots \subset \mcal{V}_d = \mcal{V}$ with corresponding Gramian matrices $G_{\mcal{V}_1},\ldots,G_{\mcal{V}_d}$ satisfying the growth conditions
\begin{equation}
\label{eq:subspace_growth}
    \mfrak{K}_{G_{\mcal{V}_1}} \le C
    \qquad\text{and}\qquad
    \mfrak{K}_{G_{\mcal{V}_{k+1}}} \le C \mfrak{K}_{G_{\mcal{V}_{k}}} .
\end{equation}
The fact that $\|\mfrak{K}_{G_{\mcal{V}_1}}\|_{L^\infty(\rho)} \le C$ allows us to estimate $G_{V_1}$ using the Monte Carlo approach~\eqref{eq:hatG_0}.
Moreover, since $V_k \subset V_{k+1}$, the growth condition~\eqref{eq:subspace_growth} ensures the uniform framing $\mfrak{K}_{G_{\mcal{V}_{k}}}\! \le \mfrak{K}_{G_{\mcal{V}_{k+1}}}\! \le C\mfrak{K}_{G_{\mcal{V}_{k}}}$.
Due to this framing, we can sequentially estimate the Gramian matrices ${G}_{\mcal{V}_{k+1}}$ by $\hat{G}_{\mcal{V}_{k+1}}$ using the estimated optimal sampling density $\mfrak{K}_{\hat G_{\mcal{V}_{k}}}$ and requiring only $C(k+1)$ sample points with high probability.
When $\rho$ is a uniform measure on a bounded domain and $\mcal{V}$ is a space of polynomials, bounds for the growth of the inverse Christoffel function have been studied in~\cite{dolbeault_2021_domains}.

\paragraph{Sampling.}
When a sufficiently good estimate $\hat{G}^{(0)}$ of $G$ is known, we can solve the least squares problem~\eqref{eq:empirical_best_approximation} using a sample $x_1,\ldots,x_n\sim w_{\hat{G}^{(0)}}^{-1}\rho$ with $w_{\hat{G}^{(0)}}$ defined as in~\eqref{eq:w_opt_est} and $n\in\mathbb{N}$ chosen according to~\eqref{eq:sample_size_bound}.
Although this procedure seems straightforward, sampling from the measure $w_{\hat{G}^{(0)}}^{-1}\rho$ may be non-trivial due to its high-dimensional and multimodal nature.
Even though sampling is not the focus of the present paper, we want to comment on several approaches.

For bounded low-dimensional domains $\mathcal{X}\subseteq\mathbb{R}^k$, sampling from $w_{\hat{G}^{(0)}}^{-1}\rho$ can be performed almost trivially by discretising the measure.
When $\mathcal{X}$ is a product domain and $B$ is a known orthonormal product basis of $\mathcal{V}$, the authors of~\cite{cohen_2017_optimal} propose a sequential sampling scheme that exploits the fact that conditional distributions can be easily computed.
This scheme can be generalised to any other case where the Lebesgue density of $w_{\hat{G}^{(0)}}^{-1}\rho$ can be well approximated by hierarchical tensor networks~\cite{Dolgov2019}.
In~\cite{adcock2022cas4dl} the authors propose to draw an i.i.d.\ $w_{\hat{G}^{(0)}}^{-1}\rho$-distributed subsample of a larger i.i.d.\ $\rho$-distributed sample.
Despite the astonishing practical performance, Theorem  \ref{thm:weighted_sample_size} in Appendix \ref{app:weighting} indicates that such a subsampled least squares approximation can only exhibit a quasi-optimal error bound~\eqref{eq:error_bound_Linf} when the larger $\rho$-distributed sample already ensures such a bound.\footnote{
This makes intuitive sense since we would expect the overall error of an approximation problem to increase when equations are removed from the (overdetermined) least squares system.
}
%\begin{minipage}{\textwidth}
%\begin{theorem}
%\label{thm:weighted_sample_size}
%\begin{revisione}[0]
%    Let $\mathcal{V}$ be chosen as in Theorem~\ref{thm:dlogd} and $\rho$ be a probability measure.
%    Suppose that $X_1,\ldots,X_m\sim\rho$ are independent and let $x_1,\ldots,x_n$ be a subsample of $X_1,\ldots,X_m$ (drawn with or without replacement according to any density).
%    Then $m\ge \tfrac12 \|\mathfrak{K}_G\|_{L^\infty(\rho)} \ln(\frac{1}p)$ initial points are necessary to ensure that
%    $$
%        \exists 
%        \,c\in(0,\infty) : \|u - u_{\mathcal{V},n}\|_{L^2(\rho)}^2
%        \le \|u - u_{\mathcal{V}}\|_{L^2(\rho)}^2 + c \|u - u_{\mathcal{V}}\|_{L^\infty(\rho)}^2
%    $$
%    %%%%%for any $c \ge 1$ 
% with probability $1-p$.
%    Moreover, $m\ge 3 \|\mathfrak{K}_G\|_{L^\infty(\rho)} \frac{c^2}{(c - C_w)^2} \ln(\frac{2d}p)$ initial points are sufficient to ensure that
%    $$
%        \|u - u_{\mathcal{V},n}\|_{L^2(\rho)}^2
%        \le \|u - u_{\mathcal{V}}\|_{L^2(\rho)}^2 + c \|u - u_{\mathcal{V}}\|_{L^\infty(\rho)}^2
%    $$
%    for any $c \ge C_w$ with probability $1-p$.
%    The precise value of $C_w \ge \frac{1}{n}$ depends on $w$,   $X_1, \ldots, X_m$, and and orthonormal basis $b$ of $\mathcal{V}$.
%    In particular, this bound requires an order of $m \in \Theta(\|\mathfrak{K}_G\|_{L^\infty(\rho)})$ initial points.
%\end{revisione}
%\end{theorem}
%\end{minipage}
% This corollary is proven in Appendix~\ref{app:weighting}, where some additional numerical illustration is provided too.
The proof of this theorem, together with a general discussion of the ineffectiveness of reweighting for least squares regression, can be found in Appendix~\ref{app:weighting}.

Theorem~\ref{thm:weighted_sample_size} shows that the sampling strategy proposed in~\cite{adcock2022cas4dl} requires the original sample to be of order $\Theta(\|\mathfrak{K}_G\|_{L^\infty(\rho)})$, which, as discussed before, may be infeasibly large.
This indicates that sampling from the inverse Christoffel function for general generating systems is an important problem and that more research is needed.

For the remainder of this paper we assume that there exists some method to draw the required sample points with a fixed cost.

% \section{\color{red}Idea/Proposal}
% \label{sec:contributions}
\section{Iterative refinement of the Gramian}
\label{sec:iterative_refinement}

The principal idea we explore in this manuscript is to use the sampling density $w_{\hat{G}^{(0)}}^{-1}\rho$ that is induced by a first estimate ${\hat{G}^{(0)}}  $ to refine the estimate ${\hat{G}^{(0)}}$ itself, pulling ourselves up by the bootstraps. 
To demonstrate why this is a sensible idea, we recall that every function $v\in\mcal{V}$ can be expanded as
$$
    v(x) = B(x)^\intercal c_v,
$$
for some $c_v\in\mbb{R}^D$.
This means that any symmetric, positive semi-definite matrix $H\in\mbb{R}^{D\times D}$ with $\operatorname{ker}(H)\supseteq\operatorname{ker}(G)$ induces a {semi}-inner product on $\mcal{V}$  given by
$$
    (v, w)_{H} := c_v^\intercal H c_w,
    \qquad
    v,w\in\mcal{V} .
$$
We denote the induced norm by $\norm{\bullet}_H$ and note that $\norm{v}_G = \norm{v}_{L^2(\rho)}$ for $v\in\mcal{V}$.
Now, let $\braces{b_1,\ldots, b_d}\subseteq\mcal{V}$ be the unique $G$-orthonormal and ${\hat{G}^{(0)}}$-orthogonal basis of $\mcal{V}$, i.e.
\begin{equation}
\label{eq:ortho-ortho-basis}
    (b_k, b_l)_G = \delta_{kl}
    \qquad\text{and}\qquad
    (b_k, b_l)_{\hat{G}^{(0)}} = \beta_k \delta_{kl} ,
\end{equation}
for some $\beta_1,\ldots,\beta_d \ge 0$.\footnote{The basis functions are given by $b_k(x) = B(x)^\intercal c_k$, where $(c_k, \beta_k)$ are the eigenpairs of the matrix pencil $(\hat{G}^{(0)}, G)$, i.e.\ $\hat{G}^{(0)} c_k = \beta_k G c_k$, with the normalisation condition $c_k^\intercal G c_k = 1$.}
Expressed in this basis, the partial norm equivalence~\eqref{eq:rip} can be equivalently written as
$$
    (1-\delta)\norm{b_k}_{G}^2 \le \norm{b_k}_{\hat{G}^{(0)}}^2,
    % (1-\delta) \le \lambda_k = \norm{b_k}_{\hat{G}^{(0)}}^2,
    % \le (1+\delta),
    % \abs{1 - \lambda_k} \le \delta,
    \qquad k=1,\ldots,d .
$$
The property~\eqref{eq:rip} hence measures how well the true $G$-norms of the basis functions are estimated by their corresponding ${\hat{G}^{(0)}}$-norm estimates.
By Proposition~\ref{prop:inverse_christoffel_definition}, the optimal sampling density $w_G^{-1} = \norm{\mfrak{K}_G}_{L^1(\rho)}^{-1}\mfrak{K}_G$ from~\eqref{eq:optimal_weight_function} is a mixture of densities
$$
    w_G(x)^{-1}
    = \frac{1}{d} \sum_{k=1}^d b_k^2(x) .
$$
Each summand in this mixture is the optimal
% \footnote{\comment{A: the following comment is not useful to rely on Bernstein, variance is zero.}\red{
%     The importance sampling density $\tilde{w}_k^{-1} = b_k^2$ minimises the uniform bound $\norm{\tilde{w}_k b_k^2}_{L^\infty(\rho)}$ as well as the variance of the estimate~\eqref{eq:MC_bk2_importance} (cf.~\cite{elvira2022advances}) and thereby maximises the rate of convergence
%     $
%         \norm{b_k}_{\hat{G}^{(0)}}^2 \xrightarrow{n\to\infty} \norm{b_k}_G^2
%     $
%     provided by Bernstein's inequality.}
% }
importance sampling density ${w}_k^{-1} := b_k^2$ for the Monte Carlo estimate for the squared norm $\norm{b_k}_G^2 = \int b_k^2 \dx[\rho] = 1$, which ensures that almost surely
\begin{equation}
\label{eq:MC_bk2_importance}
    %\norm{b_k}_G^2
    %= \int b_k^2 \dx[\rho]
    %\approx 
    \frac1{\tilde{n}} \sum_{i=1}^{\tilde{n}} {w}_k(x_{i}) b_k(x_{i})^2
    = 1,% \norm{b_k}_{\hat{G}^{(0)}}^2,
    \qquad \text{for }
    x_{1}, \ldots, x_{\tilde{n}} \sim {w}_k^{-1}\rho
    .
\end{equation}
% {\color{red}(
% Note here that Bernstein depends on an $L^\infty(\rho)$ factor
% The $L^\infty(\rho)$-factortor is optimised by this choice of $w$ for the same reason as above in the LS discussion.
% For the $L^\infty(\rho)$-norm bound note that Monte Carlo integration of $b_j^2$ is a special case of best approximation with $\mcal{V} = \operatorname{span}\braces{b_j}$. Thus, for the same reason as above, $w = \mfrak{K}^{-1} = b_j^{-2}$ is optimal.)}
%
This means that a set of independent samples from the density $w_G^{-1}$ contains, on average, for every basis function $b_k$ the same number of sample points that are optimal for estimating its norm.\footnote{Since a bounded relative integration error for $b_k^2$ is equivalent to a norm equivalence in the one-dimensional linear space $\operatorname{span}(b_k)$, it is not surprising that the corresponding optimal importance sampling densities are the Christoffel functions for $\operatorname{span}(b_k)$.
A formal explanation for why the sum of the individual basis functions remains optimal for least squares approximation in the span of all basis functions $\operatorname{span}(b_1, \ldots, b_d)$ is given in Theorem~7~(6) in~\cite{Trunschke2024}.}
Through this lens, Christoffel sampling can be viewed as a stratified importance sampling for the $L^2(\rho)$-norm.

Assuming that $\beta_k > 0$ for all $k=1,\ldots, d$, it holds that
\begin{equation}
\label{eq:mixture_of_optimals}
    % w_{\hat{G}^{(0)}}(x)^{-1} = \frac{1}{z} \sum_{k=1}^d \lambda_k^{-1} b_k^2(x),
     w^{(0)}(x)^{-1} := w_{\hat{G}^{(0)}}(x)^{-1} = \frac{1}{z^{(0)}} \sum_{k=1}^d \beta_k^{-1} b_k^2(x),
    \qquad
    z^{(0)} := \sum_{k=1}^d \beta_k^{-1} .
\end{equation}
This is a mixture of the same importance sampling densities but with different weights.
For basis functions $b_k$ for which the estimate $\beta_k = \norm{b_k}_{\hat{G}^{(0)}}$ is small, i.e.\ a worse estimate for $\norm{b_k}_G =1$, the optimal sampling densities ${w}_k^{-1} = b_k^2$ are weighted heavier.
This means that a set of independent samples from this density contains, on average, more optimal sample points for the currently worse estimated functions.
%This means that a (stratified) sample from this density contains more optimal sample points for the currently worse estimated functions.
Drawing points $x_i^{(1)}$ from the distribution $ ({w}^{(0)})^{-1}\rho$ for $i=1,\ldots,n^{(1)}$, we can update the Gramian matrix estimate $$\hat{G}^{(1)} = \frac{n^{(0)}}{n^{(0)} + n^{(1)}}\hat{G}^{(0)} + \frac{n^{(1)}}{n^{(0)} + n^{(1)}} \sum_{i=1}^{n^{(1)}} {w}^{(0)}(x_i^{(1)}) B(x_i^{(1)})B(x_i^{(1)})^\intercal,$$
where $n^{(0)}$ is the number of samples used for the first estimate $\hat G^{(0)}$, and $n^{(1)}$  the number of samples used for the correction.
Building on these insights, we propose to define a sequence of Monte Carlo estimates $\hat{G}^{(k)}$, each improving the previous estimate $\hat{G}^{(k-1)}$ by utilising the corresponding optimal sampling density $(w^{(k-1)})^{-1}\rho$.
A deeper intuition for the behaviour of this method is developed in Appendix~\ref{app:update}.

\revision[0]{A related iterative strategy was proposed for discrete $\mathcal{X}$ in~\cite{musco_musco, musco_musco_2} and for general $\mathcal{X}$ in~\cite{herremans2025refinementbasedchristoffelsamplingsquares}.
Unlike our proposed method, which updates the Gramian and uses the corresponding weight function, their approach updates the weight function directly.
While this distinction is negligible for discrete $\mathcal{X}$, our approach makes sampling arguably easier in the continuous setting.}

\revision[0]{
\subsection{Initialisation}
\label{sec:init}

The above reasoning is only valid when $\beta_k > 0$ for all $k=1,\ldots, d$.
Since this is equivalent to having $\operatorname{rk}(\hat G^{(0)}) = d$, we now briefly discuss an initialisation procedure which returns an unbiased rank-$d$ estimate $\hat G^{(0)}$ of the Gramian matrix $G$.

Let
$$
    \tilde{\mathfrak{K}}^{(0)}(x) := \Vert B(x) \Vert_2^2 = \sum_{j=1}^D B_j(x)^2,
$$
and 
$$
    \tilde w^{(0)}(x)^{-1} = (\tilde z^{(0)})^{-1} \tilde{\mathfrak{K}}^{(0)}(x),
    \qquad
    \tilde z^{(0)} = \int \tilde{\mathfrak{K}}^{(0)}(x) \dx[\rho](x).
$$
We first draw a sample $x_1 \sim (\tilde w^{(0)})^{-1} \rho$ and let 
$$
    H^{(1)}
    = \tilde w^{(0)}(x_1) B(x_1) B(x_1)^\intercal.
$$
Assuming  we have generated samples $(x_1, \dots , x_k)$, we define $U_k = \operatorname{span}\{B(x_1), \dots, B(x_k)\}$ and let 
$$
    \tilde{\mathfrak{K}}^{(k)}(x)
    := \Vert B(x) - P_{U_k} B(x) \Vert_2^2
    = \Vert B(x) \Vert_2 - \Vert P_{U_k} B(x) \Vert_2^2 ,
$$
and 
$$
    \tilde z^{(k)} = \int \tilde{\mathfrak{K}}^{(k)}(x) \dx[\rho](x).
$$
Assuming $\tilde z^{(k)}>0$, we can define a probability density  
$$
    \tilde w^{(k)}(x)^{-1} := (\tilde z^{(k)})^{-1} \tilde{\mathfrak{K}}^{(k)}(x),
$$
draw a sample $x_{k+1} \sim (\tilde w^{(k)})^{-1} \rho$ and define 
$$
    H^{(k+1)}
    = \frac{k}{k+1} H^{(k)} + \frac{1}{k+1} \tilde w^{(0)}(x_{k+1}) B(x_{k+1}) B(x_{k+1})^\intercal
    = \frac{1}{k+1} \sum_{i=1}^{k+1} \tilde w^{(0)}(x_i) B(x_i) B(x_i)^\intercal .
$$
The subsequent proposition proves 
that $(x_1 , \dots , x_k)$ has a density with respect to $\rho^{\otimes k}$ (given by a determinant) which is invariant under permutations, and has all marginals equal to $(\tilde w^{(0)})^{-1} \rho$.
\begin{minipage}{\textwidth}
\begin{proposition}
\label{prop:vol_sampling}
    Let $r\le d$ and $(x_1,\ldots, x_r)$ be such that $x_1\sim(\tilde w^{(0)})^{-1} \rho$ and $x_{k+1}\sim (\tilde w^{(k)})^{-1}\rho$, given $(x_1, \ldots, x_k)$, for all $k = 1, \ldots, r-1$.
    Then
    \begin{equation}
    \label{vol_sampling}
        (x_1,\hdots x_r) \sim \det(C(x_1, \hdots, x_r)^\intercal C(x_1, \hdots, x_r) ) \dx[\rho](x_1) \dots \dx[\rho](x_r)
    \end{equation}
    with $C(x_1, \hdots, x_r) = (B(x_1) \dots B(x_r)) \in \mathbb{R}^{D \times r}$.
    In particular, $x_k \sim (\tilde w^{(0)})^{-1} \rho$ for all $1\le k \le r$.
\end{proposition}
\end{minipage}
\begin{proof}
    The density of $(x_1 , \dots , x_r)$ w.r.t.\ $\rho^{\otimes r}$ is given by
    $$
        \tilde w^{(0)}(x_1)^{-1}  \tilde w^{(1)}(x_2)^{-1} \dots \tilde w^{(r-1)}(x_r)^{-1}
        \propto \Vert B(x_1) \Vert^2_2 \Vert B(x_2) - P_{U_1}B(x_2)  \Vert^2_2  \dots \Vert B(x_r) - P_{U_{r-1}}B(x_r)  \Vert^2_2 .
    $$
    We recognize the right-hand side as the square of the volume of the $r$-dimensional parallelotope in $\mathbb{R}^D$ generated by the vectors $(B(x_1), \hdots, B(x_r))$, given by $\det(C(x_1, \hdots, x_r)^\intercal C(x_1, \hdots, x_r) )$~\cite[Proposition~1]{Gover2010Jul}.
    This proves the first claim.
    For the second claim, we note that the density~\eqref{vol_sampling} is invariant under permutations, which implies that all marginals are equal to the distribution of $x_1$, which is $(\tilde w^{(0)})^{-1} \rho$.  
\end{proof}

In the next proposition, we prove that when $\tilde z^{(k)} = 0$ (or equivalently
$\tilde{\mathfrak{K}}^{(k)}(x)=0$ $\rho$-almost everywhere), then $\operatorname{rk}( H^{(k)}) = d$ almost surely, which provides a stopping criterium for the initialisation procedure.
Moreover, we show that this happens almost surely after $d$ steps and that the resulting Gramian $H^{(d)}$ has the same kernel as $G$.  
\begin{proposition}
\label{prop:increasing_rank}
    The distribution~\eqref{vol_sampling} is well-defined (i.e.\ $\tilde z^{(k)} > 0$) for all $k<d$.
    Moreover, almost surely $\operatorname{rk}(H^{(k)}) = k$ for any $k\le d$ and almost surely $\tilde z^{(d)} = 0$ and $\ker(H^{(d)}) = \ker(G)$. 
    % Almost surely, it holds $\operatorname{rk}(H^{(k)}) = k$ and $\operatorname{Ker}(H^{(k)}) \supset  \operatorname{Ker}(G)$ for any $k\ge d$.  It holds $\tilde z^{(k)} = 0$ if and only if $\operatorname{rk}(H^{(k)}) = d$ and $\operatorname{Ker}(H^{(k)}) = \operatorname{Ker}(G)$. If $\operatorname{rk}(H^{(k)}) < d$ then $\tilde z^{(k)} >0$, so that $(\tilde w^{(k)})^{-1} \rho$ is a well defined probability measure.
    % % and 
    % almost surely, 
    % $$
    %      \operatorname{rk}(H^{(k+1)}) = \operatorname{rk}( H^{(k)}) +1 = k.
    % $$
\end{proposition}
\begin{proof}
    Note that $
        \tilde{\mathfrak{K}}^{(k)}(x)
        = \|(I - P_{U_k}) B(x)\|_2^2
        = \inner{I - P_{U_k}, B(x)B(x)^\intercal}_{\mathrm{Fro}}
    $ and $I - P_{U_k} = P_{\ker(H^{(k)})}$ and, therefore,
    $$
        \tilde z^{(k)}
        = \int \tilde{\mathfrak{K}}^{(k)}(x) \dx[\rho(x)]
        = \inner{I - P_{U_k}, G}_{\mathrm{Fro}}
        = \operatorname{tr}((I - P_{U_k}) G (I - P_{U_k}))
        = \|(I - P_{U_k}) G (I - P_{U_k})\|_{\mathrm{nuc}} .
    $$
    This implies that $\tilde z^{(k)} = 0$ if and only if $(I - P_{U_k}) G (I - P_{U_k}) = 0$, which is equivalent to $\operatorname{Im}(I - P_{U_k}) \subseteq \ker(G)$.
    Since $\operatorname{Im}(I - P_{U_k}) = \ker(P_{U_k}) = \ker(H^{(k)})$, this condition simplifies to
    $$
        \tilde z^{(k)} = 0
        \quad\Leftrightarrow\quad
        \ker(H^{(k)}) \subseteq \ker(G) .
    $$
    Therefore, $\tilde z^{(k)} = 0$ requires $D - k \le \dim\ker(H^{(k)}) \le \dim\ker(G) = D - d$, i.e.\ $k \ge d$.
    % Let $k < d$ and note that $\tilde z^{(k)} = 0$ if and only if $\tilde{\mathfrak{K}}^{(k)}(x) = 0$ for $\rho$-almost all $x$.
    % Therefore, $\tilde z^{(k)} = 0$ if and only if $B(x) \in U_k = \operatorname{im}(H^{(k)}) = \ker(H^{(k)})^\perp$ for $\rho$-almost all $x$.
    % Since $y\in \ker(H^{(k)})$ implies $y^\intercal G y = \int (y^\intercal B(x))^2 \dx[\rho](x) = 0$, which implies $y\in \ker(G)$, we deduce that $\ker(H^{(k)}) \subseteq \ker(G)$ and, consequently, $\ker(H^{(k)})^\perp \supseteq \ker(G)^\perp$.
    % This means that $\tilde z^{(k)} = 0$ implies $B(x) \in \ker(G)^\perp$ for $\rho$-almost all $x$.
\\
    Now let $k\le d$ and recall that $(x_1,\ldots, x_k)\sim \det(C(x_1,\dots, x_k)^\intercal C(x_1,\dots, x_k))\dx[\rho](x_1) \dots \!\dx[\rho](x_k)$ by Proposition~\ref{prop:vol_sampling}.
    It holds almost surely that $\det(C(x_1,\dots, x_k)^\intercal C(x_1,\dots, x_k)) > 0$ and $\tilde w(x_i) > 0$ for all $1\le i \le k$, which implies that the vectors $(B(x_1), \hdots, B(x_k))$ are linearly independent and  the matrix $H^{(k)}$ has rank $k$.
\\
    Finally, note that $v^\intercal G v = 0$ implies $v^\intercal B(x) = 0$ for $\rho$-almost all $x$ and, in particular, $v^\intercal B(x_i) = 0$ almost surely for all $i=1,\ldots,k$.
    From this, we conclude that $v^\intercal H^{(k)} v = 0$ almost surely.
    Since almost surely $\ker(H^{(k)})\supseteq \ker(G)$, we can conclude the final claim.
    % $$
    %     \tilde z^{(k)} = 0
    %     \quad\Leftrightarrow\quad
    %     \ker(H^{(k)}) = \ker(G) .
    % $$
    % In particular, this means that almost surely $\operatorname{rk}(H^{(d)}) = d$.
    % Then for $k \le d$, it holds almost surely $\dim(\operatorname{Ker}(H^{(k)})) = D-k \ge D-d = \dim(\operatorname{Ker}(G))$.
    % Therefore, almost surely, $\tilde z^{(k)} = 0$ if and only if $\operatorname{Ker}(H^{(k)}) = \operatorname{Ker}(G)$ and $k=d$. 
\end{proof}
From the above result, we deduce an initialisation procedure which applies the above algorithm until we detect
$\tilde z^{(k)} = 0$, when $k=d$, and returns  $\hat G^{(0)} = H^{(d)}$, satisfying $\ker(\hat G^{(0)}) = \ker(G)$.
Note that detecting that $\tilde z^{(k)} = 0$ may be a difficult task, especially when the dictionary contains localised functions in high dimension.

Also note that, in practice, we have to estimate the normalisation constants $\tilde z^{(k)}$ from samples.
}

\subsection{Iterative refinement}

% {\color{red}
% \begin{equation}
% \label{eq:hatG_1}
%     \hat{G}^{(1)} := \tfrac12\hat{G}^{(0)} + \tfrac12\hat{G}^{(\sfrac{1}{2})},
%     \qquad
%     \hat{G}^{(\sfrac12)} := \frac1n \sum_{i=1}^n w_{\hat{G}^{(0)}}(x^{(1)}_{i}) B(x^{(1)}_{i}) B(x^{(1)}_{i})^\intercal,
%     \qquad
%     x^{(1)}_{1},\ldots,x^{(1)}_{n} \sim w_{\hat{G}^{(0)}}^{-1}\rho
% \end{equation}
% to improve the estimate $\hat{G}^{(0)}$.
% Starting from the initial estimate~\eqref{eq:hatG_0} and iterating the procedure delineated in equation~\eqref{eq:hatG_1}, we can define the sequence of estimates
% \begin{equation}
% \label{eq:hatG_k+1}
% \begin{aligned}
%     \hat{G}^{(k+1)} &:= \tfrac{k}{k+1}\hat{G}^{(k)} + \tfrac1{k+1}\hat{G}^{(k+\sfrac{1}{2})}, \\
%     \hat{G}^{(k+\sfrac12)} &:= \frac1n \sum_{i=1}^n w_{\hat{G}^{(k)}}(x^{(k+1)}_{i}) B(x^{(k+1)}_{i}) B(x^{(k+1)}_{i})^\intercal,
%     \qquad
%     x^{(k+1)}_{1},\ldots,x^{(k+1)}_{n} \sim w_{\hat{G}^{(k)}}^{-1}\rho
% \end{aligned}
% \end{equation}
% for all $k\ge1$.
% }

Given an initial estimate
% \comment{A: here $\hat{G}^{(0)} = \hat{G}^{(0+\sfrac12)} = \hat{G}^{(1)} $, isn't it ? Not consistant with the first algo but ok !}
\begin{equation}
    \hat{G}^{(0)}
    % := \hat{G}^{(\sfrac12)}
    := \frac{1}{n^{(0)}} \sum_{i=1}^{n^{(0)}} w^{(0)}(x^{(0)})B\pars{x^{(0)}_{i}} B\pars{x^{(0)}_{i}}^\intercal
\end{equation}
with $\ker(\hat G^{(0)}) = \ker(G)$, we now define sequences of estimates $\hat{G}^{(k)}$. 
% and $\hat{z}^{(k)}$.
%     For this, let $\hat{n}^{(k)}$ be a sequence of $\mathbb{N}$-valued random variables, where each $\hat{n}^{(k)}$ is conditionally independent of $\hat{z}^{(k-1)}$ given $\hat{G}^{(k-1)}$.
%     Moreover, let $\hat{m}$ be an $\mathbb{N}$-valued random variable independent of $\hat{n}^{(k)}$, $\hat{z}^{(k)}$ and $\hat{G}^{(k)}$ for all $k\in\mathbb{N}$.
%
%\modifan{and assume that $\hat{G}^{(0)} \neq 0$.}
For every $k\ge1$, we define the sampling measure
$$
	\mu^{(k)} := \frac{\mathfrak{K}^{(k)}}{z^{(k)}} \rho
$$
with the inverse Christoffel function and normalisation constant given by
\begin{equation}
\label{eq:def_z_z_bar}
    \mathfrak{K}^{(k)} := \mathfrak{K}_{\hat{G}^{(k)}}
    \qquad\text{and}\qquad
	z^{(k)}
	:= \int\mathfrak{K}^{(k)}\dx[\rho]
	% = \int \inner{(\hat{G}^{(k)})^+, BB^\intercal}_{\mathrm{Fro}} \dx[\rho]
	= \inner{(\hat{G}^{(k)})^+, G}_{\mathrm{Fro}} .
\end{equation}
%{\color{red}
% This only works when $\ker(\hat{G}^{(k)}) = \ker(G)$!
% When this is not the case, we utilise the spectral decomposition $\hat{G}^{(k)} = U\Lambda U^\intercal$ with $\Lambda\in\mathbb{R}^{D\times D}$, to define, analogously to~\eqref{eq:mixture_of_optimals}
% $$
%     \tilde{\mathfrak{K}}^{(k)} = \sum_{j:\Lambda_{jj}=0} (U^\intercal B)_j^2
%     \quad\text{and}\quad
%     \tilde{w}^{(k)} := (\tilde{\mathfrak{K}}^{(k)})^{-1} \int \tilde{\mathfrak{K}}^{(k)} \dx[\rho] .
% $$
% We then draw $\tilde{x}^{(k)}\sim (\tilde{w}^{(k)})^{-1}\rho$ and define $\hat{G}^{(k+1)} = \frac{k}{k+1}\hat{G}^{(k)} + \frac{1}{k+1} \tilde{w}^{(k)}(\tilde{x}^{(k)}) B(\tilde{x}^{(k)})B(\tilde{x}^{(k)})^\intercal$.
%}
The corresponding weight function is given by
$$
	w^{(k)} := \frac{z^{(k)}}{\mathfrak{K}^{(k)}} .
$$
Since the exact normalisation constant $z^{(k)}$ is inaccessible, we replace it with the Monte Carlo estimate
\begin{equation}
\label{eq:def_zhat}
    \hat{z}^{(k)} := \frac{1}{m} \sum_{i=1}^m \mathfrak{K}^{(k)}(y^{(k)}_i),
    \qquad
    y^{(k)}_1,\ldots,y^{(k)}_m \overset{\mathrm{i.i.d.}}{\sim} \rho .
\end{equation}
Doing so results in the estimated weight function
$$
	\hat{w}^{(k)} := \frac{\hat{z}^{(k)}}{\mathfrak{K}^{(k)}} .
$$
With this estimate in hand, we can define the sequence of Gramian estimates
\begin{equation}
\label{eq:Gramian_iterates_def}
    \hat{G}^{(k+1)} := \tfrac{k}{k+1}\hat{G}^{(k)} + \tfrac1{k+1}\hat{G}^{(k+\sfrac{1}{2})}, 
\end{equation}
with 
\begin{equation}
    \hat{G}^{(k+\sfrac12)} := \frac1{n^{(k+1)}} \sum_{i=1}^{n^{(k+1)}} \hat{w}^{(k)}(x^{(k+1)}_{i}) B(x^{(k+1)}_{i}) B(x^{(k+1)}_{i})^\intercal,
    \qquad
    x^{(k+1)}_{1},\ldots,x^{(k+1)}_{n^{(k+1)}} \sim \mu^{(k)},
\end{equation}
for all $k\ge1$.

\revision[0]{Note that the number $m$ of samples used for estimating the normalisation constant and the   $n^{(k)}$ may be random. In this case, we denote them by $\hat m$ and $\hat n^{(k)}$ respectively. We  assume that ${n}^{(k)}$ is conditionally independent of $\hat{z}^{(k-1)}$ given $\hat{G}^{(k-1)}$, and that $\hat{m}$ is independent of ${n}^{(k)}$, $\hat{z}^{(k)}$ and $\hat{G}^{(k)}$ for all $k\in\mathbb{N}$.} 
% {\color{red}
% To bound the variance of $\hat{z}_{\hat{G}_k}$, we assume that there exists a measure $\nu$ such that $\rho \ll \nu$ is absolutely continuous with respect to $\nu$.
% Denoting the Radon--Nikodym derivative by $\gamma := \frac{\dx[\rho]}{\dx[\nu]}$, we define
% \begin{equation}
% \label{eq:zhat}
%     \hat{z}_{H} := \frac{1}{m} \sum_{i=1}^m \gamma\pars{x_i} \mfrak{K}_{H}\pars{x_i}
%     \qquad\text{with}\qquad
%     x_i\overset{\mathrm{i.i.d.}}{\sim} \nu .
% \end{equation}
% % \begin{remark}
% %     Although the choice of $m\in\mathbb{N}$ in definition~\eqref{eq:zhat} should intuitively influence the quality of the estimate $\hat{z}_H$, we do not see this dependence in our derivation.
% %     This is due to the naive bounds that we utilise for the sake of simplicity.
% % \end{remark}
% }

% {\color{Mahogany}
% \subsection{Termination}
% \label{sec:termination}
% TODO: A brief section on termination?
% }

\section{Theory}
\label{sec:convergence_of_Gk}

In this section, we analyse the convergence of the iterates $\hat{G}^{(k)}$ to the true Gramian matrix $G$.
Instead of the close-to-optimal bounds that hold with high probability, we restrict our analysis to the easier setting of the estimate's mean and variance.
We show that the iterates converge in expected mean squared error and leave the more difficult question of tight non-asymptotic rates to future research.
% Theorem~\ref{thm:convergence_TV} ensures that the convergence of the Gramian matrix estimates induces a convergence of the sampling measures in total variation distance.
Although the resulting convergence rate seems slow at first, recall that we do not require an accurate estimate of $G$ but that it suffices to obtain an estimate that satisfies the uniform framing~\eqref{eq:framing}.

\revision[0]{Using the initialisation procedure from section~\ref{sec:init}, we can subsequently assume that $\ker(\hat G^{(k)}) = \ker(G)$ for all $k\in\mathbb{N}$ (cf.~Proposition~\ref{prop:increasing_rank}).}

\subsection{Convergence of \texorpdfstring{$\boldsymbol{\hat{G}^{(k)}}$}{G\^(k)} to \texorpdfstring{$\boldsymbol{G}$}{G}}

\revision[0]{Let $\lambda_{\mathrm{min}>0}\pars{G}$ denote the smallest positive eigenvalue of $G$, and let $\kappa_{>0}(G) := \frac{\lambda_{\mathrm{max}}(G)}{\lambda_{\mathrm{min}>0}(G)}$ be the effective condition number of $G$.} 

\begin{theorem}
\label{thm:mean_and_var}
    It holds for all $k\ge 0$ that $\mathbb{E}\big[\hat{G}^{(k)}\big] = G$.
    Moreover, for any fixed integer $N\ge 1$, if $k \ge 1$ and $\hat n^{(k+1)} \ge N\hat{z}^{(k)}$, it holds for all $\hat{m}$ that
    \begin{equation}
    \label{eq:var_zhat}
        \operatorname{Var}\big(\hat{G}^{(k)}_\alpha\big)
        \le \frac{W^2}{kN},
        \qquad
        W^2 := \kappa_{>0}(G) G_{\alpha}^2,
    \end{equation}
    for all multi-indices $\alpha \in \{1,\dots,D\}^2$.
    If $k \ge 0$, $\hat n^{(k+1)} \ge \frac{2N}{\lambda_{\mathrm{min}>0}(\hat{G}^{(k)})}$ and $\hat{m} \ge V$ with $V := \mathbb{E}\bracs*{\norm{b(x_1)b(x_1)^\intercal - I}_{2}^2} \le \|\mathfrak{K}_G\|_{L^2(\rho)}^2 - 2(d-1)$, it holds that
    \begin{equation}
    \label{eq:var_lambda_min}
    	\operatorname{Var}\big(\hat{G}^{(k)}_\alpha\big)
    	\le \frac{W^2}{kN},
        \qquad
        W^2 := \operatorname{tr}(G) G_{\alpha}^2 .
    \end{equation}
\end{theorem}

To prove the above theorem, we first prove a sequence of intermediate results.
The following four lemmas prove that the iterates are unbiased and uncorrelated, and provide bounds for their variance.
% express their variance in terms of the variance of the update matrices $\hat{G}^{(k+\sfrac12)}$.

\begin{lemma}
\label{lem:G_mean}
    It holds that $\mathbb{E}\big[\hat{G}^{(k)}\big] = \mathbb{E}\big[\hat{G}^{(k+\sfrac12)}\big] = G$ for all $k\ge 0$.
\end{lemma}
\begin{proof}
    The proof of this lemma   can be found in Appendix~\ref{app:proof:lem:G_mean}.
\end{proof}

\begin{lemma}
\label{lem:G_cov}
    Let $j\ne k\in\mathbb{N}$ and $\alpha,\beta\in \{1, \ldots, D\}^2$ be arbitrary multi-indices.
    Then $\operatorname{Cov}\big(\hat{G}^{(k+\sfrac12)}_{\alpha}, \hat{G}^{(j+\sfrac12)}_{\beta}\big) = 0$.
\end{lemma}
\begin{proof}
    The proof of this lemma  can be found in Appendix~\ref{app:proof:lem:G_cov}.
\end{proof}

\begin{lemma}
\label{lem:var_G_update}
    Let $k\in\mathbb{N}$ and $\alpha\in \{1, \ldots, D\}^2$ be an arbitrary multi-index.
    % Then $\operatorname{Var}\big(\hat{G}^{(k+1)}_{\alpha}\big) = \tfrac{k^2}{(k+1)^2} \operatorname{Var}\big(\hat{G}^{\pars{k}}_\alpha\big)
    % + \tfrac{1}{(k+1)^2} \operatorname{Var}\big(\hat{G}^{\pars{k+\sfrac12}}_\alpha\big)
    % $.
    Then
    $$
        \operatorname{Var}\big(\hat{G}^{(k+1)}_{\alpha}\big)
        = \frac{1}{(k+1)^2} \sum_{j=0}^k \operatorname{Var}\big(\hat{G}^{\pars{j+\sfrac12}}_\alpha\big) .
    $$
\end{lemma}
\begin{proof}
    The proof of this lemma   can be found in Appendix~\ref{app:proof:lem:var_G_update}.
\end{proof}

\begin{lemma}
\label{lem:var_G}
    If $\hat{n}^{(k+1)} \ge N \hat{z}^{(k)}$, it holds that
    $$
    	\operatorname{Var}\big(\hat{G}^{(k+\sfrac12)}_\alpha\big)
    	\le \frac{1}{N} \kappa_{>0}(G) G_\alpha^2 .
    $$
    Moreover, if $\hat{n}^{(k+1)}$ is independent of $\hat{z}^{(k)}$, given $\hat{G}^{(k)}$ and $\hat{m}$, it holds that
    $$
    	\operatorname{Var}\big(\hat{G}^{(k+\sfrac12)}_\alpha\,\big|\,\hat{n}^{(k+1)}\big)
    	\le \frac{1}{\hat n^{(k+1)}} \pars*{\mathbb{E}\bracs*{\pars*{\frac{V}{\hat{m}} + 1} z^{(k)} \hat{G}^{(k)}_{\alpha} \,\bigg|\, \hat{n}^{(k+1)}} G_{\alpha} - G_\alpha^2} .
    $$
    In particular, if $\hat{m} \ge V$ with $V := \mathbb{E}\bracs*{\norm{b(x_1)b(x_1)^\intercal - I}_{2}^2}$ and $\hat{n}^{(k+1)} \ge \frac{2N}{\lambda_{\mathrm{min}>0}(\hat{G}^{(k)})}$, it holds that
    $$
    	\operatorname{Var}\big(\hat{G}^{(k+\sfrac12)}_\alpha\big)
    	\le \frac{1}{N}\operatorname{tr}(G) G_{\alpha}^2 .
    $$
\end{lemma}
\begin{proof}
    The proof of this lemma can be found in Appendix~\ref{app:proof:lem:var_G}.
\end{proof}

\begin{proof}[Proof of Theorem~\ref{thm:mean_and_var}]
    \revision[0]{
    The first claim follows from Lemma~\ref{lem:G_mean}.
    Moreover, Proposition~\ref{prop:increasing_rank} ensures that almost surely $\ker(\hat{G}^{(k)}) = \ker(G)$ such that Lemmas~\ref{lem:var_G_update} and~\ref{lem:var_G} provide the claimed variance bounds.
    % Proof:
    % \begin{align}
    %     v_{k+1}
    %     &= \frac{k^2}{(k+1)^2} v_k + \frac{1}{(k+1)^2} c \\
    %     &= \frac{k^2}{(k+1)^2}\bigg(\frac{(k-1)^2}{k^2} v_{k-1} + \frac{1}{k^2} c\bigg) + \frac{1}{(k+1)^2} c \\
    %     &= \frac{(k-1)^2}{(k+1)^2} v_{k-1} + \frac{1}{(k+1)^2} c + \frac{1}{(k+1)^2} c \\
    %     &= \frac{k+1}{(k+1)^2} c \\
    %     &= \frac{1}{k+1} c
    %     .
    % \end{align}
    % }
    }
\end{proof}

\begin{remark}[Choosing $\hat{m}$ and $\hat{n}^{(k+1)}$ in Theorem~\ref{thm:mean_and_var}]
\label{rmk:choosing_m_and_n}
    \revision[0]{
    While the second bound in Theorem~\ref{thm:mean_and_var} is typically smaller than the first, it also requires larger sample sizes $\hat m$ and $\hat n^{(k+1)}$.
    The dependence on $\hat m$ is intuitively understandable, since a high variance of the normalisation constants $\hat z^{(j)}$ results in some summands carrying more weight than others.
    Although the dependence on $\hat n^{(k+1)}$ is not as intuitive, it is also not as strong.
    Since $\lambda_{\mathrm{min}>0}(\hat{G}^{(k)})^{-1}$ is typically very large, it may seem more prudent to choose $\hat{n}^{(k+1)} \ge N \hat{z}^{(k)}$ rather than $\hat{n}^{(k+1)} \ge \frac{2N}{\lambda_{\mathrm{min}>0}(\hat{G}^{(k)})}$.
    However, since $\ker(\hat{G}^{(k)}) = \ker(G)$ there exist constants $c_k$ and $C_k$ such that $\hat{G}^{(k)}$ satisfies the matrix framing
    $$
        C_k^{-1} G \preceq \hat{G}^{(k)} \preceq c_k^{-1} G .
    $$
    The tightest such constants are given by $c_k^{-1} = \lambda_{\mathrm{max}}(G^{+1/2} \hat{G}^{(k)} G^{+1/2})$ and $C_k = \lambda_{\mathrm{max}}((\hat{G}^{(k)})^{+1/2} G (\hat{G}^{(k)})^{+1/2})$.
    Since $z_{\hat{G}^{(k)}} = \operatorname{tr}((\hat{G}^{(k)})^{+1/2} G (\hat{G}^{(k)})^{+1/2})$, it follows that
    $$
        C_k \le z_{\hat{G}^{(k)}} \le d C_k .
    $$
    % Proof:
    % Note that $C^{-1}G \preceq H$ is equivalent to $C^{-1} I \preceq G^{+1/2} H G^{+1/2}$, which is equivalent to $C^{-1} \le \lambda_{\mathrm{min}}(G^{+1/2} H G^{+1/2})$.
    % Equivalently, $C^{-1}G \preceq H$ is equivalent to $H^{+1/2} G H^{+1/2} \preceq C I$, i.e.\ $C \ge \lambda_{\mathrm{max}}(H^{+1/2} G H^{+1/2})$.
    % If $C$ is chosen tight, we have equality.
    % Since $z_H = \operatorname{tr}(H^{+1/2} G H^{+1/2})$, we have, obviously
    % $$
    %     C \le z_H \le d C .
    % $$
    This means that $z^{(k)}$ is of the same order as the (potentially large) framing constant $C_k$.
    Therefore, the sample size $\hat{n}^{(k+1)} \ge N \hat{z}^{(k)}$ may still be huge.
    Nevertheless, the experiments in section~\ref{sec:experiments_intro} indicate that the choice $\hat{n}^{(k+1)} \ge \hat{z}^{(k)}$ may not be necessary.
    }
\end{remark}

\begin{remark}[Relation to standard Monte Carlo]
    \revision[0]{
    Note that the standard Monte Carlo estimate
    $$
        \hat{G}^{(\mathrm{MC}, n)}
        := \frac1n \sum_{i=1}^n B(y^{(\mathrm{MC})}_{i}) B(y^{(\mathrm{MC})}_{i})^\intercal,
        \qquad
        y^{(\mathrm{MC})}_{1},\ldots,y^{(\mathrm{MC})}_{n} \sim \rho ,
    $$
    is unbiased with $\operatorname{Var}(G^{(\mathrm{MC},n)}_{lm}) = \tfrac{1}{n} (\mathfrak{G}_{lm} - G_{lm}^2)$ where $\mfrak{G}_{lm} := \int B_l^2B_m^2 \dx[\rho]$.
    This means that Theorem~\ref{thm:mean_and_var} only reduces the variance significantly when $\mathfrak{G}_{lm}$ is large.
    Moreover, it means that the rate of convergence does not change.
    % Moreover, the sample size $\hat m\ge V$ necessary for equation~\eqref{eq:var_lambda_min} is of the same order as $\mathfrak{G}_{lm}$.
    We believe that this is only an artifact of our initial, crude analysis and present a refined analysis in an idealised setting in Appendix~\ref{app:update}.
    The numerical experiments in Section~\ref{sec:experiments_intro} indeed confirm this intuition of faster convergence.
    }
\end{remark}

\subsection{Matrix framing with high probability}

Recall that the objective of this paper is to obtain a sampling density for which the sample size required by Theorem~\ref{thm:sample_size_bounds} is minimal.
Equation~\eqref{eq:suboptimality_factor} reduces this problem to minimising the oversampling factor $\tfrac{C}{c}$ corresponding to the uniform framing~\eqref{eq:framing} and Lemma~\ref{lem:framing_equivalence} provides an equivalence of the uniform framing~\eqref{eq:framing} to a matrix framing.
We hence want to show that the estimated Gramian matrices $\hat{G}^{(k)}$ produced by our algorithm exhibit a tight matrix framing with high probability.
To do this, we employ a matrix version of Chebyshev's inequality, first published in~\cite{985947}, and recalled in the subsequent lemma.
This allows us to derive probability bounds from the variance bounds derived in the preceding sections.

% Recall that Theorem~\ref{thm:sample_size_bounds} implies that bounding the sample size required for a stable approximation requires bounding the oversampling factor $\tfrac{C}{c}$ from equation~\eqref{eq:suboptimality_factor}.
% % Let $\tfrac{C}{c}(k)$ denote the suboptimality factor in iteration $k$ of the algorithm.
% Since Lemma~\ref{lem:framing_equivalence} provides an equivalence of the uniform framing~\eqref{eq:framing} and the matrix framing, this can be reduced to finding a probability bound for the corresponding matrix framing. % in order to bound $\tfrac{C}{c}(k)$.
% First published in~\cite{985947}, the matrix version of Chebyshev's inequality, presented in the subsequent lemma, can be used to derive probability bounds for the matrix framing from the variance bounds derived so far.

\begin{lemma}[Matrix Chebyshev's inequality~\cite{Tu_2017}]
\label{lemma:mci}
    Let $A \succ 0$, and let $X$ be a symmetric random matrix.
    Then,
    $$
        \mathbb{P}\bracs{ \abs{X} \not\preceq A } \le \operatorname{tr}(\mathbb{E}[X^2] A^{-2}) \:,
    $$
    \revision[0]{where $\abs{X}$ is defined by applying $\abs\bullet$ to the eigenvalues of $X$.}
\end{lemma}

% Let Y be a random matrix and A = E[Y]. Then
% $$
% E[(Y-A)^2] = E[(Y-A)(Y-A)] = E[Y^2] - E[Y]A - AE[Y] + A^2 = E[Y^2] - A^2.
% $$
% Consequently,
% $$
% Tr[E[(Y-A)^2] A^{-2}] = Tr[E[Y^2] A^{-2}] - Tr[I] = E[Tr[A^-1Y^2A^{-1}]] - d = E[ |A^{-1}Y|_F ] - d.
% $$

Given a matrix $\hat{G}^{(k)}$, we define $C(\hat{G}^{(k)}) \in [1,\infty)$ and $c(\hat{G}^{(k)}) \in [0, 1]$ as the tightest constants satisfying the matrix framing $C(\hat{G}^{(k)})^{-1}G \preceq \hat{G}^{(k)} \preceq c(\hat{G}^{(k)})^{-1}G$.
Moreover, we define the oversampling factor
$$
    \gamma(\hat{G}^{(k)}) := \tfrac{C(\hat{G}^{(k)})}{c(\hat{G}^{(k)})} \:.
$$
This factor is bounded in the subsequent theorem.

\begin{theorem}
\label{thm:high_probability}
    \revision[0]{Let $\hat{m}$, $\hat{n}^{(k+1)}$, $N$ and $W$ be as in Theorem~\ref{thm:mean_and_var}.}
    % Suppose that $\hat{G}^{(k)}$ is an unbiased estimate of $G$ with $\operatorname{Var}(\hat{G}^{(k)}_{jl}) \le \frac{\alpha_1^2}{kn}$ for all $1 \le j,l\le D$.
    Then with $\alpha := DW \lambda_{\mathrm{min}>0}(G)^{-1} (1-p)^{-1/2}$ it holds 
    \begin{equation}
    \label{eq:gamma_bound}
        \gamma(\hat{G}^{(k)})
        \le \frac{\sqrt{kN} + \alpha}{\sqrt{kN} - \alpha}
    \end{equation}
    with a probability of at least $p$, when $\sqrt{kN} > \alpha$.
\end{theorem}
\begin{proof}
    Since $\gamma(\hat{G}^{(k)}) \ge 1$, we start by bounding the failure probability
    $$
        \mathbb{P}\bracs{\gamma(\hat{G}^{(k)}) > \gamma^*}
    $$
    for any $\gamma^* > 1$.
    The probability of the complementary event $\gamma(\hat{G}^{(k)}) \le \gamma^*$ is always larger than having $C(\hat{G}^{(k)}) \le C^*$ and $c(\hat{G}^{(k)})\ge c^*$ for some $0<c^*<1<C^*$ satisfying $\tfrac{C^*}{c^*} \le \gamma^*$, i.e.
    $$
        \mathbb{P}\bracs{\gamma(\hat{G}^{(k)}) \le \gamma^*}
        \ge \mathbb{P}\bracs{(C^*)^{-1}G \preceq \hat{G}^{(k)} \preceq (c^*)^{-1}G} \:.
    $$
    Moreover, since $0<c^*<1<C^*$, there exists (more than one) $\delta^*\in(0,1)$ such that $(C^*)^{-1} \le 1-\delta^*$ and $(c^*)^{-1} \ge 1+\delta^*$, and it holds that
    $$
        (1-\delta^*)G \preceq \hat{G}^{(k)} \preceq (1+\delta^*)G
        \quad\Rightarrow\quad
        (C^*)^{-1}G \preceq \hat{G}^{(k)} \preceq (c^*)^{-1}G .
    $$
    By the monotonicity of probability, we thus have
    $$
        \mathbb{P}\bracs{\gamma(\hat{G}^{(k)}) \le \gamma^*}
        \ge \mathbb{P}\bracs{(C^*)^{-1}G \preceq \hat{G}^{(k)} \preceq (c^*)^{-1}G}
        \ge \mathbb{P}\bracs{(1-\delta^*)G \preceq \hat{G}^{(k)} \preceq (1+\delta^*) G} .
    $$
    Using the equivalence
    $$
        % \tfrac{C}{c}(k) \le \gamma^*
        % \quad\Leftarrow\quad
        % C^{-1}G \preceq \hat{G}^{(k)} \preceq c^{-1}G
        % \quad\Leftarrow\quad
        (1-\delta^*)G \preceq \hat{G}^{(k)} \preceq (1+\delta^*)G
        \quad\Leftrightarrow\quad
        \abs{\hat{G}^{(k)}-G} \preceq \delta^* G
    $$
    together with the monotonicity of probability and applying matrix Chebyshev's inequality yields
    % {\color{red}I would rather have implications of the form delta-framing -> c(k)-C(k)-framing -> C(k)/c(k) <= x.
    % By monotonicity of probability, we thus get P[delta-framing] <= P[c(k)-C(k)-framing] <= P[C(k)/c(k) <= x], or, alternatively,
    % P[C(k)/c(k) > x] <= P[not delta-framing].
    % Since we want delta-framing -> c-C-framing, we need }
    $$
        % \mathbb{P}\big[\neg (C^{-1}G \preceq \hat{G}^{(k)} \preceq c^{-1}G)\big]
        \mathbb{P}\big[\gamma(\hat{G}^{(k)}) > \gamma^*]
        \le \mathbb{P}\big[ \abs{\hat{G}^{(k)}-G} \succ \delta^* G \big]
        \le (\delta^*)^{-2} \operatorname{tr}\pars*{ \mathbb{E}\big[(\hat{G}^{(k)}-G)^2\big] (G^{+})^2 } \:.
    $$
    Since $(\hat{G}^{(k)}-G)^2$ is positive semi-definite, we can apply Lemma~\ref{lem:inner_bounds} to bound
    \begin{align}
        \operatorname{tr}\pars*{ \mathbb{E}\big[(\hat{G}^{(k)}-G)^2\big] (G^{+})^2 }
        &\le \lambda_{\mathrm{min}>0}(G)^{-2} \operatorname{tr}\pars*{\mathbb{E}\big[(\hat{G}^{(k)}-G)^2\big]} \\
        &= \lambda_{\mathrm{min}>0}(G)^{-2}\: \mathbb{E}\bracs*{\operatorname{tr}\big((\hat{G}^{(k)}-G)^2\big)} \\
        &= \lambda_{\mathrm{min}>0}(G)^{-2}\: \mathbb{E}\bracs*{\norm{\hat{G}^{(k)}-G}_{\mathrm{Fro}}^2} ,
    \end{align}
    where the equalities follow by linearity of the trace and symmetry of $\hat{G}^{(k)}-G$, respectively.
    Now, observe that
    $$
        \mathbb{E}\bracs*{\norm{\hat{G}^{(k)}-G}_{\mathrm{Fro}}^2}
        = \sum_{j,l=1}^D \mathbb{E}\bracs*{\big(\hat{G}^{(k)}_{jl} - G_{jl}\big)^2}
        = \sum_{j,l=1}^D \operatorname{Var}\big(\hat{G}^{(k)}_{jl}\big)
        % \le \frac{D^2\alpha_1^2}{kn} \sum_{j,l=1}^D \bracs*{\pars*{1 + \tfrac{V}{m}} \tfrac{\norm{G}_*}{c_k} \mfrak{G}_{jl} - G_{jl}^2}
        % =: \frac{\alpha_1^2}{kn} \:.
        \le \frac{D^2W^2}{kN} \:.
    $$
    Combining the preceding bounds yields
    $$
        \mathbb{P}\big[\gamma(\hat{G}^{(k)}) > \gamma^*\big]
        \le (\delta^*)^{-2} \lambda_{\mathrm{min}>0}(G)^{-2} \frac{D^2W^2}{kN}
        =: \frac{\alpha_{\mathrm{tmp},1}^2}{kN(\delta^*)^2} \:.
    $$
    % For the framing to be satisfied with fixed probability $p\in(0,1)$, we thus require
    % $$
    %     \mathbb{P}\big[C^{-1}G \preceq \hat{G}^{(k)} \preceq c^{-1}G\big]
    %     \le 1-p \le \frac{\alpha_2^2}{kn(\delta^*)^2} \:.
    % $$
    % Solving for $\delta$ gives the final relation
    % {\color{red}Wait! $\delta$ has nothing to do with C and c and (directly), which has nothing to do with the suboptimality factor... This is unsatisfactory...}
    Recall that this bound is valid for any $\delta^*\in(0,1)$ satisfying $(C^*)^{-1} \le 1-\delta^*$ and $(c^*)^{-1} \ge 1+\delta^*$.
    The largest $\delta^*$ adhering to these constraints is given by
    $$
        \delta^* := \min\braces{1-(C^*)^{-1}, (c^*)^{-1}-1} \:.
    $$
    Recall, moreover, that the preceding probability bound is valid for any $0<c^*<1<C^*$ satisfying $\tfrac{C^*}{c^*} \le \gamma^*$.
    Choosing $c^*$ and $C^*$ as to maximise $\delta^*$ yields
    $$
        \delta^*
        := \max_{c^*\in(0,1)} \max_{C^*\in(1, \gamma c^*]} \min\braces{1-(C^*)^{-1}, (c^*)^{-1}-1}
        = \frac{\gamma^*-1}{\gamma^*+1} \:.
    $$
    Substituting this optimal $\delta^*$ into the derived bound, we obtain
    $$
        \mathbb{P}\big[\gamma(\hat{G}^{(k)}) > \gamma^*\big]
        \le \frac{\alpha_{\mathrm{tmp},1}^2}{kN} \Big(\frac{\gamma^*+1}{\gamma^*-1}\Big)^2
        \:.
    $$
    To guarantee $\gamma(\hat{G}^{(k)}) \le \gamma^*$ with a probability larger than $p\in(0,1)$, we set 
%    $$
%        %\mathbb{P}\bracs*{\gamma(\hat{G}^{(k)}) > \gamma^*}
%        %\le 1-p
%        %\overset{!}{\le}
%        \frac{\alpha_{\mathrm{tmp},1}^2}{kN} \Big(\frac{\gamma^*+1}{\gamma^*-1}\Big)^2
%        =  1-p .
%    $$
%    Using the assumption $\gamma^*>1$, we can solve the final equality for $\gamma^*$ and obtain
    \begin{align}
        \frac{\alpha_{\mathrm{tmp},1}^2}{kN} \Big(\frac{\gamma^*+1}{\gamma^*-1}\Big)^2 = 1-p
        \quad&\Leftrightarrow\quad
        \frac{\gamma^*+1}{\gamma^*-1} = \sqrt{(1-p)\frac{kN}{\alpha_{\mathrm{tmp},1}^2}} =: \alpha_{\mathrm{tmp},2} \\
        \quad&\Leftrightarrow\quad
        \gamma^* + 1 = \alpha_{\mathrm{tmp},2} (\gamma^* - 1) \\
        \quad&\Leftrightarrow\quad
        \gamma^* = \frac{\alpha_{\mathrm{tmp},2}+1}{\alpha_{\mathrm{tmp}.2} - 1} .
    \end{align}
    Obtaining a solution $\gamma^*>1$ requires $\alpha_{\mathrm{tmp},2} > 1$, which follows from the assumption $\sqrt{kN} > \alpha := \alpha_{\mathrm{tmp},1}(1-p)^{-1/2}$.
    Resubstituting $\alpha_{\mathrm{tmp},2}$ yields the claim since
    \begin{align}
        \gamma^*
        &=  \frac{\alpha_{\mathrm{tmp},2}+1}{\alpha_{\mathrm{tmp},2} - 1}
        = \frac{\sqrt{\pars{1-p} \frac{kN}{\alpha_{\mathrm{tmp},1}^2}}+1}{\sqrt{\pars{1-p} \frac{kN}{\alpha_{\mathrm{tmp},1}^2}} - 1}
        = \frac{\sqrt{kN} + \alpha_{\mathrm{tmp},1}(1-p)^{-1/2}}{\sqrt{kN} - \alpha_{\mathrm{tmp},1}(1-p)^{-1/2}}
        = \frac{\sqrt{kN} + \alpha}{\sqrt{kN} - \alpha} . \qedhere
    \end{align}
\end{proof}

% We want to have the matrix framing
% (1-δ)A <= Y <= (1+δ)A    <—>    |Y-A| <= δA.
% We can apply the Matrix-Chebyshev inequality
% P[ |Y-A| > δA ] <= Tr[ E[(Y-A)^2] A^{-2} ] / δ^2
% We know that
% E[ (Y - A).^2 ] <= c
% i.e. we have point-wise bounds for the elements of the matrix
% Since Y and A are symmetric,
% Tr[ (Y-A)^2 ] = |Y-A|_F^2 = \sum_{jk} (Y-A)_{jk}^2
% So
% E[Tr[ (Y-A)^2 ]] <= d^2 c.
% Using Tr[AB] <= Tr[A] lmax(B), we can thus bound
% Tr[E[(Y-A)^2] A^{-2}] <= d^2c / lmin(A)^2 = kc
% for k = (d/lmin(A))^2.
% This gives us a rate of convergence of
% P[ |Y-A| > δA ] <= k c / δ^2
% However, we want this to be satisfied with a constant probability p for different sample sizes n and the resulting c=c(n).
% Hence, we need to solve the subsequent equation for δ:
% p <= kc(n) / δ^2    <—>    δ <= sqrt(kc(n)/p)    —>    δ = sqrt(kc(n)/p)
% where the final implication follows because we want to find the largest δ.
% Finally, recall that c=c(k,n)=γ/N, with the total number of samples N=kn where k is the step in the iterative scheme and n is the sample size in each step.
% Hence
% δ = sqrt(kγ/p) / sqrt(N)
% i.e. the rate is
% δ = O(N^{-1/2})
% Note that we actually don’t care for δ but for the ratio (1+δ)/(1-δ).
\section{Experiments}
\label{sec:experiments_intro}

This section illustrates the presented approach on four numerical examples.
As discussed in section~\ref{sec:setting}, we assume that sample points can be drawn exactly, restricting our experiments to mostly one-dimensional problems:
\begin{itemize}
    \item a dictionary of Hermite polynomials and the Gaussian measure on $\mbb{R}$,
    \item a dictionary of random polynomials and the Gaussian measure on $\mbb{R}$
    \item a dictionary of step functions and the uniform measure on $[0, 1]$, and 
    \item a dictionary of monomials and a singular measure on a two-dimensional domain.
\end{itemize}

\begin{figure}[htb]
    \centering
    \includegraphics[width=\textwidth]{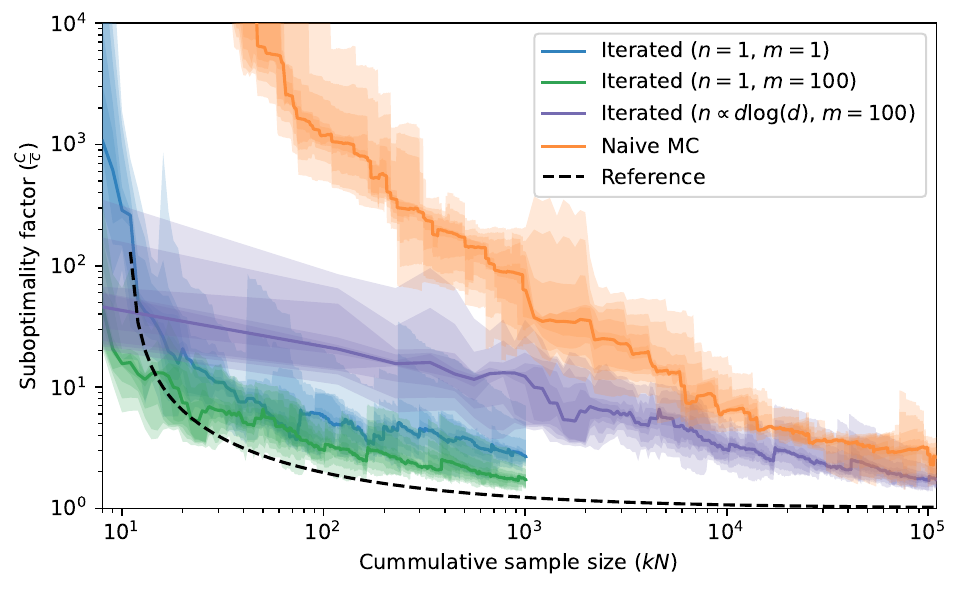}
    \caption{Experiment from section~\ref{sec:experiment_hermite}. Equispaced quantiles for the suboptimality factor $\tfrac{C}{c}$ from equation~\eqref{eq:suboptimality_factor} plotted against the cumulative sample size $kN$.}
    \label{fig:unbounded_domain_1d_hermite}
\end{figure}

\begin{figure}[htb]
    \centering
    \includegraphics[width=\textwidth]{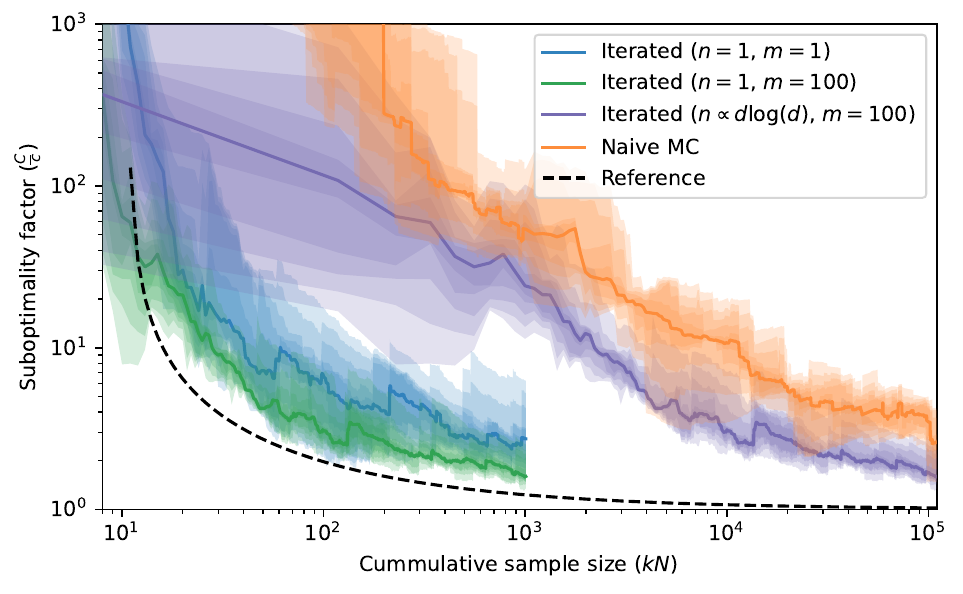}
    \caption{Experiment from section~\ref{sec:experiment_redundant_monomial}. Equispaced quantiles for the suboptimality factor $\tfrac{C}{c}$ from equation~\eqref{eq:suboptimality_factor} plotted against the cumulative sample size $kN$.}
    \label{fig:unbounded_domain_1d_random}
\end{figure}

\begin{figure}[htb]
    \centering
    \includegraphics[width=\textwidth]{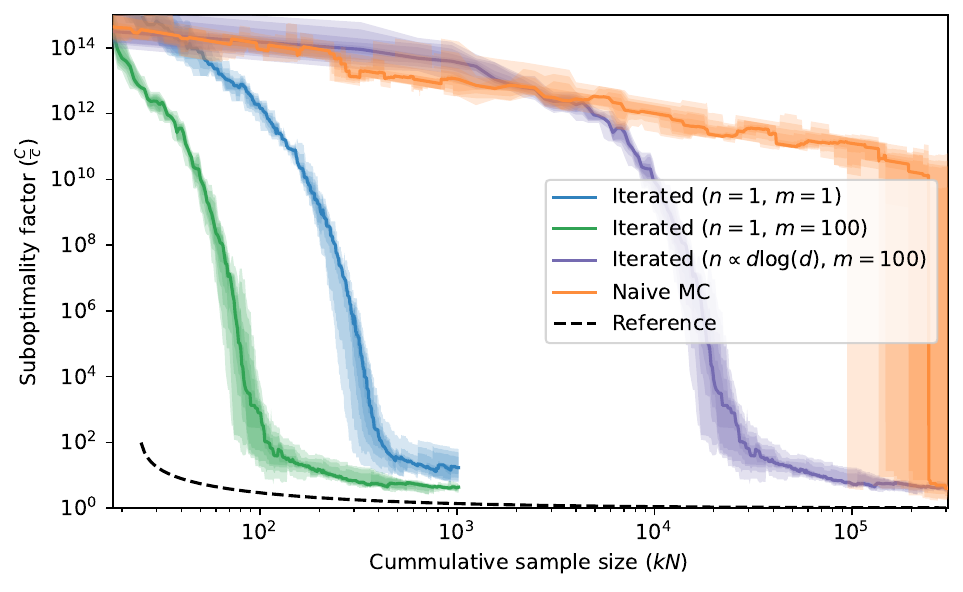}
    \caption{Experiment from section~\ref{ex:pw_const}. Equispaced quantiles for the suboptimality factor $\tfrac{C}{c}$ from equation~\eqref{eq:suboptimality_factor} plotted against the cumulative sample size $kN$.}
    \label{fig:compact_domain_1d_step}
\end{figure}

To illustrate the convergence of the Gramian matrix $\hat{G}^{(k)}$, we plot the suboptimality factor $\tfrac{C}{c}$ from equation~\eqref{eq:suboptimality_factor} against the required computational time, measured by the cumulative sample size $kN$ for $k\le1000$ iterations.
We compare the resulting values to those obtained when estimating the Gramian matrix naively via the Monte Carlo estimation
\begin{equation}
\label{eq:G_naive_MC}
\begin{aligned}
    \hat{G}^{(\mathrm{MC}, n)}
    := \frac1n \sum_{i=1}^n B(y^{(\mathrm{MC})}_{i}) B(y^{(\mathrm{MC})}_{i})^\intercal,
    \qquad
    y^{(\mathrm{MC})}_{1},\ldots,y^{(\mathrm{MC})}_{n} \sim \rho .
\end{aligned}
\end{equation}
We repeat this experiment $10$ times and compute the $10\%\,$--$\,90\%$ percentiles for the resulting $10$ suboptimality factors.
These quantiles are displayed in Figures~\ref{fig:unbounded_domain_1d_hermite}--\ref{fig:2d_approximation}.
The reference lines correspond to the bound~\eqref{eq:gamma_bound} with $\alpha=\tfrac{\sqrt{d-1}}{\sqrt{1-p}}$ for a fixed probability $p=\tfrac34$.\footnote{We make this choice, reasoning that the bound should be infinite for $kN < d$ but finite when $kN \gtrsim d$ (with a multiplicative factor depending on $p$).
While the general behaviour of the bound seems to be confirmed by the experiments, we can see that the lower bound from Theorem~\ref{thm:high_probability} is clearly too optimistic in general.}
Further details about the experiments are provided in the subsequent three paragraphs.

All experiments demonstrate the advantage of the proposed initialisation strategy and clearly display regions of (super-)exponential decay of the suboptimality factor.
While this is not ensured by the variance bounds in Theorem~\ref{thm:mean_and_var}, it is not surprising given the motivation of the proposed algorithm.
Indeed, this behaviour is predicted in a simple setting in Appendix~\ref{app:update}.

\begin{figure}[htb]
    \centering
    \includegraphics[width=\textwidth]{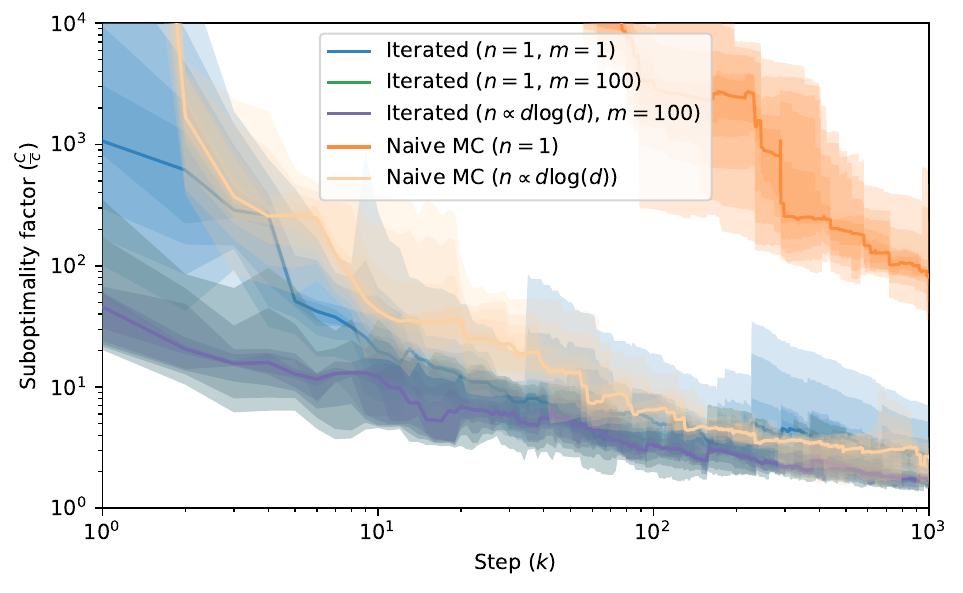}
    \caption{Experiment from section~\ref{sec:experiment_hermite}. Equispaced quantiles for the suboptimality factor $\tfrac{C}{c}$ from equation~\eqref{eq:suboptimality_factor} plotted against the number of steps $k$.}
    \label{fig:unbounded_domain_1d_hermite_step}
\end{figure}

As predicted in Theorem~\ref{thm:mean_and_var}, the experiments indicate that a larger sample size $\hat{n}^{(k)}$ has a positive influence on the convergence.
Increasing the sample size from $n=1$ to $n=4d\ln(4d)$ can reduce the suboptimality factor $\frac{C}c$ at a given iteration $k$ by up to two orders of magnitude (cf.~Figure~\ref{fig:unbounded_domain_1d_hermite_step}).
Note, however, that the cumulative sample size $kN$ is more than two orders of magnitudes larger, indicating that this advantage vanishes when considering the cumulative step size $kN$ instead of the step $k$ (cf.~Figure~\ref{fig:unbounded_domain_1d_hermite}).
Finally, note that although a larger sample size $n$ reduces the suboptimality factor, an acceptable ratio of $\frac{C}c\approx 10$ can already be achieved in all  experiments with the minimal choice $n=1$.

A similar behaviour can be observed for the parameter $\hat m$, where Theorem~\ref{thm:mean_and_var} also predicts a decreasing variance of the estimate with increasing $m$.
Notably, however, a large problem-adapted choice of $\hat m=V$ does not seem to be necessary and the algorithm already converges well for a sample size of $m=100$.

\subsection{Unbounded domain}
\label{sec:experiment_hermite}

The first experiment demonstrates the good behavior of the proposed algorithm on a standard but important reference problem.
We take $\rho$ as the Gaussian measure on $\mbb{R}$ and define the $L^2(\rho)$-orthonormal basis of normalised Hermite polynomials $B_j(x) = H_j(x)$ for $j=0,\ldots,7$.
These polynomials are commonly used in polynomial chaos expansions in finance and  uncertainty quantification.
%~\cite{farchmin2021,Bayer2023}. \comment{[A: PC expansions have been used in thousands of works for about 30 years, i would either cite a fundamental paper of Wiener or simply skip citations.]}
% Looking at Figure~\ref{fig:trivial_1d_unbounded_domain_hermite}, we see that it is extremely difficult to estimate the Gramian matrix $G$, even though $G=I$ is very well-conditioned.

\subsection{Unbounded domain with redundancies}
\label{sec:experiment_redundant_monomial}

The second experiment tests the robustness of our algorithm with respect to redundancies in the dictionary.
%This models the behaviour of tangent vectors to neural networks that are computed in natural gradient descent (cf.~Section~\ref{sec:introduction}). 
We take $\rho$ again as the Gaussian measure on $\mbb{R}$ and define the basis of monomials $b_j(x) = x^j$ for $j=0,\ldots,7$.
For $K=50$, we then draw a random matrix $S\in \mbb{R}^{Kd\times d}$, $d=8$, with independent Gaussian entries\footnote{Observe that $\operatorname{rank}(S) = d$ almost surely.} and define the vector of dictionary functions
$$
    B := S b .
$$
This random dictionary is selected once and remains fixed during all repetitions of the experiment.
% The proposed algorithm still performs quite well in Figure~\ref{fig:trivial_1d_unbounded_domain_random} even though the Gramian matrix $G$ is ill-conditioned.

\begin{remark}
    Although the first two experiments are performed on an unbounded domain, they demonstrate the advantage of the proposed method for compact domains.
    Indeed, when the basis does not have to be polynomial we can consider the quantile function of the standard Gaussian measure $q(x) := \sqrt{2} \operatorname{erf}^{-1}(2x - 1)$ and define the dictionary $B_j(x) = q(x)^j$.
    This dictionary produces the same (empirical) Gramians with respect to the uniform measure on $[0, 1]$ as the dictionary of monomials $B_j(x) = x^j$ with respect to the Gaussian measure.
\end{remark}

\subsection{Compact domain}
\label{ex:pw_const}

In the third experiment, we take $\rho$ as the uniform measure on $[0,1]$ and define a $d=18$ dimensional dictionary of step functions.
To do this, we define the tuple $\beta := \pars{0, 2^{-17}, 2^{-16}, \ldots, 2^0}\in[0,1]^{18}$ and the corresponding partitioning of the interval $[0,1]$ by $I_j := [\beta_j, \beta_{j+1}]$ for $j=1,\ldots,18$.
This partitioning gives rise to the $L^2(\rho)$-orthogonal basis
$$
    B_j := \chi_{I_j} ,
    \qquad j=1,\ldots,18 .
$$
This type of piece-wise constant basis models the result of an adaptive refinement of a mesh for approximating a function with a singularity at $0$.
%PDE with singularities.
% {\color{red}
% Although the (empirical) Gramians in this experiment are technically ill-conditioned, they are diagonal and the basis functions are all bounded in $\braces{0, 1}$.
% While the variance of the Gramian in the first two experiments is high, since the basis functions are unbounded, the variance in the final experiment is naturally bounded.
% }

\subsection{Singular measure}
\label{sec:christoffel-darboux}

This experiment illustrates the usefullness of the presented approach in an application to semi-algebraic approximation. 
Let $f:[0,1]\to[0,1]$ be a bounded measurable function and consider the associated singular measure
$$
    \dx[\rho(x,y)] := \delta_{f(x)}(y) \dx ,
$$
supported on the graph of $f$.
The moment matrix of degree $2d$ for this measure is given by the Gramian matrix $G$ with respect to the dictionary $B_{jk}(x, y) := x^j y^k$ with $j,k=0, \ldots d-1$.
Given this matrix,~\cite{Marx2021} propose to compute the approximation
$$
    f_d(x) := \argmin_{y\in[0,1]} \mfrak{K}_G(x, y) .
$$
They provide convergence guarantees for the approximation of semi-algebraic functions in the limit $d\to +\infty$.
These error bounds rely on the uniform framing~\eqref{eq:framing} (cf.~Assumption~2 and equation~(12) in~\cite{Marx2021}), and the algorithm proposed in the present paper can be used to ensure this framing with an acceptable sample size.

\begin{figure}
    \centering
    \includegraphics[width=\linewidth]{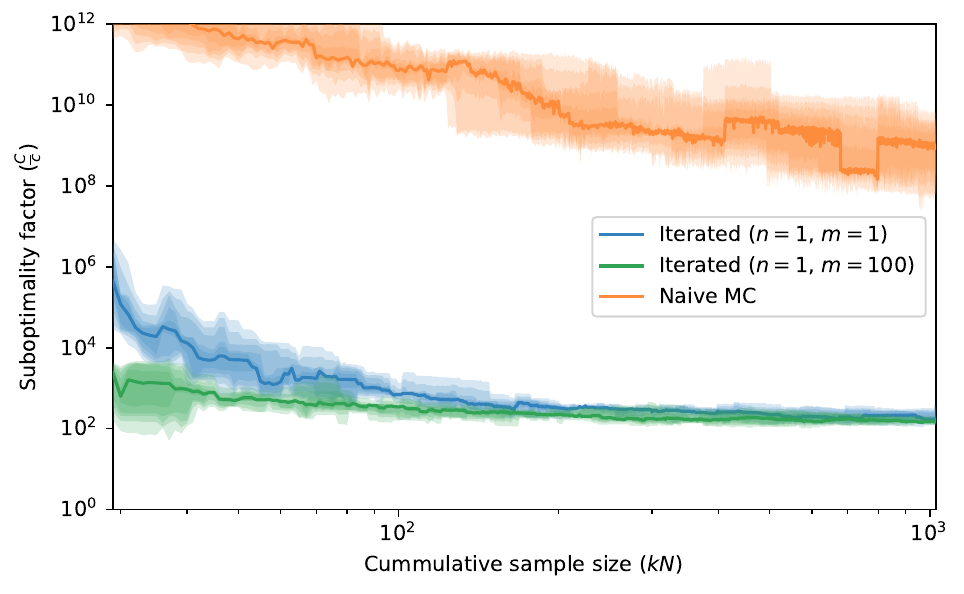}
    \caption{Experiment from section~\ref{sec:christoffel-darboux}. Equispaced quantiles for the suboptimality factor $\tfrac{C}{c}$ from equation~\eqref{eq:suboptimality_factor} plotted against the cumulative sample size $kN$.}
    % \caption{Christoffel--Darboux approximation of $f$.
    % The level sets of $\mfrak{K}_H$ ($H=\hat{G}^{(0)}$ left, $H=\hat{G}^{(10)}$ middle, $H=G$ right) are displayed in grayscale.
    % The target function $f$ is displayed in red, and the corresponding approximation $f_d$ is displayed in blue.}
    \label{fig:2d_approximation}
\end{figure}

To give a numerical illustration, let $q$ be the quantile function for the standard Gaussian distribution and define the function $f = f_{10^{-3}}$ with $f_\varepsilon:[0,1]\to[0,1]$ given by 
$$
    f_\varepsilon(x) = \frac{q(\varepsilon) - q((1 - 2\varepsilon) x + \varepsilon)}{2q(\varepsilon)} .
$$
The convergence of the suboptimality ratio $\frac{C}{c}$ is presented in Figure~\ref{fig:2d_approximation}.\footnote{
We attribute the stagnation at $\frac{C}c \approx 10^2$ to a problem of numerical stability when computing the reference density.
The $n \propto d\log(d)$-line is not plotted for the same reason.}

\section{Discussion}
\label{sec:discussion}

% {\color{red}
% \begin{itemize}
%     \item We propose an iterative method to initialise and further refine such an estimate.
%     \item In particular, we show in Propositions~\ref{prop:vol_sampling} and~\ref{prop:increasing_rank} that the proposed initialisation procedure is able to find an estimate with the same image and kernel as the sought Gramian.
%     \item We then motivate and analyse the proposed iterative refinement strategy.
%     \item The resulting theoretical bounds of Theorem~\ref{thm:mean_and_var} ensure convergence.
%     \item While these bounds are rather crude, they ensure that the proposed algorithm converges.
%     \item The performed numerical experiments demonstrate that it works remarkably better than the convergence rates would make us believe.
%     This is discussed in Appendix~\ref{app:update}.
% \end{itemize}
% }

\revision[0]{
The present paper addresses the problem of estimating the Gramian matrix associated with a potentially overcomplete dictionary.
This task is crucial in least squares approximation and the many algorithms built upon it.
To tackle it, we first introduce and analyse an iterative procedure   for constructing an initial estimate.
Propositions~\ref{prop:vol_sampling} and~\ref{prop:increasing_rank} show that this initialisation finds an estimate with the correct image and kernel.
We then motivate and analyse a refinement strategy and establish its convergence in Theorem~\ref{thm:mean_and_var}.
Given the crudeness of these bounds, it is reasonable to expect stronger practical performance, which our numerical experiments indeed confirm.
(This phenomenon is also discussed in Appendix~\ref{app:update}.)}

\paragraph{Sampling.}

Given that $\mu_{k}$ is only known up to the normalisation factor, it also seems natural to employ \emph{Markov chain Monte Carlo} (MCMC) methods to sample from this measure and estimate the Gramian $G^{(k+1)}$.
This requires an adequate truncation of the chain and post-processing~\cite{South2022}.

\revision[0]{
To ensure that the distribution of the sample points quickly converges to the desired stationary distribution $\mu_k$, a proper \emph{burn-in phase} and initialisation is essential.
Since the Gramian estimates change less with every step, using the last element of the chain in step $k$ as an initial point for the chain in step $k+1$ seems sensible.

Although each individual chain is Markov (conditioned on the corresponding Gramian $\hat{G}^{(k)}$), the concatenation of all chains no longer possess the Markov property.
Establishing a central limit type theorem for this so-called \emph{adapted Markov chain} is challenging, but sufficient conditions were analysed in~\cite{Fort2011}.
Demonstrating these sufficient conditions for a fixed chain length $n^{(k)} = n$, or even a chain length of $n^{(k)} = 1$, would be beneficial to simplify the sampling algorithm.

However, this idea of using adapted Markov chains should be approached with caution, since MCMC are known to converge slowly for highly localised probability densities and in high dimensions (precisely the regimes where the Gramian is most difficult to estimate).
\footnote{\revision[0]{For this reason, we did not consider high-dimensional examples using MCMC in our numerical experiments.
A successful experiment would require a reasonable convergence of the MCMC method, indicating that the Christoffel density we aim to approximate was relatively simple.
Conversely, if an experiment failed, it would be unclear whether the failure was caused by problems with the MCMC sampling or by inaccuracies in the Gramian estimation.}}

It is conceivable that we may use an iterative scheme similar to the one for approximating $G$ for drawing samples that discretise the measure $w^{-1}\rho$.
}

\revision[0]{\paragraph{Practical implementation.}

Computing the inverse Christoffel function $\mfrak{K}_{\hat{G}^{(k)}}$ according to equation~\eqref{eq:gramian} relies on computing the pseudo-inverse of $\hat{G}^{(k)}$.
Since $G$ may be extremely ill-conditioned some of its eigenvalues may be numerically zero.
This issue is only aggravated by $\hat{G}^{(k)}$ being an empirical estimate, where some eigenvalues may be severely underestimated
(Section~\ref{sec:experiments_intro} presents examples where this happens).
It is not clear how this influences numerical stability.

Moreover, there currently exists no simple a posteriori method to check whether~\eqref{eq:framing} holds for a given sample.
}

\paragraph{Analysis.}

\revision[0]{Many of the difficulties in proving convergences arise form the necessity to estimate the normalisation constant $z^{(k)}$ by $\hat{z}^{(k)}$.
This, in fact, is the only reason why we can not make the motivational argument from  the introduction rigorous. 
(See Appendix~\ref{app:update} for a deeper discussion.)}

Also note that the proposed algorithm is equivalent to empirically orthogonalising the dictionary $B$ with respect to $\hat{G}^{(k)}$ before computing the Gramian estimate $\hat{G}^{(k+\sfrac12)}$.
The formula for $\hat{G}^{(k+1)}$ in this orthonormal basis is then given by
$$
    \hat G^{(k+1)} = \tfrac{k}{k+1} I + \tfrac{1}{k+1} \hat{G}^{(k+1/2)} .
$$
Although the fundamental algorithm remains unchanged, analysing it in this way may bring some advantages because the Gramian estimates converge to identity.
\revision[0]{This may also help with numerical stability.}

\revision[0]{\paragraph{Relation to least squares.} Finally,} we recognise that sampling from the Christoffel function is seldom the ultimate goal.
Usually, the goal is to obtain a least squares approximation with a prescribed $L^2(\rho)$-error.
To this end, it may be more effective to perform a least squares approximation for a proper subspace of $\mcal{V} = \operatorname{span}(B)$.
An intuition for this is given by the bound
$$
	\norm{\mfrak{K}}_{L^\infty(\rho)} \le \sum_{j=1}^d \norm{b_j}_{L^\infty(\rho)}^2 ,
$$
which is true for any orthonormal basis $b:\mcal{X}\to\mbb{R}^d$ of $\mcal{V}$.
Since basis functions with large $L^\infty(\rho)$-norm have a larger oscillation ($\operatorname{osc}(b_j) := \sup_{\mcal{X}} b_j(x) - \inf_{\mcal{X}} b_j(x) \le 2 \norm{b_j}_{L^\infty(\rho)}$), they are generally less relevant for approximating functions in classical smoothness classes.
Consider, for example, approximating a smooth function using the piece-wise constant basis from the third experiment in section~\ref{sec:experiments_intro}.
In this setting, it is reasonable to disregard basis functions with smaller support since their corresponding inner products are more complicated to estimate, while their impact on the approximation error is less significant.
One way to do this is by using an algorithm based on three independent samples: the first for empirical orthogonalisation, the second for selecting a stable subspace, and the final sample for the actual approximation.
However, selecting a subspace that ensures stability with high probability while simultaneously bounding the ensuing approximation error for all functions in a particular regularity class remains an open problem.
\revision[0]{A similar argument has been put forward and rigorously analysed in~\cite{Adcock2019_frames}.}

\section*{Acknowledgements}
This project is funded by the ANR-DFG project COFNET (ANR-21-CE46-0015).
This work was partially conducted within the France 2030 framework programme, Centre Henri Lebesgue ANR-11-LABX-0020-01.

Our code makes extensive use of the Python packages \texttt{numpy}~\cite{numpy}, \texttt{scipy}~\cite{scipy}, and \texttt{matplotlib}~\cite{matplotlib}.

{
    \emergencystretch=3em
    \printbibliography
}

\appendix

\section{Discussion on the inverse Christoffel function}
\label{app:christoffel}

Let $B:\mcal{X}\to\mbb{R}^D$ be a given dictionary of functions.
This section demonstrates that the definition of the inverse Christoffel function
\begin{equation*}
%\label{eq:gramian:recall}
    \mfrak{K}_G\pars{x}
    := \inner{G^+, B(x)B(x)^\intercal}_{\mathrm{Fro}}
    \qquad\text{with}\qquad
    G_{kl} := (B_k, B_l)_{L^2(\rho)} %\tag{\ref*{eq:gramian}}
\end{equation*}%
%\noeqref{eq:gramian:recall}%
from equation~\eqref{eq:gramian} in section~\ref{sec:setting} is equivalent $\rho$-almost everywhere  to the customary definition 
\begin{equation}
\label{eq:sup_definition}
    \mfrak{K}_G\pars{x} = \sup_{v\in\mcal{V}} \frac{\abs{v(x)}^2}{\norm{v}_{L^2(\rho)}^2} .
\end{equation}

\begin{proposition}
\label{prop:inverse_christoffel_definition}
    Let $H\in\mbb{R}^{D\times D}$ be a positive semi-definite matrix and $\norm{v}_H^2 := v^\intercal H v$ the corresponding \mbox{(semi-)norm} for $v\in\mbb{R}^D$.
    Then, if $\ker\pars{H} \subseteq \ker\pars{G}$,
    $$
        \mfrak{K}_H(x)
        = \sup_{v\in\mbb{R}^D } \frac{\abs{B(x)^\intercal v}^2}{\norm{v}_{H}^2} 
    $$
    $\rho$-almost everywhere.
    Since every function $v\in\mcal{V} = \operatorname{span}(B)$ can be expanded as $v(x) = B(x)^\intercal c_v$ with some $c_v\in\mbb{R}^D$, the inclusion $\ker(H)\supseteq\ker(G)$ implies that $\norm{v}_H := \norm{c_v}_H$ is a well-defined \mbox{(semi-)norm} on $\mcal{V}$ and that $\norm{v}_G = \norm{v}_{L^2(\rho)}$ for all $v \in \mcal{V}$.
    This entails the identity~\eqref{eq:sup_definition} $\rho$-almost everywhere.
    %\begin{equation}
    %\label{eq:sup_definition:recall}
    %    \mfrak{K}_G(x)
    %    = \sup_{v\in\mcal{V}} \frac{\abs{v(x)}^2}{\norm{v}_{L^2(\rho)}^2}
    %     \tag{\ref*{eq:sup_definition}}
    %\end{equation}%
    %\noeqref{eq:sup_definition:recall}%
    Moreover, suppose that $b:\mcal{X}\to\mbb{R}^d$ is a vector of $L^2(\rho)$-orthonormal basis functions for $\mcal{V}$.
    Then
    $$
        \mfrak{K}_G(x)
        = \sum_{j=1}^d b_j(x)^2 .
    $$
\end{proposition}
Note that the condition $\ker(H)\subseteq\ker(G)$ is necessary because $\mfrak{K}_G$, as defined in equation~\eqref{eq:gramian} is always finite, while the inverse Christoffel function in definition~\eqref{eq:sup_definition} can be infinite.
Let, e.g.\ $B$ be a basis of step functions with $G = I$ and $H = I - e_1e_1^\intercal$.
Then the supremum is infinite on $\operatorname{supp}(B_1)$.

\begin{proof}[Proof of Proposition~\ref{prop:inverse_christoffel_definition}]
    First, assume that $\ker\pars{H} \subseteq \ker\pars{G}$ and let $H = U\Lambda^2 U^\intercal$ be a rank-revealing spectral decomposition of $H$ where $U\in\mbb{R}^{D\times r}$ and $\Lambda\in\mbb{R}^{r\times r}$ with $\operatorname{rank}(\Lambda) = r$ for some $r\ge d$.
    Since $I - UU^\intercal$ is an orthogonal projection onto $\ker(H) \subseteq \ker(G)$, it follows for all $v\in\mathbb{R}^D$ that
    $$
        \|B^\intercal v - B^\intercal UU^\intercal v\|_{L^2(\rho)}^2
        = v^\intercal (I - UU^\intercal) G (I - UU^\intercal) v
        = 0 .
    $$
    From this, we can deduce that for $\rho$-almost all $x\in \mcal{X}$
    \begin{align}
        \mfrak{K}_{H}\pars{x}
    	&= \inner{H^+, B\pars{x}B\pars{x}^\intercal}_{\mathrm{Fro}} \\
    	&= B\pars{x}^\intercal H^+B\pars{x} \\
    	&= (\Lambda^{-1} U^\intercal B\pars{x})^\intercal (\Lambda^{-1} U^\intercal B\pars{x}) \\
        &= \sup_{w\in\mbb{R}^r} \frac{\abs{(\Lambda^{-1} U^\intercal B\pars{x})^\intercal w}^2}{w^\intercal w} \label{eq:step:onb} \\
        &= \sup_{w\in\mbb{R}^r} \frac{\abs{B\pars{x}^\intercal Uw}^2}{w^\intercal \Lambda^2 w} \\
        &= \sup_{v\in\mbb{R}^D} \frac{\abs{B\pars{x}^\intercal UU^\intercal v}^2}{(U^\intercal v)^\intercal \Lambda^2 (U^\intercal v)} \\
    	&= \sup_{v\in\mbb{R}^D} \frac{\abs{B\pars{x}^\intercal v}^2}{v^\intercal H v} .
        \label{eq:step:w->v} 
    \end{align}
    To see the second assertion, let $H=G$ and observe that
    $$
        \mfrak{K}_{G}\pars{x}
        = \sup_{v\in\mbb{R}^D} \frac{\abs{B(x)^\intercal v}^2}{v^\intercal G v}
        = \sup_{v\in\mathcal{V}} \frac{\abs{v(x)}^2}{\|v\|_{L^2(\rho)}^2} .
    $$
    The supremum in this equation can be computed explicitly, yielding the final assertion.
\end{proof}

\section{On the optimality of ``optimal sampling''}
\label{app:dlogd}

\revision[0]{The next theorem shows that the optimal sampling distribution is indeed optimal in the case of i.i.d.~samples. This result is a generalisation of Proposition~5 in~\cite{Trunschke2024}). 
 This theorem also shows that the factor $\|w \mathfrak{K}_G\|_{L^\infty(\rho)}$ appears not only in the sufficient sample size bound~\eqref{eq:sample_size_bound} but also in the necessary sample size bound~\eqref{eq:sample_size_bound_necessary}.}
 
\begin{revisione}[0]
{
\begin{minipage}{\textwidth}
\begin{theorem}\label{thm:dlogd}
    \revision[0]{
    Let $T_1,\ldots,T_d \subseteq \mathcal{X}$ be a non-trivial partition of $\mathcal{X}$ by sets of positive measures.
    For all $k=1,\ldots,d$ let $B_k = \chi_{T_k}$ be the characteristic function of $T_k \subseteq \mathcal{X}$.
    Assume that $n$ independent points $x_1,\ldots, x_n$ are drawn according to $\nu = w^{-1}\rho$, with $w$ piecewise constant on the above defined partition, and define the random event
    $$
        S_n  :\Leftrightarrow \forall \,v\in\mathcal{V}\setminus \{0\} : \|v\|_{n}^2>0
        % \ :\Leftrightarrow\ 
        % \exists \,\delta\in(0,1)\,
        % \forall \,v\in\mathcal{V}:
        %        (1-\delta)\norm{v}_{L^2(\rho)}^2
        %        \le \norm{v}_{n}^2 .
    $$
    as well as the stopping time $N$ as the smallest number of samples $n$ required to satisfy $S_n$.
    Then
    \begin{equation}
    \label{eq:sample_size_bound_necessary}
        n > \tfrac{1}{2}\|w\mathfrak{K}_G\|_{L^\infty(\rho)}\ln\big(\tfrac{1}{p}\big)
    \end{equation}
    sample points are necessary to ensure $\mathbb{P}[S_n] \ge 1-p$ and
    \begin{equation}
    \label{eq:sample_size_bound_sufficient}
        n \ge \|w\mathfrak{K}_G\|_{L^\infty(\rho)}\ln\big(\tfrac{d}{p}\big)
    \end{equation}
    sample points are sufficient to ensure $\mathbb{P}[S_n] \ge 1-p$.
    Moreover, $\mathbb{E}[N]$ is minimised for the weight function $w = d\mathfrak{K}_G^{-1}$ for which (in Landau notation)
    $$
        \mathbb{E}[N] = d \ln(d) + \Theta(d) .
    $$
    Moreover, it holds for all $c\in\mathbb{R}$ that}
    $$
        \revision[0]{\lim_{d\to\infty} \mathbb{P}[S_{d\ln(d) + cd}] = \exp(-\exp(-c)) .}
    $$
\end{theorem}
\end{minipage}
%\addtocounter{theorem}{-1}
}
\begin{proof}
    Observe that $S_n$ is equivalent to $\lambda_{\mathrm{min}}(\hat{G})>0$ where $\hat{G}$ is the Gramian matrix with respect to the empirical (semi-)inner product
    $$
        \hat{G}_{kl} = (B_k, B_l)_n := \frac1n \sum_{i=1}^n w(x_i) B_k(x_i) B_l(x_i) .
    $$
    Since $\hat{G} = \operatorname{diag}(\|B_1\|_n^2, \ldots, \|B_d\|_n^2)$, we obtain that
    \begin{align}
        S_n
        \ &\Leftrightarrow\ 
        \forall \,k=1,\ldots,d: 
                \|B_k\|_n^2 > 0 \\
        \ &\Leftrightarrow\ 
        \forall \,k=1,\ldots,d\,\,
        \exists \,i \in \{1,\ldots,n\} :
                x_i\in T_k .
    \end{align}
    Defining $p_k := \mathbb{P}[x_1\in T_k] = \int_{T_k} w^{-1}\dx[\rho]$, and using Fr\'echet's inequality, we obtain that
    \begin{align}
        \mathbb{P}[S_n]
        &\le \min_{k=1, \ldots, d} \mathbb{P}[\exists \,i \in \{1,\ldots,n\}: x_i\in T_k] \\
        &= 1 - \max_{k=1, \ldots, d} \mathbb{P}[\forall \,i \in \{1,\ldots,n\} : x_i\not\in T_k] \\
        &= 1 - \max _{k=1, \ldots, d} (1- p_k)^n
    \intertext{and}
        \mathbb{P}[S_n]
        &\ge 1 - \sum_{k=1}^d \mathbb{P}[\forall \,i \in \{1,\ldots,n\}: x_i\not\in T_k] \\
        &\ge 1 - d \max_{k=1,\ldots,d} \mathbb{P}[\forall \,i \in \{1,\ldots,n\}: x_i\not\in T_k] \\
        &\ge 1 - d \max_{k=1,\ldots,d} (1 - p_k)^n .
    \end{align}
    Since $w$ is piecewise constant on the partition $T_1,\ldots,T_d$, we can write $w = \sum_{k=1}^d w_k B_k$ and
    $$
        p_k^{-1}
        = w_k \rho(T_k)^{-1}
        = \bigg(w \frac{B_k^2}{\|B_k\|_{L^2(\rho)}^2}\bigg)\Bigg|_{T_k}
        % = (w \mathfrak{K}_G)\big|_{T_k}
        = \|w\mathfrak{K}_G\|_{L^{\infty}(\rho|_{T_k})} .
    $$
    This implies, $\max_{k=1,\ldots,d} (1 - p_k)^n
        = \big(1 - \big(\max_{k} p_k^{-1}\big){\vbox to 0.65em {}}^{-1}\big)^n
        = \big(1 - \|w\mathfrak{K}_G\|_{L^\infty(\rho)}^{-1}\big)^n$ and, therefore,
    $$
        1 - d \big(1 - \|w\mathfrak{K}_G\|_{L^\infty(\rho)}^{-1}\big)^n
        \le \mathbb{P}[S_n]
        \le 1 - \big(1 - \|w\mathfrak{K}_G\|_{L^\infty(\rho)}^{-1}\big)^n .
    $$
    Since $\|w\mathfrak{K}_G\|_{L^\infty(\rho)} \ge \|\mathfrak{K}_G\|_{L^1(\rho)} = d \ge 2$ (the partition is non-trivial) we can further simplify
    $$
        1 - d \exp\big(-\tfrac{n}{\|w\mathfrak{K}_G\|_{L^\infty(\rho)}}\big)
        \le 1 - d \big(1 - \|w\mathfrak{K}_G\|_{L^\infty(\rho)}^{-1}\big)^n
        \le \mathbb{P}[S_n]
        \le 1 - \big(1 - \|w\mathfrak{K}_G\|_{L^\infty(\rho)}^{-1}\big)^n
        \le 1 - \exp\big(-\tfrac{2n}{\|w\mathfrak{K}_G\|_{L^\infty(\rho)}}\big) .
    $$
    Solving these inequalities for $n$ yields the first two claims
    \begin{align}
        n \ge \hphantom{\tfrac{1}{2}} \|w\mathfrak{K}_G\|_{L^\infty(\rho)} \ln(\tfrac{d}{p})
        \quad&\Rightarrow\quad \mathbb{P}[S_n] \ge 1 - p,
        \quad\text{and} \\
        n \le \tfrac{1}{2} \|w\mathfrak{K}_G\|_{L^\infty(\rho)} \ln(\tfrac{1}{p})
        \quad&\Rightarrow\quad \mathbb{P}[S_n] \le 1 - p .
    \end{align}
    % \begin{align}
    %     %1 - p
    %     %\le 
    %     \mathbb{P}[S_n]
    %     &\le \min_{k=1,\ldots,d} \mathbb{P}[
    %         \exists \,i=1,\ldots,n:
    %                 x_i\in T_k
    %     ] \\
    %     &= 1 - \max_{k=1,\ldots,d} \mathbb{P}[
    %         \forall \,i=1,\ldots,n:
    %                 x_i\not\in T_k
    %     ] \\
    %     &= 1 - \max_{k=1,\ldots,d} \mathbb{P}[x_i\not\in T_k]^n \\
    %     &= 1 - \max_{k=1,\ldots,d} (1 - \rho(T_k))^n \\
    %     &\le 1 - \max_{k=1,\ldots,d} \exp(-2n\rho(T_k)) \\
    %     &= 1 - \exp\bigg(-\frac{2n}{\max_{k=1,\ldots,d}\rho(T_k)^{-1}}\bigg) \\
    %     &= 1 - \exp\bigg(-\frac{2n}{\|\mathfrak{K}_G\|_{L^\infty(\rho)}}\bigg) .
    % \end{align}
    To show the finer bounds, observe that 
    \begin{align}
        S_n
        % \ &\Leftrightarrow\ 
        % \forall \,k=1,\ldots,d: 
        %         \|B_k\|_n^2 > 0 \\
        \ \Leftrightarrow\ 
        \forall \,k=1,\ldots,d,\,\,
        \exists \,i \in \{1,\ldots,n\} :
                x_i\in T_k 
       \ \Leftrightarrow\ N \le n .
    \end{align}
    and $N \le n$ is an instance of the generalised coupon collector's problem, where the probability to collect the $k$\textsuperscript{th} coupon type is given by $p_k$.
    The expected number of coupons required to collect a full set of $d$ different coupon types is given in~\cite{Flajolet1992} as
    $$
        \mathbb{E}[N] = \int_0^\infty \bigg(1 - \prod_{k=1}^d \big(1 - \exp(-p_k t)\big)\bigg) \dx[t] .
    $$
    Since $\mathbb{E}[N]$ depends on $w$ only via the integrals $p_k = \int_{T_k} w^{-1}(x_i) \dx[\rho(x_i)]$, $k=1,\ldots,d$,  we can assume that the minimiser $w_\star$ is piece-wise constant, i.e.,
    \begin{equation}
    \label{eq:w_star}
        w_\star^{-1} = \sum_{k=1}^d \frac{c_k}{\rho(T_k)} B_k^2 .
    \end{equation}
    Since $p_k = c_k$ for all $k=1,\ldots,d$ we can write the optimisation of $w$ as the constrained minimisation problem
    $$
        \begin{aligned}
        & \underset{p\in\mathbb{R}^d}{\text{minimise}}
        & & f(p) \\
        & \text{subject to}
        & & g_k(p) \le 0,\; k=1,\ldots,d \\
        &&& h(p) = 0
        \end{aligned}
    $$
    with $f(p) = \int_0^\infty \big(1 - \prod_{k=1}^d (1 - \exp(-p_k t))\big) \dx[t]$, $g_k(p) = -p_k$, $k=1,\ldots,d$ and $h(p) = 1 - \sum_{k=1}^d p_k$.
    Since this problem satisfies the \emph{linearity constraint qualification}, the minimiser $p_\star$ has to satisfy the \emph{Karush--Kuhn--Tucker} (KKT) conditions
    $$
        p_\star \in[0,1]^d,
        \quad
        \|p_\star\|_1 = 1,
        % \sum_{k=1}^d p_{\star,k} = 1,
        \quad
        \nabla f(p_\star) = \lambda \boldsymbol{1} + \mu
        \quad\text{and}\quad
        p_\star^\intercal \mu = 0
    $$
    for some $\lambda\in\mathbb{R}$ and $\mu\in[0,\infty)^d$.
    Since $f(p) = \infty$ when $p_k=0$ for any $k=1,\ldots,d$, the optimality conditions reduce to
    $$
        p_\star \in(0,1]^d,
        \quad
        % \sum_{k=1}^d p_{\star,k} = 1
        \|p_\star\|_1 = 1
        \quad\text{and}\quad
        \nabla f(p_\star) = \lambda \boldsymbol{1} ,
    $$
    implying $p_\star = \frac1d \boldsymbol{1}$.
    Inserting $c = p_\star$ into~\eqref{eq:w_star} yields
    $$
        w_\star
        = \bigg(\sum_{k=1}^d \frac{p_k}{\rho(T_k)} B_k^2\bigg)^{-1}
        = d \mathfrak{K}_G^{-1}.
    $$
    This weight function ensures that $p_\star = \frac1d \boldsymbol{1}$, i.e., that all coupon types occur with equal probability.
    % Since $p_{\star,k} = p_{\star,1}$ for all $k=1,\ldots,d$, this weight function 
    This simplifies the problem to the standard coupon collector's problem, for which~\cite[Section~2.4.1]{Mitzenmacher2017-bx}
    $$
        \mathbb{E}[N] = d\ln(d) + \Theta(d) .
    $$
    Moreover, the same book shows in Theorem~5.13 that for all $c\in\mathbb{R}$
    \begin{align}
        \lim_{d\to\infty} \mathbb{P}[S_{d\ln(d) + cd}]
        &= \lim_{d\to\infty} \mathbb{P}[N \le d\ln(d) + cd]
        = \exp(-\exp(-c))
        . \qedhere
    \end{align}
\end{proof}
\end{revisione}
\section{Subsampling can not improve the quasi-optimality constant}
%Proof of Theorem~\ref{thm:weighted_sample_size}
%}
\label{app:weighting}

%{
%\renewcommand{\thetheorem}%{\ref*{thm:weighted_sample_size}}
%\addtocounter{theorem}{-1}
\begin{theorem}\label{thm:weighted_sample_size}
\begin{revisione}[0]
    Let $\mathcal{V}$ be chosen as in Theorem~\ref{thm:dlogd}.
    Suppose that $X_1,\ldots,X_m\sim\rho$ are independent and let $x_1,\ldots,x_n$ be a subsample of $X_1,\ldots,X_m$ (drawn with or without replacement according to any density).
    Then $m\ge \tfrac12 \|\mathfrak{K}_G\|_{L^\infty(\rho)} \ln(\frac{1}p)$ initial points are necessary to ensure that
    $$
        \exists 
        \,c\in(0,\infty) : \|u - u_{\mathcal{V},n}\|_{L^2(\rho)}^2
        \le \|u - u_{\mathcal{V}}\|_{L^2(\rho)}^2 + c \|u - u_{\mathcal{V}}\|_{L^\infty(\rho)}^2
    $$
    %for any $c\ge 1$ 
    with probability $1-p$.
    Moreover, $m\ge 3 \|\mathfrak{K}_G\|_{L^\infty(\rho)} \frac{c^2}{(c - C_w)^2} \ln(\frac{2d}p)$ initial points are sufficient to ensure that
    $$
        \|u - u_{\mathcal{V},n}\|_{L^2(\rho)}^2
        \le \|u - u_{\mathcal{V}}\|_{L^2(\rho)}^2 + c \|u - u_{\mathcal{V}}\|_{L^\infty(\rho)}^2
    $$
    for any $c\ge C_w$ with probability $1-p$.
    The precise value of $C_w \ge \frac{1}{n}$ depends on $w$, $b$ and $X_1, \ldots, X_m$.
    In particular, this bound requires an order of $m \in \Theta(\|\mathfrak{K}_G\|_{L^\infty(\rho)})$ initial points.
\end{revisione}
\end{theorem}
%}

\revision[0]{
To investigate the quality of least squares approximation with a subsample of a given random sample, we begin by considering the slightly more general setting of weighted least squares approximation.
We assume that $x_1, \ldots, x_n\in \mathcal{X}$ are drawn independently from $\rho$ and} define for any $w\in[0,\infty)^n$ the weighted empirical (semi)norm
$$
	\norm{v}_{w}^2 := \sum_{i=1}^n w_i v(x_i)^2 ,
    \quad v\in L^2(\rho).
$$
Analogously to equation~\eqref{eq:empirical_best_approximation}, we define
$$
    u_{\mcal{V},w} := \argmin_{v\in\mcal{V}} \norm{u - v}_w .
$$
\revision[0]{In this setting, standard least squares approximation corresponds to the weight sequence $w = \sfrac{\boldsymbol{1}}n$ and least squares approximation with a subsample corresponds to a sparse weight sequence.
The approximation error can be bounded by the subsequent theorem.
}

\begin{revisione}[0]
\begin{minipage}{\textwidth}
\begin{theorem}
\label{thm:weighted_error_bound}
\begin{revisione}[0]
    Let $b : \mcal{X}\to\mbb{R}^d$ be a vector of $L^2(\rho)$-orthonormal basis functions of $\mathcal{V}$ and define $M\in\mbb{R}^{n\times d}$ by $M_{ij} := b_j(x_i)$ and the (weighted) empirical Gramian matrix
    $$
        G_w := M^\intercal W M
        \qquad\text{with}\qquad
        W := \operatorname{diag}(w) .
    $$
    Then it holds that
    \begin{align}
    	\norm{u - u_{\mcal{V},w}}_{L^2(\rho)}^2
    	&\le \norm{u - u_{\mcal{V}}}_{L^2(\rho)}^2 + \tfrac{1}{\lambda_{\mathrm{min}}(G_{w})}\norm{u - u_{\mcal{V}}}_{w}^2 \\
    	&\le \norm{u - u_{\mcal{V}}}_{L^2(\rho)}^2 + \tfrac{\|w\|_1}{\lambda_{\mathrm{min}}(G_{w})}\norm{u - u_{\mcal{V}}}_{L^{\infty}(\rho)}^2 . \label{eq:suboptimality_Linfty_weighted}
    \end{align}
\end{revisione}
\end{theorem}
\end{minipage}
The proof of this fact follows by the same arguments as the proof of~\eqref{eq:error_bound_n} and~\eqref{eq:error_bound_Linf} in~\cite{cohen_2017_optimal}.
The last inequality is sharp and is attained when $u(x_i) - u_{\mathcal{V}}(x_i) = u(x_1) - u_{\mathcal{V}}(x_1)$ for all $i=1,\ldots,n$.
While this is a worst-case scenario, it is plausible that the constant $\frac{\|v\|_w^2}{\|v\|_{L^\infty(\rho)}^2}$ approaches $\|w\|_1$ when the smoothness of $v$ (measured e.g.\ by a Sobolev norm) increases and $v(x_i)\to v(x_1)$ for all $i=1,\ldots,n$.

Given the error bound~\eqref{eq:suboptimality_Linfty_weighted}, the quality of the weight sequence $w$ over the uniform weight sequence $\sfrac{\boldsymbol{1}}{n}$ can be measured by the fraction
$$
    C_w
    := \frac{\sfrac{\|w\|_1}{\lambda_{\mathrm{min}}(G_{w})}}{\sfrac{\|\sfrac{\boldsymbol{1}}n\|_1}{\lambda_{\mathrm{min}}(G_{\sfrac{\boldsymbol{1}}n})}}
    = \|w\|_1\frac{\lambda_{\mathrm{min}}(G_{\sfrac{\boldsymbol{1}}n})}{\lambda_{\mathrm{min}}(G_w)}.
$$
This fraction can be bounded with the help of the subsequent lemma.
\begin{lemma}
\label{lem:weighted_gramian_smallest_eigenvalue}
    Given any matrix $M\in\mbb{R}^{n\times d}$ it holds for every matrix $W\in\mbb{R}^{n\times n}$ that
    $$
    	\lambda_{\mathrm{min}}(M^\intercal WM) \le \lambda_{\mathrm{max}}(W) \lambda_{\mathrm{min}}(M^\intercal M) .
    $$
\end{lemma}
\begin{proof}
    Let $U\Sigma V^\intercal = M$ be the compact SVD of $M$ and let $r$ denote the rank of $M$ (i.e.\ $U\in\mbb{R}^{n\times r}$).
    Moreover, let $\tilde{W} := U^\intercal W U$ and let $w$ be the eigenvalue corresponding to the smallest eigenvalue $\lambda_{\mathrm{min}}(\Sigma)$ of $\Sigma$.
    Then
    \begin{align}
    	\lambda_{\mathrm{min}}(M^\intercal W M)
    	&= \lambda_{\mathrm{min}}(\Sigma U^\intercal WU\Sigma)
    	= \lambda_{\mathrm{min}}(\Sigma \tilde{W}\Sigma) \\
    	&= \min_{\norm{v}_2=1} v^\intercal \Sigma \tilde{W}\Sigma v
    	\le \lambda_{\mathrm{min}}(\Sigma)^2 w^\intercal \tilde{W} w \\
    	&\le \lambda_{\mathrm{min}}(\Sigma)^2 \lambda_{\mathrm{max}}(W) \norm{w}_2^2
    	= \lambda_{\mathrm{min}}(M^\intercal M) \lambda_{\mathrm{max}}(W)
    	. \qedhere
    \end{align}
\end{proof}

\begin{lemma}
\label{lem:Cw_bound}
    It holds that $C_w \ge \frac1n$.
\end{lemma}
\begin{proof}
    Lemma~\ref{lem:weighted_gramian_smallest_eigenvalue} implies
    $$
    	\lambda_{\mathrm{min}}(G_{w})
    	\le \norm{w}_\infty \lambda_{\mathrm{min}}(G_{\boldsymbol{1}})
    	= n\norm{w}_\infty \lambda_{\mathrm{min}}(G_{\sfrac{\boldsymbol{1}}n})
    $$
    and therefore
    \begin{align}
    	C_w &\ge \frac{\|w\|_1}{n\norm{w}_\infty} \ge \frac1n. \qedhere
    \end{align}
\end{proof}

This means that weighting with $w$ can improve the error bound by, at most, a factor $C_w = \frac{1}{n}$.
Note that this argument applies to any weight vector $w$, such as those corresponding to optimal subsampling.
We are now ready to prove the claim of Theorem~\ref{thm:weighted_sample_size}.
\end{revisione}

\begin{revisione}[0]
\begin{proof}[Proof of Theorem~\ref{thm:weighted_sample_size}]
    Let $w$ be the weight vector corresponding to the subsample $x_1,\ldots,x_n$.
    Then Theorems~\ref{thm:weighted_error_bound} yields the error bound
    $$
    	\norm{u - u_{\mcal{V},w}}_{L^2(\rho)}^2
    	\le \norm{u - u_{\mcal{V}}}_{L^2(\rho)}^2 + C_w\tfrac{1}{\lambda_{\mathrm{min}}(G_{\sfrac{\boldsymbol{1}}{n}})}\norm{u - u_{\mcal{V}}}_{L^{\infty}(\rho)}^2 .
    $$
    Since $C_w\ge\frac1n>0$ by Lemma~\ref{lem:Cw_bound}, we require $\lambda_{\mathrm{min}}(G_{\sfrac{\boldsymbol{1}}{n}}) = \tfrac{C_w}{c}$ to ensure the desired bound.
    The first statement thus follows from Theorem~\ref{thm:dlogd} and the second statement follows from Theorem~\ref{thm:sample_size_bounds} with the fact that $\gamma \in (\tfrac{2}{\delta^2}, \tfrac{3}{\delta^2})$ for all $\delta\in(0,1)$.
\end{proof}
\end{revisione}

\begin{figure}
    \centering
    \begin{subfigure}[b]{0.24\textwidth}
        \centering
        \includegraphics[width=\textwidth]{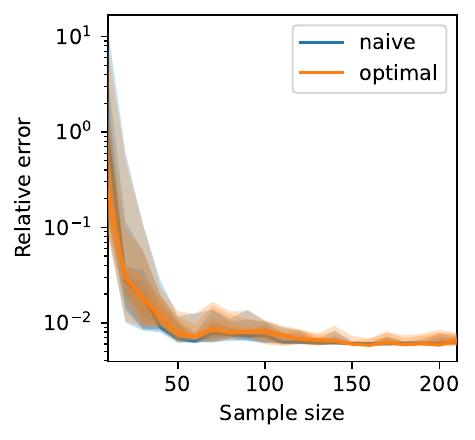}
        \caption{$u(x) := \sin(2\pi x)$}
        \label{fig:subsampling:sin}
    \end{subfigure}
    \hfill
    \begin{subfigure}[b]{0.24\textwidth}
        \centering
        \includegraphics[width=\textwidth]{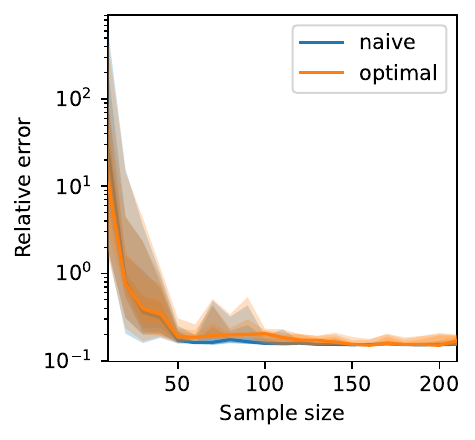}
        \caption{$u(x) := \frac{1}{1 + 25x^2}$}
        \label{fig:subsampling:runge}
    \end{subfigure}
    \hfill
    \begin{subfigure}[b]{0.24\textwidth}
        \centering
        \includegraphics[width=\textwidth]{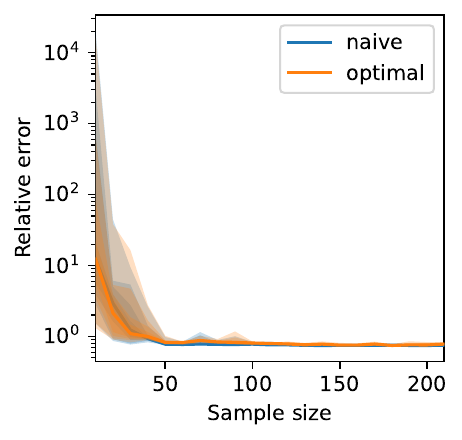}
        \caption{$u(x) := \min\braces{x^{-2}, 10^3}$}
        \label{fig:subsampling:peak}
    \end{subfigure}
    \hfill
    \begin{subfigure}[b]{0.24\textwidth}
        \centering
        \includegraphics[width=\textwidth]{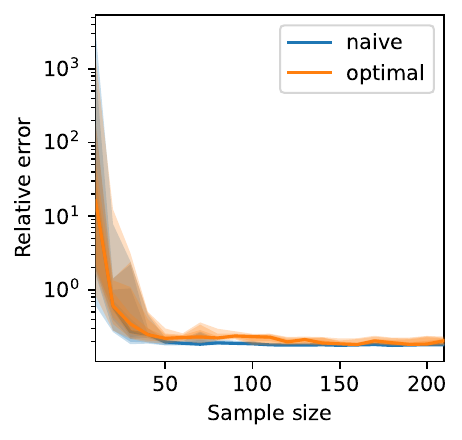}
        \caption{$u(x) := \chi_{[0, 1]}(x)$}
        \label{fig:subsampling:step}
    \end{subfigure}
    \caption{For every sample size $n$, sample points are drawn at random and used either directly (naive) or in conjunction with an optimal weighting (optimal) for a least squares regression.
    The least squares regression is performed for $\rho$ being the uniform measure on $[-1, 1]$ and with a basis of Legendre polynomials of dimension $10$.
    Four different functions are tested, and each experiment is repeated $10$ times.
    The plots show different quantiles of the relative error of the resulting approximation.
    To compute the optimal weight vector $w\in[0,\infty)^n$, we take the vector of Legendre polynomials $b:[-1,1]\to\mbb{R}^{10}$ and define the weighted empirical Gramian $G_w := \sum_{i=1}^n w_i b(x_i) b(x_i)^\intercal$.
    Then $w$ is optimised as to maximise $\lambda_{\mathrm{min}}(G_w)$ subject to $\lambda_{\mathrm{max}}(G_w) \le 2$.
    }
    \label{fig:subsampling}
\end{figure}

\revision[0]{Theorem~\ref{thm:weighted_sample_size} implies that a subsampled version of the original least squares problem can not be stable (in the sense that the approximation error is quasi-optimal) if the original problem is not stable.}
This also makes intuitive sense since we would expect the overall error of an approximation problem to increase when equations are removed from the (overdetermined) least squares system.
A numerical illustration of this fact is provided in Figure~\ref{fig:subsampling}.

\section{Proof of Lemma~\ref{lem:framing_equivalence}}
\label{app:matrix_framing}

{
\renewcommand{\thelemma}{\ref*{lem:framing_equivalence}}
\begin{lemma}
    Let $G, H \in\mbb{R}^{D\times D}$ \revision[0]{be symmetric and positive semi-definite}.
    Then
    $$
        C^{-1} G \preceq H \preceq c^{-1} G
        \qquad\Leftrightarrow\qquad
        \ker(H) = \ker(G)
        \quad\text{and}\quad
        c \mfrak{K}_G \le \mfrak{K}_{H} \le C\mfrak{K}_G .
    $$
\end{lemma}
\addtocounter{lemma}{-1}
}
\begin{proof}
    The proof of this lemma relies on the equivalence
    $$
        C^{-1} G \preceq H \preceq c^{-1} G
        \quad\Leftrightarrow\quad
        CG^+ \succeq H^+ \succeq cG^+ .
    $$
    To show that
    $$
        C^{-1} G \preceq H \preceq c^{-1} G
        \qquad\Rightarrow\qquad
        \ker(H) = \ker(G)
        \quad\text{and}\quad
        c \mfrak{K}_G \le \mfrak{K}_{H} \le C\mfrak{K}_G ,
    $$
    observe that the assertion $\ker(G) = \ker(H)$ is directly implied by the matrix framing $C^{-1}G\preceq H\preceq c^{-1}G$.
    The second assertion holds because
    $$
        c\mfrak{K}_G(x)
        = \inner{cG^+, B(x)B(x)^\intercal}_{\mathrm{Fro}} 
        \le \underbrace{\inner{H^+, B(x)B(x)^\intercal}_{\mathrm{Fro}} }_{= \mfrak{K}_H(x)}
        \le \inner{CG^+, B(x)B(x)^\intercal}_{\mathrm{Fro}} 
        = C\mfrak{K}_G(x) .
    $$
    To prove that
    $$
        \ker(H) = \ker(G)
        \quad\text{and}\quad
        c \mfrak{K}_G \le \mfrak{K}_{H} \le C\mfrak{K}_G
        \qquad\Rightarrow\qquad
        C^{-1} G \preceq H \preceq c^{-1} G,
    $$
    suppose that $\ker(H) = \ker(G)$ and observe that
    $$
        c\mfrak{K}_G \le \mfrak{K}_{H} \le C\mfrak{K}_G
        \quad\Rightarrow\quad
        \inner{cG^+ - H^+, B(x)B(x)^\intercal}_{\mathrm{Fro}} \le 0 \le \inner{CG^+ - H^+, B(x)B(x)^\intercal}_{\mathrm{Fro}} .
    $$
    This implies
    $$
        0 \le v^\intercal(CG^+ - H^+)v
        \qquad\text{and}\qquad
        0 \le v^\intercal(H^+ - cG^+)v
    $$
    \revision[0]{for all $v \in\operatorname{span}(B(\mathcal{X})) = \operatorname{span}\{B(x) : x\in\mathcal{X}\}$.
    Now let $v\in\operatorname{span}(B(\mathcal{X}))^\perp$.
    Then $B(x)^\intercal v = 0$ for all $x\in\mathcal{X}$.
    This implies $B^\intercal v = 0$ and, consequently, $v^\intercal G v = \|B^\intercal v\|_{L^2(\rho)}^2 = 0$.
    Since $v^\intercal G v = \|G^{\sfrac12}v\|_2^2$, it follows that $Gv = G^{\sfrac12} (G^{\sfrac12} v) = 0$.
    This implies $v \in \ker(G)$.
    }
    Since $\ker(G^+) = \ker(G)$ and $\ker(H^+) = \ker(H) = \ker(G)$ by assumption, we conclude that
    $$
        0 \le v^\intercal(CG^+ - H^+)v
        \qquad\text{and}\qquad
        0 \le v^\intercal(H^+ - cG^+)v
    $$
    for all $v \in\mbb{R}^D$ or, equivalently, $CG^+ \succeq H^+$ and $H^+ \succeq cG^+$.
\end{proof}

Lemma~\ref{lem:framing_equivalence} states that the matrix framing
$
    C^{-1}G \preceq \hat H \preceq c^{-1}G
$
is sufficient but not necessary to obtain the uniform framing (cf.~\eqref{eq:framing})
$$
    c \mfrak{K}_G \le \mfrak{K}_{H} \le C\mfrak{K}_G .
$$
This can be exemplified by the polynomial dictionary $B(x) := (1, x, x^2)^\intercal$ on the compact interval $\mcal{X} = [-r, r] \subseteq \mathbb{R}$ with $r \ge 1$ and the choices $G=\operatorname{diag}(1, 1, 0)$ and $H=\operatorname{diag}(1, 0, 1)$.
Since $\ker(G) \ne \ker(H)$, it is easy to see that the matrix framing is not satisfied.
However, since the domain $\mcal{X}$ is compact and the dictionary $B$ is continuous, the uniform framing is always satisfied.

\revision[0]{Indeed, since $\mathfrak{K}_G(x) = 1 + x^2$ and $\mathfrak{K}_H(x) = 1 + x^4$, it is easy to see that $c = 2\sqrt{2} - 2$ and $C = \frac{1 + r^4}{1+r^2} \ge r^2$.
This also demonstrates that, although the ratio $\frac{C}{c}$ is indeed finite, it may be arbitrarily large, depending on the domain $\mcal{X}$. (It is infinite for non-compact domains.)}

\section{Proof of Lemmas~\ref{lem:G_mean},~\ref{lem:G_cov},~\ref{lem:var_G_update}  and~\ref{lem:var_G}}

\subsection{Proof of Lemma~\ref{lem:G_mean}}
\label{app:proof:lem:G_mean}
{
\renewcommand{\thelemma}{\ref*{lem:G_mean}}
\begin{lemma}
    It holds that $\mathbb{E}\big[\hat{G}^{(k)}\big] = \mathbb{E}\big[\hat{G}^{(k+\sfrac12)}\big] = G$ for all $k\ge 0$.
\end{lemma}
}
\begin{proof}
    Observe that
    \begin{align}
    	\mathbb{E}\bracs*{\hat{G}^{(k+\sfrac{1}{2})}\,\Big|\, \hat{G}^{(k)}}
    	&= \frac{1}{n}\sum_{i=1}^n \mathbb{E}\bracs*{\hat{w}_k\pars{y_i} B\pars{y_i}B\pars{y_i}^\intercal\,\Big|\, \hat{G}^{(k)}}
    	\quad\text{where}\quad
    	y_i\overset{\mathrm{i.i.d.}}{\sim} \mu^{(k)}\\
    	&= \mathbb{E}\bracs*{\int \hat{w}_k BB^\intercal \dx[\mu^{(k)}] \,\Bigg|\, \hat{G}^{(k)}} \\
    	&= \mathbb{E}\bracs*{\int \frac{\hat{z}_{k}}{\mfrak{K}_{\hat{G}^{(k)}}} BB^\intercal \frac{\mfrak{K}_{\hat{G}^{(k)}}}{z_{\hat{G}^{(k)}}}\dx[\rho] \,\Bigg|\, \hat{G}^{(k)}} \\
    	&= \mathbb{E}\bracs*{\frac{\hat{z}_{k}}{z_{\hat{G}^{(k)}}} \,\Bigg|\, \hat{G}^{(k)}} \int BB^\intercal \dx[\rho] \\
    	&= \frac{\mathbb{E}\bracs*{\hat{z}_{k}\,\big|\, \hat{G}^{(k)}}}{z_{\hat{G}^{(k)}} } G \\
    	&= G .
    \end{align}
    Thus, by the tower property,
    $
    	\mathbb{E}\big[\hat{G}^{(k+\sfrac{1}{2})}\big]
    	= \mathbb{E}\big[
            \mathbb{E}\big[\hat{G}^{(k+\sfrac{1}{2})}\,\big|\, \hat{G}^{(k)}\big]
        \big]
    	= G
    $.
    As a direct consequence, we obtain
    \begin{align}
        \mathbb{E}\bracs*{\hat{G}^{(k)}}
        &= \frac{1}{k}\sum_{j=0}^{k-1}\mathbb{E}\bracs*{\hat{G}^{(j+\sfrac{1}{2})}}
        = G .
        \qedhere
    \end{align}
\end{proof}

\subsection{Proof of Lemma~\ref{lem:G_cov}}
\label{app:proof:lem:G_cov}
{
\renewcommand{\thelemma}{\ref*{lem:G_cov}}
\begin{lemma}
    Let $j\ne k\in\mathbb{N}$ and $\alpha,\beta\in \{1, \ldots, D\}^2$ be arbitrary multi-indices.
    Then $\operatorname{Cov}\big(\hat{G}^{(k+\sfrac12)}_{\alpha}, \hat{G}^{(j+\sfrac12)}_{\beta}\big) = 0$.
\end{lemma}
}
\begin{proof}
    Without loss of generality, we may assume that $k > j$.
    Then
    \begin{align}
    	\operatorname{Cov}\pars*{\hat{G}^{(k+\sfrac12)}_\alpha, \hat{G}^{(j+\sfrac12)}_\beta}
    	&= \mathbb{E}\bracs*{\pars*{\hat{G}^{(k+\sfrac12)}_\alpha - G_\alpha} \pars*{\hat{G}^{(j+\sfrac12)}_\beta - G_\beta}} \\
    	&= \mathbb{E}\bracs*{\hat{G}^{(k+\sfrac12)}_\alpha \hat{G}^{(j+\sfrac12)}_\beta}
    	- G_\alpha G_\beta \\
    	&= \mathbb{E}\bracs*{\mathbb{E}\bracs*{\hat{G}^{(k+\sfrac12)}_\alpha \hat{G}^{(j+\sfrac12)}_\beta\,\Big\vert\, \hat{G}^{(\sfrac12)}, \ldots,\hat{G}^{(k-1+\sfrac12)}}}
    	- G_\alpha G_\beta \\
    	&= \mathbb{E}\bracs*{\mathbb{E}\bracs*{\hat{G}^{(k+\sfrac12)}_\alpha \,\Big\vert\, \hat{G}^{(\sfrac12)}, \ldots,\hat{G}^{(k-1+\sfrac12)}} \hat{G}^{(j+\sfrac12)}_\beta}
    	- G_\alpha G_\beta \\
    	&= G_\alpha \mathbb{E}\bracs*{\hat{G}^{(j+\sfrac12)}_\beta}
    	- G_\alpha G_\beta \\
    	&= 0 . \qedhere
    \end{align}
\end{proof}

\subsection{Proof of Lemma~\ref{lem:var_G_update}}
\label{app:proof:lem:var_G_update}
{
\renewcommand{\thelemma}{\ref*{lem:var_G_update}}
\begin{lemma}
    Let $k\in\mathbb{N}$ and $\alpha\in \{1, \ldots, D\}^2$ be an arbitrary multi-index.
    % Then $\operatorname{Var}\big(\hat{G}^{(k+1)}_{\alpha}\big) = \tfrac{k^2}{(k+1)^2} \operatorname{Var}\big(\hat{G}^{\pars{k}}_\alpha\big)
    % + \tfrac{1}{(k+1)^2} \operatorname{Var}\big(\hat{G}^{\pars{k+\frac12}}_\alpha\big)
    % $.
    Then
    $$
        \operatorname{Var}\big(\hat{G}^{(k+1)}_{\alpha}\big)
        = \frac{1}{(k+1)^2} \sum_{j=-1}^k \operatorname{Var}\big(\hat{G}^{\pars{j+\sfrac12}}_\alpha\big) ,
    $$
    \revision[0]{with the convention $\hat{G}^{\pars{-\sfrac12}} = \hat G^{(0)}$.}
\end{lemma}
}
\begin{proof}
    Using Lemma~\ref{lem:G_cov}, it follows that
    \begin{align}
    	\operatorname{Var}\big(G_{\alpha}^{(k+1)}\big)
    	&= \operatorname{Cov}\big(\hat{G}^{\pars{k+1}}_\alpha, \hat{G}^{(k+1)}_\alpha\big) \\
    	&= \tfrac{k^2}{(k+1)^2} \operatorname{Cov}\big(\hat{G}^{\pars{k}}_\alpha, \hat{G}^{(k)}_\alpha\big)
        + \tfrac{2k}{(k+1)^2} \operatorname{Cov}\big(\hat{G}^{\pars{k}}_\alpha, \hat{G}^{(k+\sfrac12)}_\alpha\big)
        + \tfrac{1}{(k+1)^2} \operatorname{Cov}\big(\hat{G}^{\pars{k+\sfrac12}}_\alpha, \hat{G}^{(k+\sfrac12)}_\alpha\big) \\
    	&= \tfrac{k^2}{(k+1)^2} \operatorname{Var}\big(\hat{G}^{\pars{k}}_\alpha\big)
        + \tfrac{1}{(k+1)^2} \operatorname{Var}\big(\hat{G}^{\pars{k+\sfrac12}}_\alpha\big) .
    \end{align}
    The claim now follows by induction.
\end{proof}

\subsection{Proof of Lemma~\ref{lem:var_G}}

The proof of Lemma~\ref{lem:var_G} relies on the subsequent three lemmas.

\begin{lemma}
\label{lem:BlBm_le_KG}
\begin{revisione}[0]
    Let $H$ be symmetric and positive semi-definite with $\ker(H) \subseteq \ker(G)$. Then $B_l B_m \le \mathfrak{K}_{H} H_{lm}$.
\end{revisione}
\end{lemma}
\begin{proof}
\begin{revisione}[0]
    Proposition~\ref{prop:inverse_christoffel_definition} and the symmetry of $H$ imply that
    $$
        \mathfrak{K}_{H}(x)
        = \sup_{v\in\mathbb{R}^D} \frac{\abs{B(x)^\intercal v}^2}{v^\intercal H v}
        = \sup_{v_1, v_2\in\mathbb{R}^D} \frac{v_1^\intercal B(x)B(x)^\intercal v_2}{v_1^\intercal H v_2} .
    $$
    This means that
    \begin{align}
        B_l B_m
        &= \frac{B_l B_m}{e_l^\intercal H e_m} H_{lm}
        \le \mathfrak{K}_{H} H_{lm} .
        \qedhere
    \end{align}
\end{revisione}
\end{proof}

\begin{lemma}
\label{lem:inner_bounds}
    Let $X\in\mathbb{R}^{D\times D}$ be arbitrary and $Y\in\mathbb{R}^{D\times D}$ be symmetric and positive semi-definite.
    % Let $\lambda_{\mathrm{max}}\pars{X}$ denote the largest eigenvalue of $X$.
    Then
    $$
    	\inner{X, Y}_{\mathrm{Fro}}
    	\le \norm{X}_2\operatorname{tr}\pars{Y}.
    $$
    Suppose, moreover, that $X$ is symmetric and positive semi-definite and $\ker(X) \subseteq \ker(Y)$.
    Let $\lambda_{\mathrm{min}>0}\pars{X}$ denote the smallest positive eigenvalue of $X$, with the convention that $\lambda_{\mathrm{min}>0}(0) = 0$.
    Then
    $$
    	\lambda_{\mathrm{min}>0}\pars{X} \operatorname{tr}\pars{Y}
    	\le \inner{X, Y}_{\mathrm{Fro}} .
    $$
\end{lemma}
\begin{proof}
    \revision[0]{While this is a standard result, we provide a proof for the sake of completeness.}
    Let $Y := U\Lambda U^\intercal$ be the full spectral decomposition of $Y$ \revision[0]{with $U,\Lambda\in \mathbb{R}^{D\times D}$} and observe that
    \begin{align}
    	\inner{X,Y}_{\mathrm{Fro}}
    	= \operatorname{tr}\pars{X^\intercal Y}
    	= \operatorname{tr}\pars{\underbrace{\pars{U^\intercal X^\intercal U}}_{=:\tilde{X}}\Lambda}
    	= \sum_{j=1}^D \Lambda_{jj} e_j^\intercal \tilde{X} e_j .
    \end{align}
    Then, since $\Lambda_{jj} \ge 0$ and $e_j^\intercal \tilde{X} e_j \le \norm{\tilde{X}}_2 = \norm{X}_2$ for all $j=1,\ldots, D$, it holds that
    \begin{align}
    	\inner{X,Y}_{\mathrm{Fro}} 
    	&= \sum_{j=1}^D \Lambda_{jj} e_j^\intercal \tilde{X} e_j
    	\le \sum_{j=1}^D \Lambda_{jj} \norm{\tilde{X}}_2
        = \norm{X}_2 \operatorname{tr}(Y) .
    % \intertext{Moreover, if $X$ is positive semi-definite, then $e_j^\intercal\tilde{X}e_j \ge 0$ and}
    % 	\inner{X,Y}_{\mathrm{Fro}} 
    % 	&= \sum_{j=1}^D \Lambda_{jj} e_j^\intercal \tilde{X} e_j
    % 	\ge \sum_{j=1}^D \lambda_{\mathrm{min}>0}\pars{\Lambda} e_j^\intercal \tilde{X} e_j
    %     = \lambda_{\mathrm{min}>0}(Y) \operatorname{tr}(X) .
    \end{align}
    To prove the second claim, let $X := U\Lambda U^\intercal$ be the spectral decomposition of $X$, where the eigenvalues in $\Lambda$ are arranged in decreasing order and let $r := \operatorname{rk}(X)$.
    Then
    \begin{align}
    	\inner{X,Y}_{\mathrm{Fro}}
    	= \operatorname{tr}\pars{X^\intercal Y}
    	= \operatorname{tr}\pars{\Lambda \underbrace{\pars{U^\intercal Y^\intercal U}}_{=:\tilde{Y}}}
    	= \sum_{j=1}^r \Lambda_{jj} e_j^\intercal \tilde{Y} e_j .
    \end{align}
    The assumption $\ker\pars{X} \subseteq \ker\pars{Y}$ then implies $e_j^\intercal \tilde{Y} e_j = \pars{Ue_j}^\intercal Y \pars{Ue_j} = 0$ for all $j = r+1,\ldots, D$.
    Moreover, since $Y$ is positive semi-definite, it holds that $e_j^\intercal\tilde{Y}e_j \ge 0$ for all $j=1,\ldots,r$ and, consequently,
    \begin{align}
    	\inner{X,Y}_{\mathrm{Fro}} 
    	&= \sum_{j=1}^r \Lambda_{jj} e_j^\intercal \tilde{Y} e_j
    	\ge \sum_{j=1}^r \lambda_{\mathrm{min}>0}\pars{\Lambda} e_j^\intercal \tilde{Y} e_j
        = \lambda_{\mathrm{min}>0}(X) \operatorname{tr}(Y) .
        \qedhere
    \end{align}
\end{proof}

\begin{lemma}
\label{lem:var_z}
    It holds that
    \begin{align}
        \operatorname{Var}\big(\hat{z}^{(k)} \,\big|\, \hat{G}^{(k)},\hat{m}\big)
        \le \frac{(z^{(k)})^2 V}{\hat{m}}
    \end{align}
    with $V := \mathbb{E}\bracs*{\norm{b(x_1)b(x_1)^\intercal - I}_{2}^2} \le \|\mathfrak{K}_G\|_{L^2(\rho)}^2 - 2(d-1)$.
\end{lemma}
\begin{proof}
    To estimate the conditional variance in the second factor, 
    we compute
    \begin{align}
        \operatorname{Var}\pars*{\hat{z}^{(k)} \,\big|\, \hat{G}^{(k)}, \hat m}
        = \frac{1}{\hat m^2} \sum_{i=1}^{\hat m} \operatorname{Var}\pars*{\mfrak{K}_{\hat{G}^{(k)}}(x_i) \,\Big|\, \hat{G}^{(k)}, \hat m}
        = \frac{1}{\hat m} \operatorname{Var}\pars*{\mfrak{K}_{\hat{G}^{(k)}}(x_1) \,\Big|\, \hat{G}^{(k)}} .
    \end{align}
    Now let $G = U\Lambda U^\intercal$ be a rank-revealing spectral decomposition of $G$ where $U\in\mbb{R}^{D\times d}$ and $\Lambda\in\mbb{R}^{d\times d}$ with $\operatorname{rank}(\Lambda) = d$.
    Given this decomposition, we define the $L^2(\rho)$-orthogonal basis $b = \Lambda^{-1/2} U^\intercal B$ (but we could use any arbitrary $L^2(\rho)$-orthogonal basis of $\mcal{V}$) and 
    utilise Lemma~\ref{lem:inner_bounds} to bound
    \begin{align}
        \operatorname{Var}\bracs*{\mfrak{K}_{\hat{G}^{(k)}}(x_1) \,\Big|\, \hat{G}^{(k)}}
        &= \mathbb{E}\bracs*{\big\langle(\hat{G}^{(k)})^{+},\ B(x_1)B(x_1)^\intercal - G\big\rangle_{\mathrm{Fro}}^2 \,\Big|\, \hat{G}^{(k)}} \\
        &= \mathbb{E}\bracs*{\big\langle(\hat{G}^{(k)})^{+},\ U\Lambda^{\sfrac12} (b(x_1)b(x_1)^\intercal - I) \Lambda^{\sfrac12} U^\intercal\big\rangle_{\mathrm{Fro}}^2 \,\Big|\, \hat{G}^{(k)}} \\
        &\le \operatorname{tr}\big( \Lambda^{\sfrac12}U^\intercal (\hat{G}^{(k)})^{+} U\Lambda^{\sfrac12} \big)^2 \mathbb{E}\bracs*{\norm{b(x_1)b(x_1)^\intercal - I}_{2}^2} \\
        &= \big\langle (\hat{G}^{(k)})^+, G\big\rangle^2 \mathbb{E}\bracs*{\norm{b(x_1)b(x_1)^\intercal - I}_{2}^2} \\
        % &= z_{\hat{G}^{(k)}}^2 \mathbb{E}\bracs*{\norm{b(x_1)b(x_1)^\intercal - I}_{2}^2} 
        &= z_{\hat{G}^{(k)}}^2 V .
    \label{eq:Var_K_le_V}
    \end{align}
    Here, the last equality follows from the definition of $z_{\hat{G}^{(k)}}$ and with $V := \mathbb{E}\big[\|b(x_1)b(x_1)^\intercal - I\|_{2}^2\big]$.
    Moreover, it holds that
    % The finiteness of $V$ can be concluded by the same arguments as before, namely
    % \begin{align}
    %     V
    %     &= \mathbb{E}\bracs*{\norm{b(x_1)b(x_1)^\intercal - I}_{2}^2}
    %     \le \mathbb{E}\bracs*{\norm{b(x_1)b(x_1)^\intercal - I}_{\mathrm{Fro}}^2}
    %     = \mathbb{E}\bracs*{\norm{b(x_1)b(x_1)^\intercal}_{\mathrm{Fro}}^2} - \norm{I}_{\mathrm{Fro}}^2 \\
    %     &= \sum_{k=1}^d \norm{b_k}_{L^4(\rho)}^4 - d
    %     \le d(C_{\mcal{V}}^4 - 1) ,
    % \end{align}
    % where $C_{\mcal{V}} := \sup_{v\in\mcal{V}} \frac{\norm{v}_{L^4(\rho)}}{\norm{v}_{L^2(\rho)}}$.
    $$
        \|b(x_1)b(x_1)^\intercal - I\|_{2}^2
        = \sup_{\|v\|_2 = 1} ((b(x_1)^\intercal v)^2 - 1)^2
        = \max\{1, (\|b(x_1)\|_2^2 - 1)^2\}
        \le 1 + (\mathfrak{K}_G(x_1) - 1)^2 .
        % = \max\{1, \mathfrak{K}_G(x_1)^2 - 2\mathfrak{K}_G(x_1)^2 + 1\}
        % = 1 + \max\{0, \mathfrak{K}_G(x_1)^2 - 2\mathfrak{K}_G(x_1)\} .
    $$
    Taking the expectation yields the bound
    \begin{align}
        V
        &= \mathbb{E}\bracs*{\norm{b(x_1)b(x_1)^\intercal - I}_{2}^2}
        \le 1 + \mathbb{E}[\mathfrak{K}_G(x_1)^2] - 2d + 1
        = \|\mathfrak{K}_G\|_{L^2(\rho)}^2 - 2(d-1) . \label{eq:V_bound}
    \end{align}
    This proves the claim.
    % KKT condition
    % \begin{align}
    %     &\partial_v ((a^\intercal v)^2 - 1)^2 + \mu\partial_v(v^\intercal v - 1)
    %     = 2 ((a^\intercal v)^2 - 1) \partial_v ((a^\intercal v)^2 - 1) + 2\mu v
    %     = 4 ((a^\intercal v)^2 - 1) (a^\intercal v) a + 2\mu v
    %     \overset!= 0 \\
    %     \Leftrightarrow\quad & (a^\intercal v)^3 a = (a^\intercal v) a + \tfrac\mu2 v \\
    %     \Leftrightarrow\quad & (a\perp v \text{ and }\mu=0) \quad\text{or}\quad (a\parallel v \text{ and } \alpha^2 = 1 + \tfrac{\mu}{2}) ,  % (v = \alpha a / \|a\|_2^2)
    % \end{align}
    % where the last equivalence follows from the choice $v = \tfrac{\alpha}{\|a\|_2^2} a$.
    % Note that the condition $\|v\|_2 = 1$ implies $\alpha = \|a\|_2$.
    % Therefore,
    % $$
    %     (a\perp v \text{ and }\mu=0) \quad\text{or}\quad (a\parallel v \text{ and } \mu = 2(\|a\|_2^2 - 1)) ,
    % $$
    % and, consequently,
    % $$
    %     \sup_{\|v\|_2 = 1} ((b(x_1)^\intercal v)^2 - 1)^2
    %     = \max\{1, (\|b(x_1)\|_2^2 - 1)^2\}
    %     = \max\{1, (\mathfrak{K}_G(x_1) - 1)^2\}
    % $$
\end{proof}

\label{app:proof:lem:var_G}
{
\renewcommand{\thelemma}{\ref*{lem:var_G}}
\begin{lemma}
\begin{revisione}[0]
    If $\hat{n}^{(k+1)} \ge N \hat{z}^{(k)}$, it holds that
    $$
    	\operatorname{Var}\big(\hat{G}^{(k+\sfrac12)}_\alpha\big)
    	\le \frac{1}{N} \kappa_{>0}(G) G_\alpha^2 ,
    $$
    where $\kappa_{>0}(G) := \frac{\lambda_{\mathrm{max}}(G)}{\lambda_{\mathrm{min}>0}(G)}$.
    Moreover, if $\hat{n}^{(k+1)}$ is independent of $\hat{z}^{(k)}$, given $\hat{G}^{(k)}$ and $\hat{m}$, it holds that
    $$
    	\operatorname{Var}\big(\hat{G}^{(k+\sfrac12)}_\alpha\,\big|\,\hat{n}^{(k+1)}\big)
    	\le \frac{1}{\hat n^{(k+1)}} \pars*{\mathbb{E}\bracs*{\pars*{\frac{V}{\hat{m}} + 1} z^{(k)} \hat{G}^{(k)}_{\alpha} \,\bigg|\, \hat{n}^{(k+1)}} G_{\alpha} - G_\alpha^2} .
    $$
    In particular, if $\hat{m} \ge V$ and $\hat{n}^{(k+1)} \ge \frac{2N}{\lambda_{\mathrm{min}>0}(\hat{G}^{(k)})}$, it holds that
    $$
    	\operatorname{Var}\big(\hat{G}^{(k+\sfrac12)}_\alpha\big)
    	\le \frac{1}{N}\operatorname{tr}(G) G_{\alpha}^2 .
    $$
\end{revisione}
\end{lemma}
}
\begin{proof}
\begin{revisione}[0]
    % To simplify notation, let $\hat{n} := \hat n^{(k+1)}$ and $\hat w_k := \hat w^{(k)}$.
    % {\color{Mahogany}
    % Using Lemmas~\ref{lem:BlBm_le_KG} and~\ref{lem:framing_equivalence}, we first bound
    % $$
    %     \frac{1}{\mathfrak{K}^{(k)}}B_l B_k
    %     \le \frac{\mathfrak{K}}{\mathfrak{K}^{(k)}} G_{kl}
    %     \le c_k^{-1} G_{kl}
    % $$
    % % \begin{align}
    % % 	\operatorname{Var}\big(\hat{G}^{(k+\sfrac12)}_\alpha \,\big|\, \hat{G}^{(k)},\hat{n}^{(k+1)},\hat{m}\big)
    % % 	&= \frac{1}{\hat n^{(k+1)}} \pars*{\mathbb{E}\bracs*{\int (\hat{w}^{(k)})^2 B_l^2 B_m^2 \dx[\mu^{(k)}] \,\bigg|\, \hat{G}^{(k)},\hat{n}^{(k+1)},\hat{m}} - G_\alpha^2} \\
    % % 	&= \frac{1}{\hat n^{(k+1)}} \pars*{\frac{(\hat{z}^{(k)})^2}{z^{(k)}} \mathbb{E}\bracs*{\int \frac{1}{\mathfrak{K}^{(k)}} B_l^2 B_m^2 \dx[\rho] \,\bigg|\, \hat{G}^{(k)},\hat{n}^{(k+1)},\hat{m}} - G_\alpha^2} \\
    % % 	&\le \frac{1}{\hat n^{(k+1)}} \pars*{\frac{(\hat{z}^{(k)})^2}{z^{(k)}} \hat{G}^{(k)}_{\alpha} \mathbb{E}\bracs*{\int B_l B_m \dx[\rho] \,\bigg|\, \hat{G}^{(k)},\hat{n}^{(k+1)},\hat{m}} - G_\alpha^2} \\
    % % 	&= \frac{1}{\hat n^{(k+1)}} \pars*{\frac{(\hat{z}^{(k)})^2}{z^{(k)}} \hat{G}^{(k)}_{\alpha}G_\alpha - G_\alpha^2} .
    % % \end{align}
    % }
    Using Lemma~\ref{lem:BlBm_le_KG}, with $\alpha = (l,m)$, we can bound
    \begin{align}
    	\operatorname{Var}\big(\hat{G}^{(k+\sfrac12)}_\alpha \,\big|\, \hat{G}^{(k)},\hat{n}^{(k+1)},\hat{m}\big)
    	&= \frac{1}{\hat n^{(k+1)}} \pars*{\mathbb{E}\bracs*{\int (\hat{w}^{(k)})^2 B_l^2 B_m^2 \dx[\mu^{(k)}] \,\bigg|\, \hat{G}^{(k)},\hat{n}^{(k+1)},\hat{m}} - G_\alpha^2} \\
    	&= \frac{1}{\hat n^{(k+1)}} \pars*{\frac{(\hat{z}^{(k)})^2}{z^{(k)}} \mathbb{E}\bracs*{\int \frac{1}{\mathfrak{K}^{(k)}} B_l^2 B_m^2 \dx[\rho] \,\bigg|\, \hat{G}^{(k)},\hat{n}^{(k+1)},\hat{m}} - G_\alpha^2} \\
    	&\le \frac{1}{\hat n^{(k+1)}} \pars*{\frac{(\hat{z}^{(k)})^2}{z^{(k)}} \hat{G}^{(k)}_{\alpha} \mathbb{E}\bracs*{\int B_l B_m \dx[\rho] \,\bigg|\, \hat{G}^{(k)},\hat{n}^{(k+1)},\hat{m}} - G_\alpha^2} \\
    	&= \frac{1}{\hat n^{(k+1)}} \pars*{\frac{(\hat{z}^{(k)})^2}{z^{(k)}} \hat{G}^{(k)}_{\alpha}G_\alpha - G_\alpha^2} .
    \end{align}
    Since $\ker(\hat{G}^{(k)}) = \ker(G)$, Lemma~\ref{lem:inner_bounds} implies
    $$
        \frac{\hat{z}^{(k)}}{z^{(k)}}
        \le \frac{\lambda_{\mathrm{max}}(G) \operatorname{tr}((\hat{G}^{(k)})^+)}{\lambda_{\mathrm{min}>0}(G) \operatorname{tr}((\hat{G}^{(k)})^+)} = \kappa_{>0}(G) ,
    $$
    we can estimate
    $$
    	\operatorname{Var}\big(\hat{G}^{(k+\sfrac12)}_\alpha \,\big|\,\hat{G}^{(k)},\hat{n}^{(k+1)},\hat{m}\big)
    	\le \frac{1}{\hat n^{(k+1)}} \pars*{\frac{(\hat{z}^{(k)})^2}{z^{(k)}} \hat{G}^{(k)}_{\alpha}G_\alpha - G_\alpha^2}
    	% \le \frac{1}{\hat n} \frac{\hat{z}_{k}^2}{z_k} \hat{G}^{(k)}_{\alpha}G_\alpha
        \le \frac{\hat{z}^{(k)}}{\hat n^{(k+1)}} \kappa_{>0}(G) \hat{G}^{(k)}_{\alpha}G_\alpha .
    $$
    Choosing $\hat{n}^{(k+1)} \ge N \hat{z}^{(k)}$ and taking the expectation yields the first bound for $\operatorname{Var}\big(\hat{G}^{(k+\sfrac12)}_\alpha\big)$.

    For the second bound, we take the conditional expectation and use the tower property to obtain
    \begin{align}
    	\operatorname{Var}\big(\hat{G}^{(k+\sfrac12)}_\alpha \,\big|\, \hat{n}^{(k+1)}\big)
    	&\le \frac{1}{\hat n^{(k+1)}} \pars*{\mathbb{E}\bracs*{\frac{(\hat{z}^{(k)})^2}{z^{(k)}} \hat{G}^{(k)}_{\alpha} \,\bigg|\, \hat{n}^{(k+1)}}G_\alpha - G_\alpha^2} \\
    	% &= \frac{1}{\hat n} \pars*{\mathbb{E}\bracs*{\frac{\mathbb{E}[\hat{z}_{k}^2\,|\,\hat{G}^{(k)}, \hat{m}, \hat{n}]}{z_k} \hat{G}^{(k)}_{\alpha} \,\bigg|\, \hat{n}} G_{\alpha} - G_\alpha^2} \\
    	&= \frac{1}{\hat n^{(k+1)}} \pars*{\mathbb{E}\bracs*{\frac{\mathbb{E}[(\hat{z}^{(k)})^2\,|\,\hat{G}^{(k)}, \hat{m}]}{z^{(k)}} \hat{G}^{(k)}_{\alpha} \,\bigg|\, \hat{n}^{(k+1)}} G_{\alpha} - G_\alpha^2} \\
    	&= \frac{1}{\hat n^{(k+1)}} \pars*{\mathbb{E}\bracs*{\frac{\operatorname{Var}(\hat{z}^{(k)}\,|\,\hat{G}^{(k)}, \hat{m}) + (z^{(k)})^2}{z^{(k)}} \hat{G}^{(k)}_{\alpha} \,\bigg|\, \hat{n}^{(k+1)}} G_{\alpha} - G_\alpha^2} ,
    \end{align}
    where we have used that $\hat{z}^{(k)}$ and $\hat{n}^{(k+1)}$ are independent, given $\hat{G}^{(k)}$ and $\hat{m}$.
    Using Lemma~\ref{lem:var_z} thus yields the claimed bound,
    \begin{align}
    	\operatorname{Var}\big(\hat{G}^{(k+\sfrac12)}_\alpha \,\big|\, \hat{n}^{(k+1)}\big)
    	&\le \frac{1}{\hat n^{(k+1)}} \pars*{\mathbb{E}\bracs*{\pars*{\frac{V}{\hat{m}} + 1} z^{(k)} \hat{G}^{(k)}_{\alpha} \,\bigg|\, \hat{n}^{(k+1)}} G_{\alpha} - G_\alpha^2} .
    \end{align}
    Finally, using Lemma~\ref{lem:inner_bounds}, we can bound $z^{(k)} \le \lambda_{\mathrm{min}>0}(\hat{G}^{(k)})^{-1} \operatorname{tr}(G)$ and therefore
    \begin{align}
    	\operatorname{Var}\big(\hat{G}^{(k+\sfrac12)}_\alpha \,\big|\, \hat{n}^{(k+1)}\big)
    	&\le \frac{1}{\hat n^{(k+1)}} \pars*{\mathbb{E}\bracs*{\pars*{\frac{V}{\hat{m}} + 1} \frac{\operatorname{tr}(G)}{\lambda_{\mathrm{min}>0}(\hat{G}^{(k)})} \hat{G}^{(k)}_{\alpha} \,\bigg|\, \hat{n}^{(k+1)}} G_{\alpha} - G_\alpha^2} \\
    	&\le \mathbb{E}\bracs*{\frac{1}{\hat n^{(k+1)}} \pars*{\frac{V}{\hat{m}} + 1} \frac{\operatorname{tr}(G)}{\lambda_{\mathrm{min}>0}(\hat{G}^{(k)})} \hat{G}^{(k)}_{\alpha} \,\bigg|\, \hat{n}^{(k+1)}} G_{\alpha} .
    \end{align}
    Choosing $\hat{m} \ge V$ and $\hat{n}^{(k+1)} \ge \frac{2N}{\lambda_{\mathrm{min}>0}(\hat{G}^{(k)})}$ yields the final bound,
    \begin{align}
    	\operatorname{Var}\big(\hat{G}^{(k+\sfrac12)}_\alpha\big)
    	&\le \mathbb{E}\bracs*{\frac{1}{\hat n^{(k+1)}} \pars*{\frac{V}{\hat{m}} + 1} \frac{\operatorname{tr}(G)}{\lambda_{\mathrm{min}>0}(\hat{G}^{(k)})} \hat{G}^{(k)}_{\alpha}} G_{\alpha}
    	\le \frac{1}{N}\operatorname{tr}(G) G_{\alpha}^2 .
        \qedhere
    \end{align}
\end{revisione}
\end{proof}

\section{Refined convergence analysis}
\label{app:update}

\revision[0]{
\ifdraftmode\color{ForestGreen}\fi
To deepen our intuition for why one   iteration improves the framing constants $c_k$ and $C_k$in the relation   
$$
    C_k^{-1} G \preceq \hat{G}^{(k)} \preceq c_k^{-1} G
$$
from Lemma~\ref{lem:framing_equivalence},  
 recall from Remark~\ref{rmk:choosing_m_and_n} that
$$
    (z^{(k)})^{-1} \le C_k^{-1}
    \quad\text{and}\quad
    c_k^{-1} \le Z^{(k)}
$$
with $Z^{(k)} := \inner{\hat{G}^{(k)}, G^+}_{\mathrm{Fro}}$.
A simple example in which   the proposed update reduces $z^{(k)}$ can be provided for $n=1$ when assuming that $B$ is a basis and $\hat{z}^{(k)} = z^{(k)}$.
In this setting, let $\hat{G} := \hat{G}^{(k)}$ be the current Gramian estimate and define $\hat\Delta := w(x) B(x)B(x)^\intercal := w^{(k)}(x_{1}^{(k+1)}) B(x_{1}^{(k+1)})B(x_{1}^{(k+1)})^\intercal$.
To make the dependence on the basis explicit, we add a corresponding subscript to the matrices $G$, $\hat G$ and $\hat \Delta$.
Similarly, we define $\mathfrak{K}_{H,\varphi} = \inner{H^+,  \varphi\varphi^\intercal}_{\mathrm{Fro}}$.
Using this notation, we can conclude that $\mathfrak{K}_{\hat{G}_\varphi,\varphi} = \mathfrak{K}_{\hat{G}_\psi,\psi}$ for all dictionaries $\varphi$ and $\psi$ and that $\varphi(x)^\intercal G_\varphi \varphi(x) = \psi(x)^\intercal G_\psi \psi(x)$ for all dictionaries $\varphi$ and corresponding $\psi = U^\intercal \varphi$ for an orthogonal matrix $U$.
Finally, let $\alpha := \frac{k}{k+1}$ and $\bar{\alpha} := 1 - \alpha$ and define
\begin{align}
    z^{\mathrm{old}} := \inner{\hat{G}_\varphi^+, G_\varphi}_{\mathrm{Fro}}
    \quad&\text{and}\quad
    z^{\mathrm{new}} := \inner{(\alpha\hat{G}_\varphi + \bar\alpha\hat{\Delta}_\varphi)^+, G_\varphi}_{\mathrm{Fro}}
\intertext{as well as}
    Z^{\mathrm{old}} := \inner{\hat{G}_\varphi, G_\varphi^+}_{\mathrm{Fro}}
    \quad&\text{and}\quad
    Z^{\mathrm{new}} := \inner{\alpha\hat{G}_\varphi + \bar\alpha\hat{\Delta}_\varphi, G_\varphi^+}_{\mathrm{Fro}} .
\end{align}

Given these definitions, we start by analysing the behaviour of $Z^{(k)}$.
For this, let $\varphi$ satisfy $G_\varphi = I$ and observe that
$$
    Z^{\mathrm{new}} = \alpha Z^{\mathrm{old}} + \bar\alpha \operatorname{tr}(\hat\Delta_\varphi) .
$$
Since
$$
    \operatorname{tr}(\hat\Delta_\varphi)
    = \|\hat\Delta_\varphi\|_2
    = w(x)\|\varphi(x)\|_2^2
    = \frac{z^{\mathrm{old}}}{\mathfrak{K}_{\hat{G}_\varphi, \varphi}(x)}\|\varphi(x)\|_2^2
    = z^{\mathrm{old}} ,
$$
we conclude that
$$
    Z^{\mathrm{new}} = (\alpha + \bar\alpha\gamma) Z^{\mathrm{old}}
    \quad\text{with}\quad
    \gamma = \frac{z^{\mathrm{old}}}{Z^{\mathrm{old}}} .
$$
Since $\alpha + \bar\alpha\gamma < 1$ for $\gamma < 1$, this ensures decay of $Z^{(k)}$ as long as $z^{(k)} < Z^{(k)}$.
This is satisfied, e.g., when $\hat{G}^{(k)} \succ G$.

We now turn to the analysis of $z^{(k)}$.
For this, let $\varphi$ be satisfy $\hat{G}_\varphi = I$ and $\hat\Delta_\varphi = \|\hat\Delta_\varphi\|_2\ e_1^{\vphantom{\intercal}}e_1^\intercal$ and observe that
$$
    \alpha\hat G_\varphi + \bar\alpha\hat\Delta_\varphi = \begin{pmatrix}
        \alpha + \bar\alpha z^{\mathrm{old}} & 0 \\
        0 & \alpha I
    \end{pmatrix} .
$$
Therefore, 
\begin{align}
    z^{\mathrm{new}}
    &= \inner{(\alpha\hat{G}_\varphi + \bar\alpha\hat{\Delta}_\varphi)^+, G_\varphi}_{\mathrm{Fro}} \\
    &= \frac1\alpha\bigg(\frac{\alpha}{\alpha + \bar\alpha z^{\mathrm{old}}}G_{\varphi,11} +\sum_{k\ge2} G_{\varphi,kk} \bigg) \\
    &= \frac1\alpha\bigg(\Big(\frac{\alpha}{\alpha + \bar\alpha z^{\mathrm{old}}} - 1\Big)G_{\varphi,11} +\sum_{k\ge1} G_{\varphi,kk} \bigg) \\
    &= \frac1\alpha\bigg(\frac{-\bar\alpha z^{\mathrm{old}}}{\alpha + \bar\alpha z^{\mathrm{old}}}G_{\varphi,11} + z^{\mathrm{old}} \bigg) \\
    &= \frac1\alpha\bigg(1 - \frac{\bar\alpha G_{\varphi,11}}{\alpha + \bar\alpha z^{\mathrm{old}}}\bigg) z^{\mathrm{old}}
    .
\end{align}
Since $\alpha = \frac{k}{k+1}$ and $\bar\alpha = \frac1{k+1}$, the equation simplifies to
$$
    z^{\mathrm{new}}
    = \frac{k+1}{k}\bigg(1 - \frac{G_{\varphi,11}}{k + z^{\mathrm{old}}}\bigg) z^{\mathrm{old}} ,
$$
% The factor is bounded by $1-\varepsilon$, if
% % \begin{align}
% %     \frac{k+1}{k}\bigg(1 - \frac{G_{\varphi,11}}{k + z^{\mathrm{old}}}\bigg)
% %     \le (1 - \varepsilon)
% %     \quad&\Leftrightarrow\quad
% %     1 - \frac{G_{\varphi,11}}{k + z^{\mathrm{old}}}
% %     \le \frac{k}{k+1}(1 - \varepsilon)\\
% %     &\Leftrightarrow\quad
% %     \frac{G_{\varphi,11}}{k + z^{\mathrm{old}}}
% %     \ge \frac{1 + k\varepsilon}{k+1}\\
% %     &\Leftrightarrow\quad
% %     G_{\varphi,11}
% %     \ge \frac{1 + k\varepsilon}{k+1} (k + z^{\mathrm{old}}) .
% % \end{align}
% \begin{align}
%     \frac{k+1}{k}\bigg(1 - \frac{G_{\varphi,11}}{k + z^{\mathrm{old}}}\bigg)
%     \le 1 - \varepsilon
%     \quad&\Leftrightarrow\quad
%     (k+1)(k+z^{\mathrm{old}})\bigg(1 - \frac{G_{\varphi,11}}{k + z^{\mathrm{old}}}\bigg)
%     \le (1 - \varepsilon)k(k+z^{\mathrm{old}}) \\
%     \quad&\Leftrightarrow\quad
%     (k+1)(k+z^{\mathrm{old}}) - (k+1)G_{\varphi,11}
%     \le (1 - \varepsilon)k(k+z^{\mathrm{old}}) \\
%     \quad&\Leftrightarrow\quad
%     k^2
%     + (z^{\mathrm{old}} - G_{\varphi,11} + 1)k
%     + (z^{\mathrm{old}} - G_{\varphi,11})
%     \le (1 - \varepsilon)k^2 + (1 - \varepsilon)z^{\mathrm{old}}k \\
%     \quad&\Leftrightarrow\quad
%     \varepsilon k^2
%     + (z^{\mathrm{old}} - G_{\varphi,11} + 1 - (1 - \varepsilon)z^{\mathrm{old}})k
%     + (z^{\mathrm{old}} - G_{\varphi,11})
%     \le 0 \\
%     % \quad&\Leftrightarrow\quad
%     % \varepsilon k^2
%     % + (\beta + 1 - (1 - \varepsilon)z^{\mathrm{old}})k
%     % + \beta
%     % \le 0 \\
% \end{align}
which ensures a geometric decay if
$$
    \frac{k+1}{k}\bigg(1 - \frac{G_{\varphi,11}}{k + z^{\mathrm{old}}}\bigg)
    \le c < 1 .
$$
A decrease of $z^{(k)}$ is guaranteed whenever
\begin{align}
    \frac{k+1}{k}\bigg(1 - \frac{G_{\varphi,11}}{k + z^{\mathrm{old}}}\bigg)
    < 1
    \quad&\Leftrightarrow\quad
    \frac{k+1}{k}\big(k + z^{\mathrm{old}} - G_{\varphi,11}\big)
    < k + z^{\mathrm{old}} \\
    \quad&\Leftrightarrow\quad
    (k+1)\big(k + z^{\mathrm{old}} - G_{\varphi,11}\big)
    < k^2 + z^{\mathrm{old}}k \\
    \quad&\Leftrightarrow\quad
    k^2 + k + k\big(z^{\mathrm{old}} - G_{\varphi,11}\big) + \big(z^{\mathrm{old}} - G_{\varphi,11}\big)
    < k^2 + z^{\mathrm{old}}k \\
    \quad&\Leftrightarrow\quad
    k\big(1 - G_{\varphi,11}\big) + \big(z^{\mathrm{old}} - G_{\varphi,11}\big)
    < 0 .
\end{align}

To find out when this inequality is satisfied, we compute the expectation of $G_{\varphi,11}$ given $\hat{G}_\varphi$.
Since $\hat\Delta_\varphi = w(x) \varphi(x)\varphi(x)^\intercal$, by definition, and $\hat\Delta_\varphi = z^{\mathrm{old}}\ e_1^{\vphantom{\intercal}}e_1^\intercal$, by the choice of $\varphi$, we can conclude that $e_1 = \sqrt{\frac{w(x)}{z^{\mathrm{old}}}} \varphi(x)$.
Using this identity, we can write
$$
    G_{\varphi,11}
    = e_1^\intercal G_\varphi e_1^{\vphantom{\intercal}}
    = \frac{w(x)}{z^{\mathrm{old}}} \varphi(x)^\intercal G_\varphi \varphi(x)
    = \frac{w(x)}{z^{\mathrm{old}}} \inner{G_\varphi, \varphi(x)\varphi(x)^\intercal}_{\mathrm{Fro}} .
$$
This implies
$$
    \mathbb{E}\big[G_{\varphi,11} \,\big|\, \hat G_\varphi\big]
    = \frac{1}{z^{\mathrm{old}}} \int \inner{G_\varphi, \varphi(y)\varphi(y)^\intercal}_{\mathrm{Fro}} \dx[\rho](y)
    = \frac{1}{z^{\mathrm{old}}} \|G_\varphi\|_{\mathrm{Fro}}^2
    \ge \frac{z^{\mathrm{old}}}{d} ,
$$
where the final bound follows from $\|G_\varphi\|_{\mathrm{Fro}} \ge d^{-\sfrac12}\|G_\varphi\|_{\mathrm{nuc}}$ and
$\|G_\varphi\|_{\mathrm{nuc}} = \operatorname{tr}(G_\varphi) = z^{\mathrm{old}}$.
% Since
% $$
%     \frac{z^{\mathrm{old}}}{d}
%     \ge \frac{1 + k\varepsilon}{k+1} (k + z^{\mathrm{old}})
%     \quad\Leftrightarrow\quad
%     z^{\mathrm{old}} - d\frac{1 + k\varepsilon}{k+1} z^{\mathrm{old}}
%     \ge d\frac{1 + k\varepsilon}{k+1} k
% $$
We thus conclude that
% \begin{align}
%     \mathbb{E}\bigg[\frac{k+1}{k}\bigg(1 - \frac{G_{\varphi,11}}{k + z^{\mathrm{old}}}\bigg) \,\bigg|\, \hat{G}_{\varphi}\bigg]
%     < 1
%     \quad&\Leftrightarrow\quad
%     k\big(1 - \mathbb{E}\big[G_{\varphi,11}\,\big|\, \hat{G}_{\varphi}\big] \big) + \big(z^{\mathrm{old}} - \mathbb{E}\big[G_{\varphi,11}\,\big|\, \hat{G}_{\varphi}\big] \big)
%     < 0 \\
%     \quad&\Leftarrow\quad
%     k\Big(1 - \frac{z^{\mathrm{old}}}{d} \Big) + \Big(z^{\mathrm{old}} - \frac{z^{\mathrm{old}}}{d} \Big)
%     < 0 \\
%     \quad&\Leftrightarrow\quad
%     k\big(d - z^{\mathrm{old}}\big) + \big(d - 1\big) z^{\mathrm{old}}
%     < 0
%     % \quad&\Leftrightarrow\quad
%     % z^{\mathrm{old}} > d \quad\text{and}\quad k > \frac{z^{\mathrm{old}}}{z^{\mathrm{old}} - d}(d - 1) .
% \end{align}
\begin{align}
    \mathbb{E}\bigg[\frac{k+1}{k}\bigg(1 - \frac{G_{\varphi,11}}{k + z^{\mathrm{old}}}\bigg) \,\bigg|\, \hat{G}_{\varphi}\bigg]
    < 1
    \quad&\Leftrightarrow\quad
    \frac{k + z^{\mathrm{old}}}{k+1} < \mathbb{E}\big[G_{\varphi,11} \,\big|\, \hat{G}_{\varphi}\big]
    \\
    \quad&\Leftarrow\quad
    \frac{k + z^{\mathrm{old}}}{k+1} < \frac{z^{\mathrm{old}}}{d}
    \\
    \quad&\Leftrightarrow\quad
    % (d-1)z^{\mathrm{old}} - k z^{\mathrm{old}} + kd < 0
    (d-1)z^{\mathrm{old}} < (z^{\mathrm{old}} - d) k.
\end{align}
This bound already holds for $k \ge 2(d-1)$ and $z^{\mathrm{old}} > 2d$, 
% Assuming that $z^{\mathrm{old}} = cd$ for some $c>1$, we see that this bound already holds for
% $$
%     (d-1)cd < (c-1)d k
% $$
% $z^{\mathrm{old}}\ge 2d$ and $k\ge 2d - 1$,
which is reasonable since $z_G = d$.

The potential geometric decay of $z^{(k)}$ and $Z^{(k)}$ is indeed observed in our numerical experiments.
A proof in a different but related setting is provided in~\cite{herremans2025refinementbasedchristoffelsamplingsquares}.
}
% \input{content/appendix/a_inner_bounds}
% \input{content/mcmc_proofs}

% \cleardoublepage
% \input{content/notes}

\end{document}